# Non-commutative $L^p$ spaces and Grassmann stochastic analysis

FRANCESCO DE VECCHI
Department of Mathematics,
University of Pavia,
Via Adolfo Ferrata 5 27100, Pavia, Italy
*Email:* francescocarlo.devecchi@unipv.it

LUCA FRESTA
Institute for Applied Mathematics &
Hausdorff Center for Mathematics,
University of Bonn,
Endenicher Allee 60 53115, Bonn, Germany
*Email:* fresta@iam.uni-bonn.de

MARIA GORDINA
Department of Mathematics
University of Connecticut
Storrs, CT 06269, *USA*
*Email:* maria.gordina@uconn.edu

MASSIMILIANO GUBINELLI
Mathematical Institute,
University of Oxford,
Woodstock Road OX2 6GG Oxford,
United Kingdom
*Email:* gubinelli@maths.ox.ac.uk

**Abstract**

We introduce a theory of non-commutative $L^p$ spaces suitable for non-commutative probability in a non-tracial setting and use it to develop stochastic analysis of Grassmann-valued processes, including martingale inequalities, stochastic integrals with respect to Grassmann Itô processes, Girsanov's formula and a weak formulation of Grassmann SDEs. We apply this new setting to the construction of several unbounded random variables including a Grassmann analog of the $\Phi_2^4$ Euclidean QFT in a bounded region and weak solution to singular SPDEs in the spirit of the early work of Jona-Lasinio and Mitter on the stochastic quantisation of $\Phi_2^4$.

# Table of contents

# 1 Introduction

We continue to develop Grassmann probability and the associated stochastic analysis initiated in [ABDG22, DFG22] and motivated by applications to stochastic quantization of Euclidean fermionic quantum field theories. The starting point is the observation in [ABDG22] that the analysis of the renormalized Wick powers of Grassmann Euclidean free fields does not fit naturally in the $C^*$-algebraic framework of non-commutative probability. Intuitively, this corresponds to the fact that while linear functionals in Grassmann algebras can be represented as bounded operators on a Hilbert space, this does not seem possible for higher Wick monomials which instead behave as unbounded operators very much like their bosonic counterparts. A similar observation has recently been made also by Chandra–Hairer–Peev [CHP23] where the authors develop a setting for "almost sure" analysis of these unbounded random variables via the notion of *locally $C^*$-algebras*. While this localization is a useful concept, applications require to be able to take averages of these unbounded elements and to estimate or compute their moments. More generally, our goal is to develop Grassmann probability and stochastic analysis analogously to the bosonic counterparts.

In this paper we introduce appropriate notions of Banach spaces of non-commutative random variables whose $p$-th powers are integrable, i.e., the equivalent of the standard $L^p$ spaces associated to a classical probability space $(\Omega, \mathscr{F}, \mathbb{P})$ based on a probability measure $\mathbb{P}$ on a measurable space $(\Omega, \mathscr{F})$. Segal [Seg53] recognized quite early the possibility to introduce a useful generalization of a measure to the *non-commutative* setting of a $C^*$-algebra endowed with a *tracial* state $\tau$ (a *gage*). See also the lucid account by Nelson [Nel74]. Unfortunately, the tracial setting does not allow for the embedding of non-trivial Grassmann random variables, since the tracial nature of the state forces the covariance to be zero due to the contrasting requirements of symmetry (by the tracial condition) and antisymmetry (by their Grassmann nature). Following the early work of Osterwalder–Schrader [OS73], Grassmann Gaussian variables of [ABDG22] are constructed over the Fock state of a pair of canonical anti-commutation relation (CAR) algebras. This state is also not suitable to develop an appropriate $L^p$ theory since it is not faithful for the underlying $C^*$-algebra. We are bound, therefore, to look for appropriate states elsewhere with the basic desiderata that they need to represent the expectation of Grassmann Gaussian random variables (and processes or, more generally, fields) and be faithful. Our proposal is to realize Grassmann Gaussian probability over a non-commutative probability space given by a large class of non-Fock (and non-tracial) quasi-free states of the CAR algebra, and in particular by the Araki–Wyss factors [AW64]. These algebras are obtained in the GNS representation associated with a quasi-free state [Ara70, BR97], so that we can extend quite naturally the Osterwalder–Schrader embedding to them in order to represent arbitrary Grassmann Gaussians. We also notice that they are factors of type III, so that our construction proves to be a further significant application of such factors in mathematical physics [Kad04b, Yng05].

The associated Tomita–Takesaki modular automorphism gives rise, via a standard construction due to Haagerup [Haa79a], to an appropriate notion of non-commutative $L^p(\mathscr{M})$ spaces in this non-tracial setting. Non-tracial $L^p$ spaces have been extensively studied in the literature [Pis], in particular standard martingale inequalities are available, see, e.g., the work of Pisier–Xu [PX97], Junge [Jun02], Junge–Xu [JX03], and hypercontractivity of the Araki–Wyss factors has been proved by Lee–Richard [LR11] in the greater generality of the $q$-Gaussian setting introduced by Hiai [Hia03]. However, despite this extensive literature, the use of these non-tracial $L^p$ spaces in non-commutative probabilistic *applications* (other than proving inequalities) has remained underdeveloped mainly because of their abstract nature. While Haagerup's theory is quite efficient and tries to bring the analysis to the tracial setting (via a crossed-product construction originally due to Connes [Con73]), this happens at the price of the appearance of many complications. Various reformulations of the $L^p$ theory have been proposed, in particular by Connes–Hilsum [Hil81], Araki–Masuda [AM82] and more concretely, with applications to Markovian semigroups in quantum statistical mechanics by Majewski–Zegarliski [MZ96].

The main difficulty is that there is no canonical identification of elements of the von Neumann algebra (vNa) $\mathscr{M}$ in the $L^p$ space and as a consequence there is no canonical way to extend the algebra product into a continuous bilinear map compatible with the expected Hölder inequality. Another difficulty is related to the abstract nature of these spaces, whereas in non-commutative probability it would be appropriate to identify normal unbounded elements affiliated with the vNa $\mathscr{M}$ with unbounded complex-valued random variables (i.e. operators on a Hilbert space). No such identification is possible, to our knowledge, in general for the normal elements of the $L^p$ spaces. Some early ideas on how to associate concrete closed unbounded operators on a Hilbert space to elements of $L^2$ spaces of general quasi-free states of CAR and CCR algebras are present in the pioneering work of Barnett–Streater–Wilde [BSW83a], which follows their systematic development of the $L^2$-based stochastic calculus in the tracial CAR setting [BSW82]. See also the follow-up work of Wilde [Wil88]. In [BSW83a] the authors identify some additional conditions which allow one to associate a closed operator to a subset of elements of $L^2$. However, the work of these authors is restricted to the $L^2$ setting and due to the lack of a natural product structure, it does not allow them to construct a fully fledged stochastic calculus.

Our main contributions can therefore be summarized as follows.

a) We define Banach spaces $\mathbb{L}^p$ (twisted $L^p$ spaces) in such a way that they densely contain the algebra $\mathscr{M}_a$ of analytic elements of $\mathscr{M}$. These spaces are constructed via Haagerup's $L^p$ spaces and for which the product in $\mathscr{M}_a$ extends canonically to a continuous bilinear map $\mathbb{L}^p \times \mathbb{L}^q \to \mathbb{L}^r$ for the standard Hölder relation among exponents $1/p + 1/q = 1/r$ (Section 2.3).

b) We extend standard non-commutative probabilistic results to these spaces, including martingale inequalities and hypercontractivity of Gaussian random variables (in particular Grassmann-valued).

c) Given a filtered probability space supporting a Grassmann Brownian motion, we construct natural notions of $\mathbb{L}^p$-valued Itô stochastic integrals and Itô processes (Section 4.1) and prove an Itô formula for polynomial functionals of Itô processes (Section 4.2).

d) We prove a Grassmann version of Lévy's martingale characterization of Brownian motion, Girsanov's theorem (Section 4.3) and give a notion of weak solutions to SDEs (Section 4.4)

e) We use the above tools to construct various examples of unbounded random variables in the Grassmann context, in particular:

   - A Grassmann analog of the $\Phi_2^4$ measure in a bounded domain (Section 5.1);
   - A Grassmann analog of a Langevin SPDE driven by a renormalized cubic polynomial (Section 5.2);

f) We identify conditions under which elements in $L^p$ can be canonically interpreted as unbounded operators on a Hilbert space $\mathcal{H}$ (here $\mathcal{H}$ can be taken to be $L^2(\mathcal{M})$ without loss of generality) (Appendix B).

Note that, even with these $\mathbb{L}^p$ spaces we are not able to provide a positive answer to the original problem formulated in [ABDG22], namely to solve non-linear singular renormalized Grassmann SPDE with additive noise. The main difficulty is the lack of an appropriate Banach space where we can find estimates of solutions via fixed-point methods or, alternatively, via appropriate coercive estimates. Informally put "the estimates do not close". Yet, we find that our investigation uncovered a very interesting and largely unexplored area of non-commutative stochastic analysis. Another unexpected byproduct is the identification of quasi-free modular states (i.e. thermal) as the natural arena for non-commutative stochastic analysis. We plan to address in future work the consequences of this observation outside the Grassmann setting, both for CAR and CCR algebras.

**Acknowledgments.** We would like to thank A. Chandra and M. Peev for interesting conversations on the subject of Grassmann SPDEs and for sharing their preprint [CHP23]. We are deeply grateful to S. Albeverio for his constant support and detailed feedback. The work of F.D. has been partially supported by Istituto Nazionale di Alta Matematica (INdAM, gruppo di ricerca GNAMPA). The work of L.F. has been partially supported by the Swiss National Science Foundation (SNSF), grant number 200160. M. Gordina's research was supported in part by NSF grant DMS-2246549. This work has been partially funded by the German Research Foundation (DFG) under Germany's Excellence Strategy - GZ 2047/1, project-id 390685813.

# 2 Non-commutative $L^p$ spaces

In this section, we briefly recall Haagerup's construction of the non-commutative $L^p$ spaces associated with a general vNa [Haa79a], see also [Ter81] for a detailed exposition. This construction is based on embedding $\mathcal{M}$ in a larger semifinite vNa obtained as a crossed product with the group $\mathbb{R}$, in which non-trivial elements of the $L^p$ spaces are intrinsically identified with affiliated unbounded operators. Haagerup's $L^p$ spaces possess desirable properties, but still lack some crucial features needed for applications to stochastic calculus. In Subsection 2.3 we introduce some novel twisted spaces that address this difficulty and will serve our purposes throughout the paper. In the final part of this section, we discuss conditional expectation and martingales inequalities for $L^p$ spaces and their twisted counterpart.

## 2.1 Modular theory and crossed products

The modular theory by Tomita and Takesaki [Tom67, Tak70, Hia03, Tak02b, KR97, Str81] is the key that opened the way to understanding factors without traces, that is, of type III, see also [Mas18]. Let $\mathcal{M}$ be a $\sigma$-finite vNa acting on $\mathcal{H}$, that is, $\mathcal{M}$ has a faithful normal state $\omega$ and, via the GNS representation, it is isomorphic to a vNa that has a cyclic and separating vector [BR87][2.1]. Without loss of generality, we can simply assume that $\omega(\cdot) = \langle \Omega, \cdot\, \Omega \rangle$, with $\Omega \in \mathcal{H}$ cyclic and separating for $\mathcal{M}$, so that $x \mapsto x\,\Omega$ establishes a bijection between $\mathcal{M}$ and the dense subspace $\mathcal{M}\,\Omega \subset \mathcal{H}$. The analysis of the operator acting on $\mathcal{H}$ associated with the involution $^*$ is the starting point of the Tomita–Takesaki modular theory. We let $S_\omega$ be the closure of the anti-linear operator on $\mathcal{M}\,\Omega$ defined by $S_\omega x\,\Omega := x^*\,\Omega$ and let $S_\omega =: J_\omega\,\Delta_\omega^{1/2}$ be its polar decomposition. Here $J_\omega$ is an anti-unitary operator called *modular conjugation* with the property $J_\omega^2 = \mathbb{1}$, whereas $\Delta_\omega$ is a positive non-singular self-adjoint operator called *modular operator*. If $\omega$ is a trace, then $\|S_\omega x\,\Omega\|_{\mathcal{H}}^2 = \|x\,\Omega\|_{\mathcal{H}}^2$, that is, $S_\omega$ is an isometry, so that $S_\omega = J_\omega$ and $\Delta_\omega = \mathbb{1}$. More generally, the modular operator reflects the non-tracial character in some way, as will be made precise below.

2.1. In the GNS construction $(\mathcal{H}_\omega, \pi_\omega, \Omega_\omega)$ $\omega(\cdot) = \langle \Omega_\omega, \cdot\, \Omega_\omega \rangle$ this means that $\Omega_\omega$ is cyclic and separating for $\pi_\omega(\mathcal{M})$.

A crucial result is Tomita's fundamental theorem [Tom67, Tak70].

**Theorem 2.1.** *(Tomita 1967)* $J_\omega \mathscr{M} J_\omega = \mathscr{M}'$ *and* $\Delta_\omega^{\mathrm{i}t} \mathscr{M} \Delta_\omega^{-\mathrm{i}t} = \mathscr{M}$ *for any* $t \in \mathbb{R}$.

As a consequence, one can introduce a one-parameter group of $^*$-automorphism $\mathbb{R} \ni t \mapsto \sigma_t^\omega$ called the *modular automorphism group* associated with the pair $(\mathscr{M}, \omega)$ and defined by

$$\sigma_t^\omega(x) = \Delta_\omega^{\mathrm{i}t} x \Delta_\omega^{-\mathrm{i}t}.$$

If $\Delta_\omega$ is bounded, then one has the stronger condition $\Delta_\omega^z \mathscr{M} \Delta_\omega^{-z} = \mathscr{M}$ for any $z \in \mathbb{C}$. Otherwise, it is useful to introduce the set of analytic elements in the vNa.

**Definition 2.2.** *The set of entire elements* $\mathscr{M}_a \subset \mathscr{M}$ *is such that* $x \in \mathscr{M}_a$ *iff* $z \mapsto \sigma_z^\omega(x)$ *extends to an entire* $\mathscr{M}$*-valued function.*

Note that $\mathscr{M}_a$ is a $w^*$-dense and $\sigma^\omega$-invariant $^*$-sub-algebra of $\mathscr{M}$. The modular automorphism group makes the non-traciality of $\omega$ quantitative precise via the KMS condition (or, more precisely, the $\beta = -1$ KMS condition [BR97]):

$$\omega(x \sigma_{-\mathrm{i}}^\omega(y)) = \omega(yx), \qquad \forall x, y \in \mathscr{M}_a,$$

showing again that $\omega$ is a trace iff $\Delta_\omega = \mathbb{1}$.

We now introduce the cross product of the vNa with the modular automorphism group [Tak73, Con73, BR87]. This larger vNa turns is at the heart of Haagerup's construction of the non-commutative $L^p$ spaces associated with an arbitrary $\sigma$-finite vNa. Such non-commutative integration spaces were first introduced by Dixmier, Segal and Kunz in the tracial setting [Seg53, Dix53, Kun58] and later extended to an arbitrary vNa by Haagerup [Haa79a]. See also the work of Kosaki [Kos80], of Araki–Masuda [AM82] for alternative definitions and the work of Trunov–Sherstnev [TS78a, TS78b] for preliminary work on $L^1$ spaces for general vNas.

Given the pair $(\mathscr{M}, \omega)$ with $\omega$ modular, the crossed product $\tilde{\mathscr{M}} := \mathscr{M} \otimes_{\sigma^\omega} \mathbb{R}$ is the vNa acting on $L^2(\mathbb{R}, \mathscr{H})$ and generated by the operators $(\pi(x))_{x \in \mathscr{M}}$ and $(\lambda(s))_{s \in \mathbb{R}}$ defined by

$$(\pi(x)\xi)(t) = \sigma_{-t}^\omega(x)\xi(t), \qquad (\lambda(s)\xi)(t) = \xi(t-s), \qquad \forall \xi \in L^2(\mathbb{R}, \mathscr{H}), \quad a.e. \quad t \in \mathbb{R}.$$

Note that $\pi$ is a normal faithful representation of $\mathscr{M}$ on $L^2(\mathbb{R}, \mathscr{H})$ and $\lambda$ gives the modular automorphism group

$$\pi(\sigma_t^\omega(x)) = \lambda(t)\pi(x)\lambda(t)^*, \qquad x \in \mathscr{M}, \quad t \in \mathbb{R}, \tag{2.1}$$

We will henceforth identify $\mathscr{M} \equiv \pi(\mathscr{M}) \subset \mathcal{R}$. One can also define a dual automorphism (dual via Fourier transform): since $(\lambda(s))_{s \in \mathbb{R}}$ is the group of translations, one lets $(W(t))_{t \in \mathbb{R}}$ be the unitary representation of $\mathbb{R}$ in $L^2(\mathbb{R}, \mathscr{H})$ given by $W(t)\xi(s) = \mathrm{e}^{-\mathrm{i}ts}\xi(s)$ and define

$$\hat{\sigma}_t(x) = W(t) x W(t)^*, \qquad x \in \mathscr{R}, t \in \mathbb{R}.$$

The following relations hold

$$\hat{\sigma}_t(x) = x, \qquad \hat{\sigma}_t(\lambda(s)) = \mathrm{e}^{-ist}\lambda(s), \qquad x \in \mathscr{M}, \ s, t \in \mathbb{R}$$

which likewise determine the action of $\hat{\sigma}_t$ uniquely. In particular this allows one to identify $\mathscr{M} \equiv \pi(\mathscr{M})$ as the subset of $\tilde{\mathscr{M}}$ invariant under $(\hat{\sigma}_t)_{t \in \mathbb{R}}$. For further details on the cross product construction for a general $W^*$-dynamical system, see [BR87].

It is well-known [Tak73, Con73] that the crossed product $\tilde{\mathscr{M}}$ is semifinite and it has a normal semifinite faithful trace $\tau$ such that $\tau \circ \hat{\sigma}_t = \mathrm{e}^{-t}\tau$ for any $t \in \mathbb{R}$. Let $\overline{P}(\mathscr{M})$ and $\overline{P}(\tilde{\mathscr{M}})$ denote the set of normal semifinite weights on $\mathscr{M}$ and $\tilde{\mathscr{M}}$ respectively [Hia20]. Any $\varphi \in \overline{P}(\mathscr{M})$ induces a dual $\tilde{\varphi} \in \overline{P}(\tilde{\mathscr{M}})$ which admits a Radon–Nikodym derivative $h_\varphi$ with respect to $\tau$ [Haa79b],

$$\tilde{\varphi}(\cdot) = \tau(h_\varphi \cdot).$$

The mapping $\varphi \mapsto \tilde{\varphi}$ is a bijection of the set of normal semifinite weights on $\mathcal{M}$ to the subset of normal semifinite weights on $\tilde{\mathcal{M}}$ such that $\tilde{\varphi} \circ \hat{\sigma}_t = \tilde{\varphi}$, $\forall t \in \mathbb{R}$, see [Ter81]. This allows one to obtain a mapping $\overline{P}(\mathcal{M}) \ni \varphi \mapsto h_\varphi \in \tilde{\mathcal{M}}$ with the following properties[2.2]:

**Theorem 2.3.** *The mapping $\varphi \mapsto h_\varphi$ extends to a bijection from the predual $\mathcal{M}_*$ to the subspace $\{h \in \tilde{\mathcal{M}} \mid \forall t \in \mathbb{R} \quad \hat{\sigma}_t(h_\varphi) = \mathrm{e}^{-t} h_\varphi\}$ and satisfies, for any $\varphi \in \mathcal{M}_*$ and $x, y \in \mathcal{M}$:*

$$h_{x\varphi y} = x h_\varphi y, \qquad h_{\varphi^*} = h_\varphi^*, \qquad |h_\varphi| = h_{|\varphi|}. \tag{2.2}$$

**Remark 2.4.** Note that for any $\varphi \in \mathcal{M}_*$, one can always define its modulus $|\varphi| \in \mathcal{M}_*$ by duality, see [Tak02a].

In particular, the dual state $\tilde{\omega}$ has Radon–Nikodym derivative with respect to $\tau$ which, following standard notation in the literature, we denote by $D = h_\omega$ (and exclusively reserve $D$ for this purpose) so that on $\tilde{\mathcal{M}}_+$ one has[2.3]

$$\tilde{\omega}(\cdot) = \tau(D \cdot).$$

It is important to recall that $D$ is a positive invertible self-adjoint operator on $L^2(\mathbb{R}, \mathcal{H})$ affiliated with $\tilde{\mathcal{M}}$ and

$$\lambda(t) = D^{\mathrm{i}t}, \qquad t \in \mathbb{R}$$

which follows from the fact that $\tilde{\mathcal{M}}$ is semifinite.

## 2.2 Haagerup's $L^p$ spaces

Because $\tau$ is a trace on $\tilde{\mathcal{M}}$, one can form the topological algebra $L^0(\tilde{\mathcal{M}}, \tau)$ of $\tau$-measurable operators affiliated with $\tilde{\mathcal{M}}$, see [Nel74, Ter81]. Note that $L^0(\tilde{\mathcal{M}}, \tau)$ also includes unbounded operators. Haagerup's spaces $L^p(\mathcal{M}, \omega)$ are defined as follows [Haa79a]:

**Definition 2.5.** *(Haagerup's $L^p$ spaces). Set*

$$\begin{aligned} L^p(\mathcal{M}, \omega) &:= \left\{x \in L^0(\tilde{\mathcal{M}}, \tau) \mid \hat{\sigma}_t(x) = \mathrm{e}^{-\frac{t}{p}} x, \quad t \in \mathbb{R}\right\} \qquad p \in [1, \infty), \\ L^\infty(\mathcal{M}, \omega) &:= \left\{x \in L^0(\tilde{\mathcal{M}}, \tau) \mid \hat{\sigma}_t(x) = x, \quad t \in \mathbb{R}\right\}. \end{aligned}$$

**Remark 2.6.** Note that $L^p(\mathcal{M}, \omega)$ is a $^*$-invariant linear subspace of $L^0(\tilde{\mathcal{M}}, \tau)$. By definition, we also have that $L^p(\mathcal{M}, \omega) \cap L^q(\mathcal{M}, \omega) = \{0\}$ if $p \neq q$.

By the considerations in the previous section we have that $L^\infty(\mathcal{M}, \omega) = \mathcal{M}$ and $L^1(\mathcal{M}, \omega) \cong \mathcal{M}_*$, so that we can equip $L^1(\mathcal{M}, \omega)$ with the norm of $\mathcal{M}_*$. It is then evident that, up to an isometry, for $p = 1, \infty$ $L^p(\mathcal{M}, \omega)$ do not depend on $\omega$, hence it can be dropped from the notation. This can actually be done for any $p \in [1, \infty]$.

**Proposition 2.7.** *Let $x \in L^0(\mathcal{R}, \tau)$ with the polar decomposition $x = u|x|$ and let $p \in [1, \infty)$. Then, $x \in L^p(\mathcal{M}, \omega)$ iff $u \in L^\infty(\mathcal{M})$ and $|x|^p \in L^1(\mathcal{M})$.*

We shall henceforth write $L^p(\mathcal{M}) \equiv L^p(\mathcal{M}, \omega)$ for any $p \in [1, \infty]$. As a consequence of Proposition 2.7, the norm on $\mathcal{M}_*$ then induces a norm on the $L^p(\mathcal{M})$ spaces; this can also be done via Haagerup's trace.

2.2. Recall that the polar decomposition of $\varphi \in \mathcal{M}_*$ is obtained via duality, see [Tak02a].

2.3. In general, the expression should be intended as a suitable regularization of $\tilde{\omega}(x) = \tau\left(D^{\frac{1}{2}} x D^{\frac{1}{2}}\right)$, for any $x \in \tilde{\mathcal{M}}^+$, see [PT73], unless $\omega$ is finite, that is, $\omega(\mathbb{1}) < \infty$.

**Definition 2.8.** *(Haagerup's trace). On $L^1(\mathcal{M})$ define the following (positive contractive) linear functional:*

$$\mathrm{tr_H}(h_\varphi) := \varphi(\mathbb{1}).$$

By the arguments presented in the previous section, we can recover our distinguished state $\omega$ as follows

$$\omega(x) = \mathrm{tr_H}(Dx) \qquad \forall x \in \mathcal{M},$$

where $D = h_\omega$. This functional is suggestively called a trace because for any $x \in L^p(\mathcal{M})$ and $y \in L^q(\mathcal{M})$ with $p, q \in [1, \infty]$ such that $1/p + 1/q = 1$, we have $xy, yx \in L^1(\mathcal{M})$ and

$$\mathrm{tr_H}(xy) = \mathrm{tr_H}(yx).$$

In particular we also have that for $x \in L^1(\mathcal{M})$ and $u \in L^\infty(\mathcal{M})$ unitary $\mathrm{tr_H}(uxu^*) = \mathrm{tr_H}(x)$, and for $y \in L^2(\mathcal{M})$ $\mathrm{tr_H}(y^*y) = \mathrm{tr_H}(yy^*)$.

The property in (2.2) implies that $\mathrm{tr_H}(|h_\varphi|) = \mathrm{tr_H}(h_{|\varphi|}) = |\varphi|(\mathbb{1}) \equiv \|\varphi\|_{\mathcal{M}_*}$, so that $x \mapsto \mathrm{tr_H}(|x|)$ is precisely the norm on $L^1(\mathcal{M})$ induced by $\|\cdot\|_{\mathcal{M}_*}$. We can then define a norm on $L^p(\mathcal{M})$ as follows.

**Definition 2.9.** *On $L^p(\mathcal{M})$ for $p \in [1, \infty)$ define*

$$\|x\|_{L^p(\mathcal{M})} := (\mathrm{tr_H}(|x|^p))^{\frac{1}{p}} \qquad \forall x \in L^p(\mathcal{M}).$$

*We also set $\|\cdot\|_{L^\infty(\mathcal{M})} := \|\cdot\|$.*

If no confusion arises we abridge $\|\cdot\|_{L^p(\mathcal{M})}$ to $\|\cdot\|_p$. One can prove the following:

**Theorem 2.10.** *(Hölder's inequality). Let $p, q, r \in [1, \infty]$ such that $1/p + 1/q = 1/r$ and let $x \in L^p(\mathcal{M})$ and $y \in L^q(\mathcal{M})$. Then, $xy \in L^r(\mathcal{M})$ and*

$$\|xy\|_r \leqslant \|x\|_p \, \|y\|_q.$$

*Furthermore, for $p \in [1, \infty)$ $(x, y) \mapsto \mathrm{tr_H}(xy)$ defines a duality between $L^p(\mathcal{M})$ and $L^q(\mathcal{M})$, $(L^p(\mathcal{M}))^* = L^q(\mathcal{M})$ isometrically.*

In particular, this duality allows one to prove Minkowski's inequality so that $\|\cdot\|_p$ is indeed a norm. The space $L^p(\mathcal{M})$ is complete in the said norm.

**Theorem 2.11.** *$(L^p(\mathcal{M}), \|\cdot\|_p)$ is Banach space for any $p \in [1, \infty]$. In particular, $L^2(\mathcal{M})$ is a Hilbert space with scalar product $\langle x, y\rangle_{L^2(\mathcal{M})} := \mathrm{tr_H}(x^*y)$.*

Haagerup's $L^p$ spaces are uniformly convex [Kos84], see also [Ter81] for a proof of Clarkson's inequality.

**Proposition 2.12.** *$(L^p(\mathcal{M}), \|\cdot\|_p)$ is uniformly convex for any $p \in (1, \infty)$.*

We have learnt that the only non-trivial elements in $L^p(\mathcal{M})$ are unbounded operators if $p \neq \infty$. The operator $D$ is unbounded and it is clear that $D^{\frac{1}{p}} \in L^p(\mathcal{M})$. We also have the following stronger result, see [JX03].

**Lemma.** *For any $p \in [1, \infty)$ $\mathcal{M}_a D^{\frac{1}{p}}$ is dense in $L^p(\mathcal{M})$ and we have $D^{\frac{1}{2p}+\tau} \mathcal{M}_a D^{\frac{1}{2p}-\tau} = \mathcal{M}_a D^{\frac{1}{p}}$, for any $\tau \in \mathbb{R}$.*

**Proof.** By Hölder's inequality $\mathcal{M}_a D^{\frac{1}{p}} \subset L^p$. The density is proved by duality: let $y \in L^q$ with $1 = 1/p + 1/q$ be such that $\mathrm{tr}_{\mathrm{H}}\left(x D^{\frac{1}{p}} y\right) = 0$ for any $x \in \mathcal{M}_a$. By Theorem 2.10, $D^{\frac{1}{p}} y \in L^1$ and $\mathcal{M} = (L^1)^*$. Since $\mathcal{M}_a$ is $w^*$-dense in $\mathcal{M}$ we have $D^{\frac{1}{p}} y = 0$ and thus $y = 0$ by the invertibility of $D$. To prove the identity, we note that if $\forall x \in \mathcal{M}_a$ then $[x]_{-\frac{1}{2p}-\tau} \in \mathcal{M}_a$ for any $\tau \in \mathbb{R}$ and thus $D^{\frac{1}{2p}+\tau}\, \mathcal{M}_a D^{\frac{1}{2p}-\tau} \ni D^{\frac{1}{2p}+\tau}\, [x]_{-\frac{1}{2p}-\theta} D^{\frac{1}{2p}-\tau} = x D^{\frac{1}{p}}$, that is, $\mathcal{M}_a D^{\frac{1}{p}} \subset D^{\frac{1}{2p}+\tau}\, \mathcal{M}_a D^{\frac{1}{2p}-\tau}$. The other inclusion is proved in the same way. □

The dense subspaces $\mathcal{M}_a D^{\frac{1}{p}}$ will play a crucial role in our non-commutative stochastic calculus and, more generally, in the calculus of unbounded operators, see Appendix B. This is not at all surprising, given their importance for actual computations, see, e.g. [Jun02, JX03, JX07].

Before delving into the details of our construction, let us conclude this section by recalling the special case in which $\mathcal{M}$ is semifinite. In this case, $\mathcal{M}$ has a faithful trace, see, e.g., [Hia20], so that one can introduce the tracial $L^p$ spaces [Seg53, Dix53, Kun58]. Let $\tau$ denote the faithful trace on $\mathcal{M}$. Recall that $\mathcal{M} \otimes_{\sigma^\tau} \mathbb{R}$ acts on $L^2(\mathbb{R}, \mathcal{H}) \cong \mathcal{H} \otimes L^2(\mathbb{R})$. Let $\mathcal{F}: L^2(\mathbb{R}) \to L^2(\mathbb{R})$ denote the Fourier transform and, for any $f$ measurable function on $\mathbb{R}$, let $m(f)$ denote the multiplication operator by $f$ on $L^2(\mathbb{R})$. We have

$$\mathcal{M} \otimes_{\sigma^\tau} \mathbb{R} = \mathcal{M} \otimes \mathcal{F}^{-1} m(L^\infty(\mathbb{R}))\mathcal{F}.$$

**Proposition 2.13.** *Let $\mathcal{M}$ be semifinite, let $\tau$ be the faithful trace on it and let $L^p(\mathcal{M}, \tau)$ denote the tracial $L^p$ space for $p \in [1, \infty)$. Then, we have the isometry*

$$L^p(\mathcal{M}) \cong L^p(\mathcal{M}, \tau) \otimes \mathcal{F}^{-1} m\left(\mathrm{e}^{\frac{\cdot}{p}}\right)\mathcal{F}.$$

Note that for $p \in [1, \infty)$ $\mathrm{e}^{\frac{\cdot}{p}}$ is not bounded and in fact $L^p(\mathcal{M})$ contains only unbounded non-vanishing operators on $\mathcal{H} \otimes L^2(\mathbb{R})$.

## 2.3 Twisted $L^p$ spaces

Even though Haagerup's $L^p$ spaces are very general and elegant, they turn out not to be so convenient for constructing a general non-commutative stochastic calculus. This is ultimately due to the following problems. First of all, unlike the tracial $L^p$ spaces, a scale is completely missing: in fact, we have $L^p(\mathcal{M}) \cap L^q(\mathcal{M}) = \{0\}$ if $p \neq q$. This is intrinsically due to the crossed-product structure, as can be seen directly when $\mathcal{M}$ is semifinite, see Proposition 2.13.

One natural way to recover such a property is to consider the closure of the injection of $\mathcal{M}_a$ into $L^p$, e.g., via $\mathcal{M}_a D^{\frac{1}{p}}$. In fact, if we have $(a_n)_n \subset \mathcal{M}_a$ such that $a_n D^{\frac{1}{p}} \to A$ in $L^p$, then by Hölder's inequality it is straightforward to see that $a_n D^{\frac{1}{q}}$ is convergent in $L^q$ for any $q < p$. Note that in [MZ96] a similar construction was adopted by considering the symmetric injection $\mathcal{M}_a \to D^{\frac{1}{2p}}\, \mathcal{M}_a D^{\frac{1}{2p}}$ into $L^p$.

Yet, this simple Ansatz is inconvenient when considering products. To see this, let $(a_n), (b_n) \subset \mathcal{M}_a$ be such that $a_n D^{\frac{1}{p}}$ and $b_n D^{\frac{1}{q}}$ converge in the $L^p$ and $L^q$ topologies respectively. But $a_n b_n D^{\frac{1}{r}} = a_n D^{\frac{1}{p}}\, \sigma_{\mathrm{i}\frac{1}{p}}(b_n)\, D^{\frac{1}{q}}$, where $\sigma_z$ is the analytic continuation of the automorphism group on $\mathcal{M}_a$, so that the product sequence can be controlled only if we can control the twisted sequence $\sigma_{\mathrm{i}\frac{1}{p}}(b_n)\, D^{\frac{1}{q}}$ in the $L^q$ topology. The same issue appears when taking the expectation (and conditional expectations) as well : by the definition of Haagerup's trace, one has $\omega(a_n b_n) = \mathrm{tr}_{\mathrm{H}}\left(a_n D^{\frac{1}{p}}\, \sigma_{\mathrm{i}\frac{1}{p}}(b_n)\, D^{\frac{1}{q}}\right)$ with $1 = 1/p + 1/q$, hence the control of a twisted sequence is needed once more.

This discussion justifies the introduction of some novel “twisted” $L^p$ spaces. First of all, let us denote the analytic continuation of the automorphism group on $\mathcal{M}_a$ by $[\cdot]_t$, that is, $[x]_t := \sigma_{-\mathrm{i}t}(x)$. Then, we introduce the curve of embeddings $T_t^{(p)}: \mathcal{M}_a \to L^p(\mathcal{M})$ defined by

$$T_t^{(p)}(x) := D^{\frac{1}{2p}}[x]_t D^{\frac{1}{2p}}, \tag{2.3}$$

satisfying $T_t^{(p)}(x^*) = (T_{-t}^{(p)}(x))^*$. This is indeed an embedding by Lemma 2.12:

**Corollary 2.14.** *The map $T_t^{(p)}$ has kernel $\{0\}$ and $T_t^{(p)}(\mathcal{M}_a)$ is dense in $L^p(\mathcal{M})$.*

The triviality of the kernel is a straightforward consequence of the invertibility of $D$. We can now introduce the twisted $L^p$ spaces, denoted by $\mathbb{L}^p(\mathcal{M}, \omega)$, as follows:

**Definition 2.15.** *Let $\mathcal{M}$ be a vNa, let $\omega$ be a nsf state on $\mathcal{M}$ and let $T_t^{(p)}$ be the embedding of $\mathcal{M}_a$ in $L^p(\mathcal{M})$ as defined in (2.3). We denote by $\mathbb{L}^p(\mathcal{M}, \omega)$ the completion of $\mathcal{M}_a$ with respect to the norm*

$$\|x\|_{\mathbb{L}^p(\mathcal{M},\omega)} := \sup_{|t| \leqslant 1-\frac{1}{2p}} \|T_t^{(p)}(x)\|_{L^p(\mathcal{M})}.$$

*We abridge the notation to $\mathbb{L}^p \equiv \mathbb{L}^p(\mathcal{M}, \omega)$ and $\|\cdot\|_{\mathbb{L}^p} = \|\cdot\|_{\mathbb{L}^p(\mathcal{M},\omega)}$ if no confusion arises.*

This definition is clearly meaningful because $T_t^{(p)}$ is an embedding and because $\|\cdot\|_{L^p(\mathcal{M})}$ is a norm. With abuse of notation, we henceforth extend the map $T_t^{(p)}$ to the the $\mathbb{L}^p$ spaces as well, $T_t^{(p)}: \mathbb{L}^p(\mathcal{M}, \omega) \to L^p(\mathcal{M})$. Therefore, for any $x_n \to x$ in $\mathbb{L}^p$, we let

$$T_t^{(p)}(x) := L^p - \lim_{n\to\infty} T_t^{(p)}(x_n),$$

for any $|t| \leqslant 1 - \frac{1}{2p}$.

We note the following properties.

**Lemma 2.16.** *We have:*

*i. For any $p \in [1, \infty]$ $\mathbb{L}^p(\mathcal{M})$ is a $^*$-invariant subspace of $L^0(\mathcal{M}, \omega)$ and $\|x\|_{\mathbb{L}^p} = \|x^*\|_{\mathbb{L}^p}$, $\forall x \in \mathbb{L}^p(\mathcal{M})$.*

*ii. Let $p, q, r \in [1, \infty]$ such that $1/p + 1/q = 1/r$. Then, for any $x, y \in \mathcal{M}_a$ we have*

$$T_t^{(r)}(xy) = T_{t+\frac{1}{2q}}^{(p)}(x)\, T_{t-\frac{1}{2p}}^{(q)}(y), \tag{2.4}$$

*and*

$$\|xy\|_{\mathbb{L}^r} \leqslant \|x\|_{\mathbb{L}^p} \|y\|_{\mathbb{L}^q}. \tag{2.5}$$

**Proof.** The first claim follows from the fact that Haagerup's $L^p$ spaces have that property and $T_t^{(p)}(x^*) = (T_{-t}^{(p)}(x))^*$, whereas (2.4) is a straightforward computation:

$$T_t^{(r)}(xy) = D^{\frac{1}{2q}+\frac{1}{2p}+t} x D^{-\frac{1}{2q}+\frac{1}{2p}-t} D^{\frac{1}{2q}-\frac{1}{2p}+t} y D^{\frac{1}{2q}+\frac{1}{2p}-t}.$$

Finally, we note that if $|t| \leqslant 1 - \frac{1}{2r}$, then $-1 + \frac{1}{2r} + \frac{1}{2q} \leqslant t + \frac{1}{2q} \leqslant 1 - \frac{1}{2p}$, and $-1 + \frac{1}{2q} \leqslant t - \frac{1}{2p} \leqslant 1 - \frac{1}{2r} - \frac{1}{2p}$, that is,

$$\left|t + \frac{1}{2q}\right| \leqslant 1 - \frac{1}{2p}, \qquad \left|t - \frac{1}{2p}\right| \leqslant 1 - \frac{1}{2q},$$

and thus also (2.5) holds true by Hölder's inequality for Haagerup's $L^p$ spaces. □

One would like to extend (2.5) to the whole twisted spaces. To this end, by Lemma 2.16 we extend the product on $\mathcal{M}_a \times \mathcal{M}_a \to \mathcal{M}_a$ to a product on $\mathbb{L}^p \times \mathbb{L}^q$.

**Definition 2.17.** *Define the product* $\cdot : \mathbb{L}^p \times \mathbb{L}^q \to \mathbb{L}^r$, *for any* $p,q,r \in [1,\infty]$ *with* $1/p + 1/q = 1/r$

$$(x,y) \mapsto x \cdot y := \mathbb{L}^r - \lim_{n\to\infty} x_n y_n,$$

*for any* $(x_n), (y_n) \subset \mathcal{M}_a$ *such that* $x = \mathbb{L}^p - \lim_{n\to\infty} x_n$ *and* $y = \mathbb{L}^q - \lim_{n\to\infty} y_n$.

**Remark 2.18.** This product is well-defined thanks to Lemma 2.16. In fact, it is straightforward to check that $x_n y_n$ is a Cauchy sequence in $\mathbb{L}^r$ and that furthermore the limit does not depend on the sequences $(x_n), (y_n) \subset \mathcal{M}_a$.

**Proposition 2.19.** *We have:*

*i. Let* $p,q,r \in [1,\infty]$ *such that* $1/p + 1/q = 1/r$. *Then, for any* $x \in \mathbb{L}^p(\mathcal{M})$ *and* $y \in \mathbb{L}^q(\mathcal{M})$, $x \cdot y \in \mathbb{L}^r(\mathcal{M})$ *and Hölder's inequality holds true*

$$\|x \cdot y\|_{\mathbb{L}^r} \leqslant \|x\|_{\mathbb{L}^p} \|y\|_{\mathbb{L}^q}. \tag{2.6}$$

*ii.* $\mathbb{L}^p(\mathcal{M}) \subset \mathbb{L}^q(\mathcal{M})$ *for any* $1 \leqslant q \leqslant p \leqslant \infty$ *and* $\|x\|_{\mathbb{L}^q} \leqslant \|x\|_{\mathbb{L}^p}$ *for any* $x \in \mathbb{L}^p(\mathcal{M})$.

**Proof.** The first claim is a direct consequence of Lemma 2.16 and Definition 2.17. The second statement follows from (2.6) since $\|\mathbb{1}\|_{\mathbb{L}^p(\mathcal{M})} = (\mathrm{tr}_{\mathrm{H}}(D))^{\frac{1}{p}} = \omega(\mathbb{1}) = 1$ for any $p \in [1,\infty]$. □

**Remark 2.20.** Note that if $x \in \mathbb{L}^p(\mathcal{M})$ and $y \in \mathbb{L}^q(\mathcal{M})$, then their product $x \cdot y$ can be made sense of as the $\mathbb{L}^r$-limit for any $1/r \geqslant 1/p + 1/q$, that is, for any $x_n \to x$ and $y_n \to y$ the sequence $x_n y_n$ is convergent in any such $\mathbb{L}^r$ topology. Furthermore, the embedding map $T_t^{(p)}$ can actually be extended to a map from $\mathbb{L}^q(\mathcal{M}) \to L^p(\mathcal{M})$ for any $q \geqslant p$.

Thanks to Proposition 2.19, we can also extend the expectation $\omega$ to any $\mathbb{L}^p$. In fact, the following lemma holds true.

**Lemma 2.21.** *For any* $p \in [1,\infty]$, *if* $x \in \mathbb{L}^p$ *the limit* $\lim_{n\to\infty} \omega(x_n)$ *exists for any* $(x_n)_n \subset \mathcal{M}_a$ *such that* $x = \mathbb{L}^p - \lim_{n\to\infty} x_n$ *and is independent of the sequence* $(x_n)_n$.

**Proof.** We have $|\omega(x_n) - \omega(x_m)| \leqslant \|x_n - x_m\|_1 \leqslant \|x_n - x_m\|_{\mathbb{L}^1} \leqslant \|x_n - x_m\|_{\mathbb{L}^p}$, where in the last step we used Proposition 2.19, so that $\omega(x_n)$ is a Cauchy sequence and thus converges. If $(y_n)_n \subset \mathcal{M}_a$ is another sequence such that $y_n \to x$ in the $\mathbb{L}^p$ topology, then $|\omega(x_n) - \omega(y_n)| \leqslant \|x_n - y_n\|_{\mathbb{L}^p}$, so that the limit does not depend on the sequence. □

## 2.4 Conditional expectation and $L^p$ processes

Let $\mathcal{M}$ be a $\sigma$-finite vNa, $\omega$ a modular state on it and $\sigma_t = \sigma_t^\omega$ its modular automorphism group. Given any vN subalgebra $\mathcal{N}$, a conditional expectation can be defined as a positive linear map $\omega_{\mathcal{N}} : \mathcal{M} \to \mathcal{N}$ such that

$$\omega_{\mathcal{N}}(\mathbb{1}_{\mathcal{M}}) = \mathbb{1}_{\mathcal{N}}, \qquad \omega_{\mathcal{N}}(x y z) = x\, \omega_{\mathcal{N}}(y)\, z,$$

for any $x, z \in \mathcal{N}$, see [Kad04a]. One can prove that such a map is a norm one projection and that it is completely positive (in particular it is a Schwarz map $\omega_{\mathcal{N}}(x^* x) \geqslant \omega_{\mathcal{N}}(x^*)\, \omega_{\mathcal{N}}(x)$), see [Tom57, Hia20]. Note that more general definitions exist, but they are beyond the scope of our work, see, e.g., [AC82].

It is a well-known result by Takesaki [Tak72] that for any vN subalgebra $\mathcal{N} \subset \mathcal{M}$ that is invariant under the modular automorphism group, that is, $\sigma_t(\mathcal{N}) = \mathcal{N}$ $\forall t \in \mathbb{R}$, there exists a conditional expectation $\omega_{\mathcal{N}} : \mathcal{M} \to \mathcal{N}$ and it satisfies

$$\omega \circ \omega_{\mathcal{N}} = \omega, \qquad \sigma_t \circ \omega_{\mathcal{N}} = \omega_{\mathcal{N}} \circ \sigma_t. \tag{2.7}$$

One can extend the conditional expectation to a positive contractive projection from $L^p(\mathscr{M})$ to $L^p(\mathscr{N})$ such that $\omega_{\mathscr{N}}(xyz) = x\,\omega_{\mathscr{N}}(y)\,z$ for any $x,z\in\mathscr{N}$ and $y\in L^p(\mathscr{M})$, see [JX03]. This is a general result that can be applied to any completely positive contraction and is also relevant for proving hypercontractivity for suitable random variables. See [JX07] for a general construction of such extensions and [LR11] for its application to the proof of hypercontractivity within the $q$-deformed Araki–Woods factors, that is, algebras interpolating between the Araki–Wyss and the Araki–Woods ones as the deformation parameter $q$ goes from $-1$ to 1[2.4]. For the sake of completeness, we report these results here. Let $P\colon\mathscr{M}\to\mathscr{N}$ be a completely positive contraction that is also state preserving (that is, $\omega\circ P=\omega$). Since $\mathscr{M}D^{\frac{1}{p}}$ is dense in $L^p$, one can introduce the densely-defined operators, for $p\in[1,\infty)$

$$P^{(p)}\colon \mathscr{M}D^{\frac{1}{p}}\to\mathscr{N}D^{\frac{1}{p}}\qquad xD^{\frac{1}{p}}\mapsto P(x)\,D^{\frac{1}{p}}. \tag{2.8}$$

which extend to a contraction on $L^p(\mathscr{M})$, see [JX07, Section 7] for the proof.

One can apply this result to the extension of the conditional expectation to Haagerup's and to the twisted $L^p$ spaces. On $D^{\frac{1}{2p}+\tau}\,\mathscr{M}_a D^{\frac{1}{2p}-\tau}$ one can define

$$\omega_{\mathscr{N}}^{(p)}\left(D^{\frac{1}{2p}+\tau}xD^{\frac{1}{2p}-\tau}\right)=D^{\frac{1}{2p}+\tau}\,\omega_{\mathscr{N}}(x)\,D^{\frac{1}{2p}-\tau}. \tag{2.9}$$

Since $\omega_{\mathscr{N}}$ maps $\mathscr{M}_a$ into the entire elements $\mathscr{N}_a$, which can be seen from (2.7), one sees that because of Lemma 2.12 and the state preserving property (2.7) $\omega_{\mathscr{N}}^{(p)}$ does not depend on $\tau$, which was omitted from the notation on purpose. One can extend the map in (2.9) to the $L^p$ spaces and prove that such a map satisfies the usual properties [JX03].

**Proposition 2.22.** *The map $\omega_{\mathscr{N}}^{(p)}$ extends to a contractive projection, still denoted by $\omega_{\mathscr{N}}^{(p)}$, from $L^p(\mathscr{M})$ to $L^p(\mathscr{N})$. Furthermore, the following properties hold true:*

- *i. Let $p\in[1,\infty]$. Then, $\forall x\in L^p(\mathscr{M})$ $(\omega_{\mathscr{N}}^{(p)}(x))^*=\omega_{\mathscr{N}}^{(p)}(x^*)$ and $x\geqslant 0\Rightarrow\omega_{\mathscr{N}}^{(p)}(x)\geqslant 0$.*
- *ii. Let $p\in[2,\infty]$. Then, $\forall x\in L^p(\mathscr{M})$ $(\omega_{\mathscr{N}}^{(p)}(x))^*\,\omega_{\mathscr{N}}^{(p)}(x)\leqslant\omega_{\mathscr{N}}^{(p)}(x^*x)$.*
- *iii. Let $p,q,r\in[1,\infty]$ such that $1/p+1/q+1/r\leqslant 1$. Then, $\omega_{\mathscr{N}}^{(p)}(yxz)=y\,\omega_{\mathscr{N}}^{(p)}(x)\,z$ for any $x\in L^p(\mathscr{M})$, $y\in L^q(\mathscr{N})$ and $z\in L^r(\mathscr{N})$.*

**Remark 2.23.** In particular, since $T_t^{(p)}(\omega_{\mathscr{N}}(x))=\omega_{\mathscr{N}}^{(p)}(T_t^{(p)}(x))$ for $x\in\mathscr{M}_a$ this result allows us to extend the conditional expectation to the $\mathbb{L}^p$ spaces as well, extension which we shall denote by $\omega_{\mathscr{N}}$ with abuse of notation, that is $\omega_{\mathscr{N}}(x)=\mathbb{L}^p-\lim_n\omega_{\mathscr{N}}(x_n)$, for any $x_n\to x$ in the $\mathbb{L}^p$ topology. It follows that, see Remark 2.20,

$$T_t^{(p)}(\omega_{\mathscr{N}}(x))=\omega_{\mathscr{N}}^{(p)}(T_t^{(p)}(x))\qquad x\in\mathbb{L}^q,\ q\geqslant p.$$

**Remark 2.24.** Note the following duality relation, for $x\in L^p(\mathscr{M})$ and $y\in L^q(\mathscr{M})$ with $1/p+1/q=1$

$$\mathrm{tr}_{\mathrm{H}}(y\,\omega_{\mathscr{N}}^{(p)}(x))=\mathrm{tr}_{\mathrm{H}}(\omega_{\mathscr{N}}^{(q)}(y)\,x), \tag{2.10}$$

which is a consequence of the tower property and the density of $\mathscr{M}_aD^{\frac{1}{p}}$ in $L^p(\mathscr{M})$. In fact, we have

$$\mathrm{tr}_{\mathrm{H}}\left(D^{\frac{1}{q}}a\omega_{\mathscr{N}}^{(p)}\left(bD^{\frac{1}{p}}\right)\right)=\omega(a\,\omega_{\mathscr{N}}(b))=\omega(\omega_{\mathscr{N}}(a)\,b)=\mathrm{tr}_{\mathrm{H}}\left(\omega_{\mathscr{N}}^{(q)}\left(D^{\frac{1}{q}}a\right)bD^{\frac{1}{p}}\right),$$

from which (2.10) is obtained by continuity.

More generally, we will consider filtrations of vN subalgebras $(\mathscr{M}_t)_t$ invariant under the modular automorphism. In this case, it is useful to introduce the following nomenclature.

2.4. In Section 3 we will introduce Grassmann random variables within the Araki–Wyss factors.

**Definition 2.25.** *(Filtered modular space). The triple* $(\mathcal{M}, \omega, (\mathcal{M}_t)_{t\in\mathbb{R}_+})$ *consisting of a vNa* $\mathcal{M}$*, a modular state* $\omega$ *and a filtration of vNa* $\mathcal{M}_t$ *such that* $\bigcup_{t\geqslant 0}\mathcal{M}_t$ *is w*-dense in* $\mathcal{M}$ *is a filtered modular space if* $\mathcal{M}_t$ *are invariant under the modular group. In such a case, denote by* $\omega_t$: $\mathcal{M}\to\mathcal{M}_t$ *the conditional expectation and its extension to* $L^p(\mathcal{M}_t)$ *by* $\omega_t^{(p)}$.

In a filtered modular space, one can define adapted processes and martingales in the usual way. Of course, thanks to Proposition 2.22 these definition extends naturally to the $L^p$ and $\mathbb{L}^p$ spaces in a natural way, the former including $\mathcal{M}$ as a special case.

**Definition 2.26.** *Let* $(\mathcal{M}, \omega, (\mathcal{M}_t)_{t\in\mathbb{R}_+})$ *be a filtered modular space. We say that* $(x_t)_t$ *is an* $L^p$ *process if* $x_t\in L^p(\mathcal{M})$*, that is an adapted* $L^p$ *process if* $x_t\in L^p(\mathcal{M}_t)$ *for any* $t\in\mathbb{R}_+$ *and that it is an* $L^p$ *martingale if additionally* $\omega_s^{(p)}(x_t)=x_s$ *for any* $s\leqslant t$*. An* $L^p$ *process is bounded if* $\sup_{t\in\mathbb{R}_+}\|x_t\|_p<\infty$*. An* $L^p$ *martingale is finite if there exists* $t_0$ *such that* $x_t=x_{t_0}$ *for any* $t\geqslant t_0$*.*

**Remark 2.27.** The same definitions hold in the case of processes valued in $\mathcal{M}_a$ or in $\mathbb{L}^p$ with the obvious substitutions.

If $(t_n)_{n\geqslant 0}\subset\mathbb{R}_+$ is an unbounded set, having no accumulation point, such that $0=t_0<\cdots<t_n$, with abuse of notation, we will denote by $(\mathcal{M}_n)_{n\in\mathbb{N}_0}$ the filtration of vNa's $(\mathcal{M}_{t_n})_{n\in\mathbb{N}_0}$. We can then introduce simple adapted processes in $L^p$: $F$ is an $L^p$ simple adapted process if there exists an unbounded set $(t_n)_{n\geqslant 0}\subset\mathbb{R}_+$ as before, such that

$$F_t=\sum_{n\geqslant 0}F_{t_n}\mathbf{1}_{[t_n,t_{n+1})}(t),\qquad F_{t_n}\in L^p(\mathcal{M}_n). \tag{2.11}$$

Let us conclude this section by considering bounded $L^p$ martingales with respect to a discrete filtration $(\mathcal{M}_n)_{n\in\mathbb{N}_0}$.

**Lemma 2.28.** *If* $p\in(1,\infty)$ *any bounded* $L^p$ *martingales* $x=(x_n)_{n\in\mathbb{N}_0}$ *has a limit* $x_\infty$ *in* $L^p$ *and* $x_n=\omega_n(x_\infty)$.

**Remark 2.29.** Conversely, note that if $x_\infty\in L^p$, $1\leqslant p<\infty$ then $x=(x_n)_n$ with $x_n=\omega_n^{(p)}(x_\infty)$ defines a bounded $L^p$ martingale and $x_n\to x_\infty$ in the $L^p$ topology. Thus, can identify a bounded $L^p$ martingale with $x_\infty$ and viceversa.

In [JX07] the proof is hinted, but we prefer to spell out the details, which differ significantly from the standard proof in the commutative setting via the Doob's theorem.

**Proof.** Since $x=(x_n)_n$ is bounded in $L^p$ it has a weakly convergent subsequence in $L^p$, $x_{n_k}\rightharpoonup X$. For any $n$ and any $Q_n\in L^q(\mathcal{M}_n)$ by (2.10) and by the fact that $x$ is a martingale we have

$$\mathrm{tr}_{\mathrm{H}}(Q_n(\omega_n^{(p)}(X)-x_n))=\mathrm{tr}_{\mathrm{H}}(Q_n(X-x_n))=\lim_{k\to\infty}\mathrm{tr}_{\mathrm{H}}(Q_n(X-x_{n_k}))=0.$$

Since this is true for all $Q_n\in L^q(\mathcal{M}_n)$, we have

$$x_n=\omega_n^{(p)}(X),\qquad\forall n. \tag{2.12}$$

Let us now assume that there is another subsequence $x_{m_k}\rightharpoonup Y$, with $Y\neq X$. Without loss of generality we may assume $\|X\|_p=1$ and $\|Y\|_p\leqslant 1$. If $\|X-Y\|\geqslant\varepsilon$, by uniform convexity of the $L^p$ spaces, see Proposition 2.12, we have that there exists $\delta>0$ such that $\|(X+Y)/2\|_p\leqslant 1-\delta$. Then for all $n>0$, by Proposition 2.22

$$\|x_n\|_p=\left\|\omega_n^{(p)}\left(\frac{X+Y}{2}\right)\right\|_p\leqslant\left\|\frac{X+Y}{2}\right\|_p\leqslant 1-\delta.$$

On the other hand, by duality we have weak lower semicontinuity

$$1=\|X\|_p\leqslant\liminf_k\|x_{n_k}\|_p\leqslant 1-\delta$$

which is a contradiction, therefore $X = Y$, that is, the weak limit is unique. Then, by (2.12) we have

$$\limsup_n \|x_n\|_p \leqslant \|X\|_p,$$

and again by the weak lower semicontinuity $1 = \|X\|_p \leqslant \liminf_n \left\|\frac{X+x_n}{2}\right\|_p$. Thus, $\left\|\frac{X+x_n}{2}\right\|_p \to 1$ and by uniform convexity $\|x_n - X\|_p \to 0$. □

## 2.5 Martingale inequalities

Martingale inequalities play a crucial role in defining the Itô integral for a sufficiently large class of integrands, see Section 4.1. In the tracial setting, they have been established in [PX97] and applied to the construction of the Itô–Clifford integral for suitable $L^p$ processes, extended the seminal result in $L^2$ [BSW82]. Because non-trivial Grassmann martingales cannot be represented in the tracial setting, one has to resort to martingale inequalities in Haagerup's $L^p$ spaces, which were established in [JX03]. Below, we recall these results and adapt some of the details to the twisted $\mathbb{L}^p$ setting.

Let $(\mathcal{M}, \omega, (\mathcal{M}_n)_{n\in\mathbb{N}_0})$ be a filtered modular space, see Definition 2.25, and denote by $\omega_n^{(p)}$ the extension of the conditional expectation $\omega_n$ to the spaces $L^p(\mathcal{M})$. We shall now introduce some Hardy spaces of non-commutative martingales as follows. Consider finite sequences $x := (x_n)_{n\geqslant 0} \subset L^p(\mathcal{M})$, with the topologies [PX97, JX03]

$$\begin{aligned}
\|x\|_{L^p(\mathcal{M};\ell_c^2)} &:= \left\|\left(\sum_{n\geqslant 0} |x_n|^2\right)^{1/2}\right\|_p, &\quad \|x\|_{L^p_{\mathrm{cond}}(\mathcal{M};\ell_c^2)} &:= \left\|\left(\sum_{n\geqslant 0} \omega_{n-1}^{(\frac{p}{2})}(|x_n|^2)\right)^{1/2}\right\|_p,\\
\|x\|_{L^p(\mathcal{M};\ell_r^2)} &:= \left\|\left(\sum_{n\geqslant 0} |x_n^*|^2\right)^{1/2}\right\|_p, &\quad \|x\|_{L^p_{\mathrm{cond}}(\mathcal{M};\ell_r^2)} &:= \left\|\left(\sum_{n\geqslant 0} \omega_{n-1}^{(\frac{p}{2})}(|x_n^*|^2)\right)^{1/2}\right\|_p,
\end{aligned} \tag{2.13}$$

with the identification $\omega_{-1}^{(p)} = \omega_0^{(p)}$. To be precise, in the case of $L^p_{\mathrm{cond}}(\mathcal{M};\ell_\sharp^2)$ for $1 \leqslant p < 2$ one should consider finite sequences in $\mathcal{M}_a D^{\frac{1}{p}}$, since the extension theorem for the conditional expectation does not hold for those values of $p$. The quantities in (2.13) define a norm on the space of finite sequences in $L^p(\mathcal{M})$ (or $\mathcal{M}_a D^{\frac{1}{p}}$ if needed) for any $p \geqslant 1$, see [Jun02, JX03]. Additionally, consider the space $\ell^p(L^p(\mathcal{M}))$ of sequences $x = (x_n)_n$ in $L^p(\mathcal{M})$ such that

$$\|x\|_{\ell^p(L^p(\mathcal{M}))} := \left(\sum_{n\geqslant 0} \|x_n\|_p^p\right)^{\frac{1}{p}} < \infty.$$

Define the difference sequence $\delta x := (\delta x_n)_{n\geqslant 0}$ with $\delta x_n := x_n - x_{n-1}$, where $x_{-1} \equiv 0$. For any finite $L^p$ martingale introduce the norms

$$\|x\|_{\mathcal{H}^p_\sharp(\mathcal{M})} := \|\delta x\|_{L^p(\mathcal{M},\ell_\sharp^2)}, \qquad \|x\|_{h^p_\sharp} := \|\delta x\|_{L^p_{\mathrm{cond}}(\mathcal{M},\ell_\sharp^2)}, \qquad \sharp = c, r, \qquad \|x\|_{h^p_d} := \|\delta x\|_{\ell^p(L^p(\mathcal{M}))}. \tag{2.14}$$

The closure of the space of finite $L^p$ martingales under the norms in (2.14) are denoted respectively by $\mathcal{H}^p_\sharp(\mathcal{M})$, $h^p_\sharp(\mathcal{M})$ and $h^p_d(\mathcal{M})$. The Hardy spaces of non-commutative martingales are then defined as follows.

**Definition 2.30.** *(Hardy spaces). For $p \in [1, 2)$ set*

$$\mathcal{H}^p(\mathcal{M}) := \mathcal{H}^p_c(\mathcal{M}) + \mathcal{H}^p_r(\mathcal{M}), \qquad h^p(\mathcal{M}) := h^p_d(\mathcal{M}) + h^p_c(\mathcal{M}) + h^p_r(\mathcal{M})$$

*with norms*

$$\begin{aligned}
\|x\|_{\mathcal{H}^p(\mathcal{M})} &:= \inf_{\substack{x = x_c + x_r\\ x_c \in \mathcal{H}^p_c(\mathcal{M}), x_r \in \mathcal{H}^p_r(\mathcal{M})}} (\|x_c\|_{\mathcal{H}^p_c(\mathcal{M})} + \|x_r\|_{\mathcal{H}^p_r(\mathcal{M})}),\\
\|x\|_{h^p(\mathcal{M})} &:= \inf_{\substack{x = x_d + x_c + x_r\\ x_d \in h^p_d(\mathcal{M}), y \in h^p_c(\mathcal{M}), z \in h^p_r(\mathcal{M})}} \left(\|x_r\|_{h^p_d(\mathcal{M})} + \|x_c\|_{h^p_c(\mathcal{M})} + \|x_r\|_{h^p_r(\mathcal{M})}\right).
\end{aligned}$$

*For* $p\in[2,\infty)$ *set*

$$\mathcal{H}^p(\mathscr{M}) := \mathcal{H}_c^p(\mathscr{M}) \cap \mathcal{H}_r^p(\mathscr{M}), \qquad h^p(\mathscr{M}) := h_d^p(\mathscr{M}) \cap h_c^p(\mathscr{M}) \cap h_r^p(\mathscr{M})$$

*with norms*

$$\|x\|_{\mathcal{H}^p(\mathscr{M})} := \max_{\sharp=c,r} \|x\|_{\mathcal{H}_\sharp^p(\mathscr{M})}, \qquad \|x\|_{h^p(\mathscr{M})} := \max_{\sharp=c,r,d} \|x\|_{h_\sharp^p(\mathscr{M})}.$$

The next theorem extends the Burkholder–Gundy and Burkholder's inequalities to the noncommutative setting of Haagerup's $L^p$ spaces, see [JX03].

**Theorem 2.31.** *(Junge–Xu, 2003). Let* $(\mathscr{M},\omega,(\mathscr{M}_n)_n)$ *be a filtered modular space. Let* $1<p<\infty$ *and let* $x=(x_n)_{n\geqslant 0}$ *be an* $L^p$*-martingale. Then,* $x$ *is bounded iff* $x\in\mathcal{H}^p(\mathscr{M})$ *and iff* $x\in h^p(\mathscr{M})$*. In this case, there exist constants* $\alpha_p,\beta_p,\delta_p,\eta_p>1$ *such that*

$$\alpha_p^{-1}\,\|x\|_{\mathcal{H}^p(\mathscr{M})} \leqslant \|x\|_p \leqslant \beta_p\,\|x\|_{\mathcal{H}^p(\mathscr{M})}, \tag{2.15}$$

$$\delta_p^{-1}\,\|x\|_{h^p(\mathscr{M})} \leqslant \|x\|_p \leqslant \eta_p\,\|x\|_{h^p(\mathscr{M})}, \tag{2.16}$$

*where the quantity in the middle identifies* $x\equiv x_\infty\in L^p(\mathscr{M})$*.*

In particular, bearing in mind the identification $x\equiv x_\infty$ for bounded martingales, see Remark 2.29, one has the following corollary.

**Corollary 2.32.** *For any* $p\in(1,\infty)$ *we have* $L^p(\mathscr{M})=\mathcal{H}^p(\mathscr{M})=h^p(\mathscr{M})$ *with equivalent norms.*

Stein's inequality [JX03] is a further corollary of Theorem 2.31.

**Corollary 2.33.** *Let* $(\mathscr{M},\omega,(\mathscr{M}_n)_n)$ *be a filtered modular space and let* $p\in(1,\infty)$*. On finite sequences* $x=(x_n)_n\subset L^p(\mathscr{M})$ *define the map* $Q(x)=(\omega_n^{(p)}(x_n))_n$*. Then, there exists a constant* $\theta_p>1$ *such that*

$$\|Q(x)\|_{L^p(\mathscr{M},\ell_\sharp^2)} \leqslant \theta_p\,\|x\|_{L^p(\mathscr{M},\ell_\sharp^2)}, \qquad \sharp=c,r$$

*that is* $Q$ *extends to a bounded projection on* $L^p(\mathscr{M},\ell_\sharp^2)$ *for* $\sharp=c,r$*.*

For constructing the Itô integral, some further spaces of adapted processes are needed. We denote by $S_{\mathrm{ad}}^p$ the linear space of all $L^p$ simple adapted processes and by $\mathbb{S}_{\mathrm{ad}}$ the simple adapted processes taking values in $\mathscr{M}_a$. Following [PX97], we introduce the following norms on $S_{\mathrm{ad}}^p$

$$\|F\|_{\mathcal{H}_c^p([0,t])} := \left\| \left( \int_0^t |F_s|^2 \mathrm{d}s \right)^{\frac{1}{2}} \right\|_p, \qquad \|F\|_{\mathcal{H}_r^p([0,t])} := \left\| \left( \int_0^t |F_s^*|^2 \mathrm{d}s \right)^{\frac{1}{2}} \right\|_p. \tag{2.17}$$

From these norm, one can introduce Hardy spaces as was done in the discrete setting, see Definition 2.30 and details in [PX97]. For applications to the twisted setting, we consider instead the following twisted Hardy spaces.

**Definition 2.34.** *For* $\sharp=c,r$ *let* $\mathbb{H}_\sharp^p([0,t])$ *be the completion of* $\mathbb{S}_{\mathrm{ad}}$ *with respect to the twisted norms*

$$\|F\|_{\mathbb{H}_\sharp^p([0,t])} := \sup_{|\tau|\leqslant 1-\frac{1}{2p}} \|T_\tau^{(p)}(F)\|_{\mathcal{H}_\sharp^p([0,t])}.$$

*For* $p\in[2,\infty)$ *set*

$$\mathbb{H}^p([0,t]) := \mathbb{H}_c^p([0,t]) \cap \mathbb{H}_r^p([0,t]), \qquad \|\cdot\|_{\mathbb{H}^p([0,t])} := \max_{\sharp=c,r} \|\cdot\|_{\mathbb{H}_\sharp^p([0,t])},$$

*and let* $\mathbb{H}_{\mathrm{loc}}^p(\mathbb{R}_+)$ *denote the space of processes on* $\mathbb{R}_+$ *whose restriction to* $[0,t]$ *is in* $\mathbb{H}^p([0,t])$*.*

In the rest of this section, we prove the density of $\mathbb{S}_{\mathrm{ad}}$ in $\mathbb{H}^p$. Whereas such a property holds by construction in $\mathbb{H}^p_\sharp$ for $\sharp = c, r$, this is not a priori true for $\mathbb{H}^p$ because of the intersection topology. As shown in [PX97] in the tracial setting, it is straightforward to project any element of $\mathbb{H}^p_{\mathrm{loc}}(\mathbb{R}_+)$ onto an approximating sequence in $\mathbb{S}_{\mathrm{ad}}$. To this end, we introduce the following bounded projections.

**Lemma 2.35.** *Let $\sigma = (t_j)_{j\geqslant 0}$ be a subdivision of $\mathbb{R}_+$ with $t_0 = 0$. On $\mathbb{S}_{\mathrm{ad}}$ let*

$$(Q_\sigma x)(t) := \frac{1}{t_{k+1} - t_k} \int_{t_k}^{t_{k+1}} \omega_k(x(s))\mathrm{d}s, \qquad t_k \leqslant t < t_{k+1}. \tag{2.18}$$

*Then, for any $p \in [2, \infty)$, $t \geqslant 0$ and $\sharp = c, r$ $Q_\sigma$ extends to a bounded projection on $\mathbb{H}^p_\sharp([0, t])$.*

**Proof.** The proof is an application of Stein's inequality in Theorem 2.33, see [PX97] for details, and the operator inequality for a family of bounded operators $(a_j)$ and the finite sequence $(r_j) \subset \mathbb{R}_+$, $\sum_j r_j = 1$

$$\Big|\sum_j r_j a_j\Big|^2 \leqslant \sum_j r_j |a_j|^2, \tag{2.19}$$

which is a consequence of convexity. The proof of [PX97] for the tracial case carries over: for fixed $\tau$ and $x \in \mathbb{S}_{\mathrm{ad}}$ we have

$$\begin{aligned}
\|T_\tau^{(p)}(Q_\sigma x)\|_{\mathcal{H}_c^p([0,t])} &= \Big\| \Big(\int_0^t |T_\tau^{(p)}(Q_\sigma x)|^2 \mathrm{d}s\Big)^{\frac{1}{2}} \Big\|_p \\
&\leqslant \Big\| \Big( \sum_k \sum_{t_k \leqslant s_j < t_{k+1}} (s_{j+1} - s_j) |T_\tau^{(p)}(\omega_k(x(s_j)))|^2 \Big)^{\frac{1}{2}} \Big\|_p \\
&\leqslant \theta_p \Big\| \Big( \sum_k \sum_{t_k \leqslant s_j < t_{k+1}} (s_{j+1} - s_j) |T_\tau^{(p)}(x(s_j))|^2 \Big)^{\frac{1}{2}} \Big\|_p \\
&= \theta_p \|T_\tau^{(p)}(x)\|_{\mathcal{H}_c^p([0,t])},
\end{aligned}$$

where the second line follows by (2.19) with $r_j = (s_{j+1} - s_j) / (t_{k+1} - t_k)$ and the definition of $Q_\sigma$, whereas the third line by Stein's inequality. □

One can also prove that, for any $x \in \mathbb{H}^p_{\sharp,\mathrm{loc}}(\mathbb{R}_+)$, by choosing $\sigma \to 0$, in the sense that $\sup_j (t_{j+1} - t_j) \to 0$,

$$\lim_{\sigma \to 0} Q_\sigma x = x,$$

in $\mathbb{H}^p_\sharp([0, t])$, for $\sharp = c, r$ and for any $t \geqslant 0$. In other words, one can produce an approximating sequence in $\mathbb{S}_{\mathrm{ad}}$ starting from $\mathbb{H}^p_\sharp([0, t])$. If $x \in \mathbb{H}^p([0, t])$ the approximating sequence converges in both $\{\mathbb{H}^p_\sharp([0, t])\}_{\sharp = c,r}$ and thus in the $\mathbb{H}^p([0, t])$ topology. We emphasize this result in the following corollary.

**Corollary 2.36.** *For any $p \in [2, \infty)$ and $t \geqslant 0$ $\mathbb{S}_{\mathrm{ad}}([0, t])$ is dense in $\mathbb{H}^p([0, t])$.*

# 3 Grassmann random variables

## 3.1 Araki–Wyss factors

While the tracial setting is not suitable for hosting non-trivial Grassmann martingales, the Araki–Wyss factors [AW64, Hia20] represent a convenient modular setting that serve our purposes. Inspired by the isomorphism of the $q$-deformed Araki–Woods factors with the $q$-deformed Baby Fock models shown in [Nou06], we introduce the Araki–Wyss factors as follows.

**Definition 3.1.** *(Araki–Wyss factors). Let $\mathcal{H}$ be a complex separable Hilbert space with conjugation $\Theta$. Let $0<\varrho<1$ be an operator acting on $\mathcal{H}$ such that $\Theta\,\varrho\,\Theta=\varrho$. Let $\Gamma_a(\mathcal{H}\oplus\mathcal{H})$ be the fermionic Fock space and let $a^*$, $a$ denote the creation and annihilation operators. Set*

$$\gamma_\varrho(f):=a(\varrho f\oplus 0)+a^*(0\oplus\varrho^{-1}\Theta f),\qquad f\in\mathcal{H},$$

*and*

$$\mathscr{A}(\mathcal{H},\varrho):=\{\gamma_\varrho(f)|f\in\mathcal{H}\}''.$$

**Remark 3.2.** One could more generally introduce the $q$-deformed Araki–Woods factors by considering $q$-Fock spaces instead, but this is beyond the scope of this paper.

To be precise, $\mathscr{A}(\mathcal{H},\varrho)$ is isomorphic to the Baby Fock model with $I=\mathbb{Z}_0$, and with $(\mu_j^{-1})_{j\in\mathbb{N}}$ being the spectrum of $\varrho$, see [Nou06, LR11].

We could equivalently work with the representation introduced by Shlyakhtenko [Shl97], see also [Hia03]. However, considering the Araki–Wyss factors in the present form makes the isomorphism with a baby Fock model manifest, so that the result of [LR11] on hypercontractivity carries over more easily, see Subsection 3.3. The following isomorphism holds.

**Proposition 3.3.** *In the same setting of Definition 3.1, $\mathscr{A}(\mathcal{H},\varrho)\cong\Gamma_{-1}(\mathcal{H}_{\mathbb{R}},(U_t^{(\varrho)})_t)$, where $\Gamma_q$ is the q-deformed Araki–Woods factor, where $\mathcal{H}_{\mathbb{R}}$ is the real Hilbert space such that $\mathcal{H}$ is complexification and where $U_t^{(\varrho)}$ is the following orthogonal transformation*

$$U_t^{(\varrho)}:=\begin{pmatrix}\cos(t\log\varrho^4)\,\mathbf{1}_{\mathcal{H}_{\mathbb{R}}} & -\sin(t\log\varrho^4)\,\mathbf{1}_{\mathcal{H}_{\mathbb{R}}}\\ \sin(t\log\varrho^4)\,\mathbf{1}_{\mathcal{H}_{\mathbb{R}}} & \cos(t\log\varrho^4)\,\mathbf{1}_{\mathcal{H}_{\mathbb{R}}}\end{pmatrix}. \tag{3.1}$$

**Remark 3.4.** Note that orthogonal transformations of the type in (3.1) are already considered in [Hia03, Nou06].

The operators $\gamma_\varrho$'s satisfy the following anticommutation relations

$$\{\gamma_\varrho(f),\gamma_\varrho(g)\}=0,\qquad\{\gamma_\varrho^*(f),\gamma_\varrho(g)\}=\langle g,(\varrho^2+\varrho^{-2})f\rangle_{\mathcal{H}},\qquad\forall f,g\in\mathcal{H} \tag{3.2}$$

One can check that the Fock vacuum $\omega(\cdot):=\langle\Omega,\cdot\,\Omega\rangle$ is a gauge-invariant quasi-free state and that

$$\omega(\gamma_\varrho^*(f)\,\gamma_\varrho(g))=\langle g,\varrho^{-2}f\rangle_{\mathcal{H}},\qquad\omega(\gamma_\varrho(g)\,\gamma_\varrho^*(f))=\langle g,\varrho^2 f\rangle_{\mathcal{H}}. \tag{3.3}$$

For convenience, we set $\mathcal{K}:=\mathcal{H}\oplus\mathcal{H}$, let $\hat{\Theta}$ be the conjugation on $\mathcal{K}$ such that $\hat{\Theta}f\oplus g=\Theta g\oplus\Theta f$ and define the linear operators $\beta_T:\mathcal{K}\to\mathscr{A}(\mathcal{H},\varrho)$

$$\begin{aligned}\beta_\varrho(f\oplus f') &:= \gamma_\varrho^*(\varrho^{-1}f)+\gamma_\varrho(\Theta\,\varrho f')\\ &= a^*(f\oplus f')+a((\varrho^2\oplus\varrho^{-2})\,\hat{\Theta}\,(f\oplus f')),\end{aligned} \tag{3.4}$$

which likewise generate $\mathscr{A}(\mathcal{H},\varrho)$, that is, $\mathscr{A}(\mathcal{H},\varrho)=\{\beta_\varrho(f)|f\in\mathcal{K}\}''$. For any $f,g\in\mathcal{K}$ we have

$$\omega(\beta_\varrho(f)\beta_\varrho(g))=\langle\hat{\Theta}f,(\varrho^2\oplus\varrho^{-2})\,g\rangle_{\mathcal{K}},\qquad\omega(\beta_\varrho^*(f)\beta_\varrho(g))=\langle f,g\rangle_{\mathcal{K}}. \tag{3.5}$$

We then introduce Wick's polynomials.

**Definition 3.5.** *Wick's Polynomials in the $\beta_\varrho$'s are defined recursively by setting $[\![\beta_\varrho(f)]\!]:=\beta_\varrho(f)$ for any $f\in\mathcal{K}$ and for $f,f_1,\ldots f_n\in\mathcal{K}$*

$$\begin{aligned}&[\![\beta_\varrho(f)\,\beta_\varrho(f_1)\cdots\beta_\varrho(f_n)]\!]:=\\ &\qquad\beta_\varrho(f)\,[\![\beta_\varrho(f_1)\cdots\beta_\varrho(f_n)]\!]+\sum_{j=1}^{b}(-)^j\,\omega(\beta_\varrho(f)\,\beta_\varrho(f_j))\,\big[\!\big[\beta_\varrho(f_1)\cdots\not\beta_\varrho(f_j)\cdots\beta_\varrho(f_n)\big]\!\big].\end{aligned}$$

Furthermore, for any $F = f_1 \wedge \cdots \wedge f_n \in \mathcal{K}^{\wedge n}$, $n \in \mathbb{N}$, we set $\beta_\varrho(F) := \beta_\varrho(f_1) \cdots \beta_\varrho(f_n)$ together with $a^*(F) := a^*(f_1) \cdots a^*(f_n)$ and extend this by linearity to polynomials in $\Gamma_a(\mathcal{K})$. By straightforward computations, one can prove that for any $F$ polynomial in $\Gamma_a(\mathcal{K})$

$$[\![\beta_\varrho(F)]\!]\,\Omega = a^*(F)\,\Omega, \tag{3.6}$$

so that $\Omega$ is cyclic for $\mathscr{A}(\mathcal{H},\varrho)$. As a simple consequence of (3.6), we can show that $\Omega$ is modular for $\mathscr{A}(\mathcal{H},\varrho)$ and thus $(\mathscr{A}(\mathcal{H},\varrho),\omega)$ is in its GNS. To this end, one argues as follows. Consider the vNa

$$\tilde{\mathscr{A}}(\mathcal{H},\varrho) := \{a(\varrho^{-1} f \oplus 0) - a^*(0 \oplus \varrho\,\Theta f) | f \in \mathcal{H}\}''.$$

Then, $\tilde{\mathscr{A}}(\mathcal{H},\varrho)$ is in the commutant of $\mathscr{A}(\mathcal{H},\varrho)$, that is, $\tilde{\mathscr{A}}(\mathcal{H},\varrho) \subset \mathscr{A}(\mathcal{H},\varrho)'$. We have proved that $\Omega$ is cyclic for $\mathscr{A}(\mathcal{H},\varrho)$ and likewise one can prove that $\Omega$ is cyclic for $\tilde{\mathscr{A}}(\mathcal{H},\varrho)$. Accordingly, $\Omega$ is cyclic also for $\mathscr{A}(\mathcal{H},\varrho)'$ and thus separating for $\mathscr{A}(\mathcal{H},\varrho)$, that is, $\Omega$ is modular for the latter. From (3.6) one can also deduce the identity

$$\omega([\![\beta_\varrho(F)]\!]^*\,[\![\beta_\varrho(G)]\!]) = \langle F, G\rangle_{\Gamma_a(\mathcal{K})}. \tag{3.7}$$

Since $\Omega$ is modular, one can compute the modular automorphism group and obtain the explicit formulas, compare with [Hia03]:

$$\Delta = \Gamma(\varrho^{-4} \oplus \varrho^4), \qquad \sigma_t(\beta_\varrho(f)) = \beta_\varrho((\varrho^{-4\mathrm{i}t} \oplus \varrho^{4\mathrm{i}t})f), \qquad \forall f \in \mathcal{K}, \tag{3.8}$$

where for any bounded operator $B: \mathcal{K} \to \mathcal{K}$, $\Gamma(B)$ denotes the operator such that $\Gamma(B) f_1 \otimes \cdots \otimes f_n = B f_1 \otimes \cdots \otimes B f_n$. The simplest way to compute $\Delta$ is by inspection of (3.5): in fact since $\omega$ is modular we have

$$\langle(\varrho^{-2} \oplus \varrho^2)\,\hat{\Theta} g, f\rangle_{\mathcal{K}} = \omega(\beta_\varrho(f)\,\beta_\varrho(g)) = \omega(\beta_\varrho(g)\,\sigma_{-\mathrm{i}}(\beta_\varrho(f))) = \langle(\varrho^2 \oplus \varrho^{-2})\hat{\Theta}\,g, f'\rangle_{\mathcal{K}},$$

so that $f' = (\varrho^{-4} \oplus \varrho^4) f$ and by the cyclic property obtain the expression for $\sigma_t(\beta_\varrho(f))$.

We conclude this section with the following classification result, which is a corollary of [Hia03, Theorem 3.1].

**Corollary 3.6.** *Let $G$ be the closed multiplicative subgroup of $\mathbb{R}_+$ generated by the spectrum of $\varrho$. Then, $\mathcal{A}(\mathcal{H},\varrho)$ is a factor of type*

$$\begin{cases} \mathrm{III}_\lambda & \text{if } G = \{\lambda^n | n \in \mathbb{Z}, \lambda \in (0,1)\}, \\ \mathrm{III}_1 & \text{if } G = \mathbb{R}_+. \end{cases}$$

## 3.2 Grassmann Brownian martingales

We begin with some general definitions about Grassmann algebras and Grassmann fields, compare with [ABDG22, DFG22]. First we recall the definition of Grassmann algebra.

**Definition 3.7.** *A Grassmann algebra is a unital complex algebra that has countably many generators $(\Psi_j)_j \subset \mathcal{A}$, which satisfy anticommutation relations*

$$\Psi_i\,\Psi_j + \Psi_j\,\Psi_i = 0, \qquad \forall i,j.$$

*If $\mathcal{A}, \mathcal{B} \subset \mathcal{C}$ are two Grassmann algebras which are subalgebras of some algebra $\mathcal{C}$, we say that $\mathcal{A}$ and $\mathcal{B}$ are compatible if* $\mathrm{span}(\mathcal{A},\mathcal{B})$*, i.e. the subalgebra generated by $\mathcal{A}, \mathcal{B}$ in $\mathcal{C}$, is a Grassmann algebra.*

**Remark 3.8.** Note that a Grassmann algebra is a *superalgebra*, that is, it can be split into two subspaces $\mathcal{A} = \mathcal{A}_+ \oplus \mathcal{A}_-$ such that for $\sigma, \sigma' \in \{\pm\}$ we have $\mathcal{A}_\sigma \mathcal{A}_{\sigma'} \subset \mathcal{A}_{\sigma\sigma'}$ and if $a \in \mathcal{A}_\sigma, a' \in \mathcal{A}_{\sigma'}$

$$a\,a' = (-1)^{\sigma\sigma'}\,a'\,a.$$

The two subspaces $\mathcal{A}_+$ and $\mathcal{A}_-$ will be referred to as even and odd elements respectively.

In our setting, we work in some filtered modular space $(\mathscr{M}, \omega, (\mathscr{M}_t)_t)$ and consider the Grassmann algebras generated by Grassmann fields indexed in a Hilbert space $\mathfrak{h}$, that is, $\Psi: \mathfrak{h} \to \mathscr{M}$ such that

$$\{\Psi(f), \Psi(g)\} = 0, \qquad \forall f, g \in \mathfrak{h}.$$

Therefore the $\mathrm{span}(\Psi(f), f \in \mathfrak{h})$, i.e. the algebra generated by the elements $\Psi(f)$, is a Grassmann subalgebra of $\mathscr{M}$ in the sense of Definition 3.7. As remarked in [ABDG22, DFG22], if $\Psi^{(1)}$ and $\Psi^{(2)}$ are two Grassmann fields, they need not be jointly Grassmann, that is,

$$\{\Psi^{(1)}(f), \Psi^{(2)}(g)\} = 0, \qquad \forall f, g \in \mathfrak{h}, \tag{3.9}$$

could fail to hold true. If (3.9) does hold we say that $\Psi^{(1)}$ and $\Psi^{(2)}$ are compatible, since in fact $\mathrm{span}(\Psi^{(1)}(f), f \in \mathfrak{h})$ and $\mathrm{span}(\Psi^{(2)}(f), f \in \mathfrak{h})$ are compatible Grassmann algebras.

A process $\Psi_t$ is said to be a Grassmann field or simply a Grassmann process if $\Psi_t: \mathfrak{h} \to \mathscr{M}$ and if satisfies anticommutation relations $\{\Psi_t(f), \Psi_s(g)\} = 0$ for any $s, t$ and any $\forall f, g \in \mathfrak{h}$. In a similar way, one defines Grassmann fields and processes valued in the $\mathbb{L}^p$ spaces.

It is also useful to introduce the notation of independent Grassmann algebras, or more generally independent subalgebras of $\mathscr{M}$.

**Definition 3.9.** *Let $(\mathscr{M}, \omega)$ be a modular space and let $\mathcal{A}_1, \mathcal{A}_2 \subset \mathscr{M}$ be two subalgebras. We say that $\mathcal{A}_1, \mathcal{A}_2$ are two independent algebras if for any $a_1 \in \mathcal{A}_1, a_2 \in \mathcal{A}_2$ we have*

$$\omega(a_1 a_2) = \omega(a_2 a_1) = \omega(a_1)\, \omega(a_2). \tag{3.10}$$

If $\mathcal{A}$ and $\mathcal{A}'$ are two Grassmann algebras, it is useful to be able to construct a new algebra containing $\mathcal{A}$ and $\mathcal{A}'$ in which $\mathcal{A}$ and $\mathcal{A}'$ which are two compatible and independent Grassmann algebras.

**Lemma 3.10.** *Let $(\mathscr{M}, \omega)$ and $(\mathscr{M}', \omega')$ be two modular spaces and let $\mathcal{A} \subset \mathscr{M}_a$ and $\mathcal{A}' = \mathscr{M}'_a$ be two Grassmann algebras which are subalgebras of the analytic elements in $\mathscr{M}$ and $\mathscr{M}'$ respectively. Suppose further that the canonical representations $C, C'$ of $(\mathscr{M}, \omega)$ and $(\mathscr{M}', \omega')$ respectively, admit some involutions $R, R'$, leaving some cyclic vectors invariant, anticommuting with $C(\Psi_j)$ and $C'(\Psi'_{j'})$, where $\Psi_j$ and $\Psi'_{j'}$ are the generators of the Grassmann algebras $\mathcal{A}$ and $\mathcal{A}'$ respectively.*

*Then, it is always possible to build a modular space $(\mathscr{M}_{\mathrm{tot}}, \omega_{\mathrm{tot}})$ such that $\mathscr{M} \hookrightarrow \mathscr{M}_{\mathrm{tot}}$ and $\mathscr{M}' \hookrightarrow \mathscr{M}_{\mathrm{tot}}$, and $\omega_{\mathrm{tot}}|_{\mathscr{M}} = \omega$ and $\omega_{\mathrm{tot}}|_{\mathscr{M}'} = \omega'$. Furthermore, if we identify $\mathcal{A}$ and $\mathcal{B}$ with their images in $\mathscr{M}_{\mathrm{tot}}$, with respect to the embeddings mentioned above, then $\mathcal{A}$ and $\mathcal{B}$ are independent compatible Grassmann algebras with respect to $\omega_{\mathrm{tot}}$ in the sense of Definition 3.9*

**Proof.** The proof can be found in [ABDG22, Lemma 5] in the case of $C^*$-algebras. The case of modular spaces, discussed in the lemma, is a straightforward generalisation. □

**Remark 3.11.** The hypotheses of Lemma 3.10 are satisfied by all finite-dimensional (or more generally countably-generated) Grassmann algebras and algebras generated by Gaussian random fields. These two examples of Grassmann algebras cover all the cases of Grassmann algebras considered in the rest of the paper.

In the following we will consider the notion of *independent increments process*.

**Definition 3.12.** *(Independent increment process) Let $X: \mathbb{R}_+ \times \mathfrak{h} \to \mathscr{M}$ be a stochastic process on $(\mathscr{M}, \omega, (\mathscr{M}_t)_t)$ indexed in $\mathfrak{h}$. We say that $X$ has independent increments if for any $s \leqslant t \in \mathbb{R}_+$, for any $k \in \mathbb{N}$ and for any $f_1, \ldots, f_k \in \mathfrak{h}$ we have*

$$\omega_s((X_t(f_1) - X_s(f_1)) \cdots (X_t(f_k) - X_s(f_k))) = \omega((X_t(f_1) - X_s(f_1)) \cdots (X_t(f_k) - X_s(f_k))).$$

**Remark 3.13.** A consequence of Definition 3.12 and the tower property for the conditional expectation, is that if the process $X_t$ has independent increments then, for any $s \leqslant t$, the algebra $\mathscr{M}_s$ and the algebra $\operatorname{span}\{(X_t(f_1) - X_s(f_1)) \cdots (X_t(f_k) - X_s(f_k)),\ k \in \mathbb{N}$ and $f_1, \ldots, f_k \in \mathfrak{h}\}$ are independent in the sense of Definition 3.9.

An important class of Grassmann fields are Grassmann Brownian martingales (GBM) [ABDG22, DFG22].

**Definition 3.14.** *Let $(G_t)_t$ a family of bounded and $\Theta$-antisymmetric operators on $\mathfrak{h}$. A Grassmann Brownian martingale on $(\mathscr{M}, \omega, (\mathscr{M}_t)_t)$ indexed in $\mathfrak{h}$ and with covariance $(G_t)_t$ is the adapted centered Gaussian Grassmann process $X = (X_t)_t$ with independent increments, that is, $\omega_s((X_t - X_s)(f)) = 0$ for any $s \leqslant t$ and $f \in \mathfrak{h}$, and such that*

$$\omega(X_t(f) X_s(g)) = \langle \Theta f, G_{t \wedge s} g \rangle_{\mathfrak{h}} \qquad \forall s, t, \forall f, g \in \mathfrak{h}. \tag{3.11}$$

We shall now show how non-trivial GBMs are constructed within Araki–Wyss factors. This is basically the Osterwalder–Schrader construction [Ost73] already implemented in [ABDG22, DFG22], although in such cases a Fock state was used instead.

Let $h$ be a complex separable Hilbert space, and let $\Theta = \begin{pmatrix} 0 & \mathbb{1} \\ \mathbb{1} & 0 \end{pmatrix} \kappa$ be a conjugation on $\mathfrak{h} := h \oplus h$, $\kappa$ denoting complex conjugation. We set

$$\mathcal{H}_t := L^2([0, t]; \mathfrak{h}),$$

with the identification $\mathcal{H} \equiv \mathcal{H}_\infty$ and for $0 < \mu < 1$ we introduce the vNa's

$$\mathscr{M}_t^{(\mu)} := \mathscr{A}(\mathcal{H}_t, \mu \mathbb{1}_{\mathcal{H}_t}) \subset \mathscr{M}^{(\mu)} := \mathscr{A}(\mathcal{H}, \mu \mathbb{1}_{\mathcal{H}}).$$

One can easily see that $(\mathscr{M}_t^{(\mu)})_{t \geqslant 0}$ is a filtration of vN subalgebras of $\mathscr{M}^{(\mu)}$ invariant under the modular automorphism group, see (3.8), in the sense that $\sigma_s(\mathscr{M}_t^{(\mu)})$ for any $t \geqslant 0$ and $s \in \mathbb{R}$. Thus, there exists the conditional expectation $\omega_t \colon \mathscr{M}^{(\mu)} \to \mathscr{M}_t^{(\mu)}$, so that $(\mathscr{M}^{(\mu)}, \omega, (\mathscr{M}_t^{(\mu)})_t)$ is a filtered modular space, see Definition 2.25.

**Lemma 3.15.** *Let $(G_t)_t$ be a differentiable family of bounded $\Theta$-symmetric operators on $\mathfrak{h}$. Let $\dot{G}_t = C_t^2 U_t$ be the polar decomposition of $\dot{G}_t$, with $U_t$ unitary and let $0 < \mu < 1$. Then, for any $t \geqslant 0$ the linear operators $\mathscr{C}_t, \tilde{\mathscr{C}}_t \colon \mathfrak{h} \to \mathcal{H}_t$ defined by*

$$\mathscr{C}_t f := (\mu^{-2} - \mu^2)^{-\frac{1}{2}} \int_0^t \delta_s \otimes C_s U_s f \, \mathrm{d}s, \qquad \tilde{\mathscr{C}}_t f := (\mu^{-2} - \mu^2)^{-\frac{1}{2}} \int_0^t \delta_s \otimes C_s \Theta f \, \mathrm{d}s,$$

*are bounded,*

$$\|\mathscr{C}_t f\|_{\mathcal{H}}^2 \lesssim_\mu \int_0^t \|C_s U_s f\|_{\mathfrak{h}}^2 \mathrm{d}s, \qquad \|\tilde{\mathscr{C}}_t f\|_{\mathcal{H}}^2 \lesssim_\mu \int_0^t \|C_s \Theta f\|_{\mathfrak{h}}^2 \mathrm{d}s,$$

*and furthermore*

$$\langle \tilde{\mathscr{C}}_s g, \mathscr{C}_t f \rangle_{\mathcal{H}} = -\langle \tilde{\mathscr{C}}_t f, \mathscr{C}_s g \rangle_{\mathcal{H}} = (\mu^2 - \mu^{-2})^{-1} \langle \Theta g, G_{s \wedge t} f \rangle_{\mathfrak{h}}.$$

We abridge our notation to $\gamma \equiv \gamma_{\varrho = \mu \mathbb{1}_{\mathcal{H}}}$ and introduce GBMs as follows.

**Definition 3.16.** *Let $(\mathscr{M}^{(\mu)}, \omega, (\mathscr{M}_t^{(\mu)})_{t \geqslant 0})$ be as above and, let $(G_t)_t$ be a differentiable family of bounded $\Theta$-symmetric operators on $\mathfrak{h}$. Letting $\mathscr{C}_t$ and $\tilde{\mathscr{C}}_t$ be as in Lemma 3.15, we define the process $X_t \colon \mathfrak{h} \to \mathscr{M}_t^{(\mu)}$*

$$X_t(f) := \gamma^*(\mathscr{C}_t f) - \gamma(\tilde{\mathscr{C}}_t f). \tag{3.12}$$

**Proposition 3.17.** *The process $X_t$ of Definition 3.16 is a GBM valued in $\mathscr{M}_a^{(\mu)} \cap \mathscr{M}_t^{(\mu)}$, with covariance $G_t$ and satisfies*

$$\|X_{s,t}(f)\|^2 \lesssim_\mu \int_s^t (\|C_r U_r f\|_{\mathfrak{h}}^2 + \|C_r \Theta f\|_{\mathfrak{h}}^2) \mathrm{d}r. \tag{3.13}$$

**Proof.** That $X$ is adapted is clear by inspection. We check anticommutation relations by (3.2) and by Lemma 3.15

$$\begin{aligned}\{X_t(f), X_s(g)\} &= -(\{\gamma^*(\mathscr{C}_t f), \gamma(\tilde{\mathscr{C}}_s g)\} + \{\gamma(\tilde{\mathscr{C}}_t f), \gamma^*(\mathscr{C}_s g)\})\\ &= -(\mu^2+\mu^{-2})\,(\langle\tilde{\mathscr{C}}_s g, \mathscr{C}_t f\rangle_{\mathcal{H}} + \langle\tilde{\mathscr{C}}_t f, \mathscr{C}_s g\rangle_{\mathcal{H}}) = 0.\end{aligned}$$

Furthermore, $X_t$ are Gaussian random variables because $\omega$ is quasi-free on $\gamma$'s. We compute the covariance with the aid of (3.3) and Lemma 3.15

$$\begin{aligned}\omega(X_t(f)\,X_s(g)) &= -(\omega(\gamma^*(\mathscr{C}_t f)\,\gamma(\tilde{\mathscr{C}}_s g) + \omega(\gamma(\tilde{\mathscr{C}}_t f)\,\gamma^*(\mathscr{C}_s g)))\\ &= -\mu^{-2}\langle\tilde{\mathscr{C}}_s g, \mathscr{C}_t f\rangle_{\mathcal{H}} - \mu^2\langle\tilde{\mathscr{C}}_t f, \mathscr{C}_s g\rangle_{\mathcal{H}}\\ &= \langle\Theta g, G_{s\wedge t} f\rangle_{\mathfrak{h}}\end{aligned}$$

Finally, (3.13) is obtained by Lemma 3.15 and by recalling that $\|a(f)\|, \|a^*(f)\| \leqslant \|f\|_{\mathcal{H}}$, for any $f \in \mathcal{H}$. □

An important special case of GBM is the Grassmann Brownian motion, corresponding to the choice $G_t = tC^2 U$ with $U = \mathbb{1} \oplus -\mathbb{1}$ and some bounded $C \geqslant 0$. In Section 4, we will only consider this special case, while more general GBMs will be used only in Section 5.1. In the former case, we note the following.

**Lemma 3.18.** *If $X_t$ is as in Definition 3.16, then $[X_t(f)]_\tau = \mu^{-4\tau} X_t(f)$. Assume furthermore that $G_t = tC^2 U$ with $[C, \Theta] = [C, U] = 0$. Then, there exist constants $c_{r,r'}$ and $c'_{r,r'}$ depending on $\mu$ and locally bounded in $r, r' \in \mathbb{R}$ such that*

$$\begin{aligned}\omega([X_t(f)]_r[X_s(g)]_{r'}) &= c_{r,r'}(s\wedge t)\langle\Theta f, G g\rangle_{\mathfrak{h}},\\ \omega([X_t^*(f)]_r[X_s(g)]_{r'}) &= c'_{r,r'}(s\wedge t)\langle f, C^2 g\rangle_{\mathfrak{h}}.\end{aligned} \tag{3.14}$$

*Moreover, setting $X_{s,t} := X_t - X_s$ we have $\|X_{s,t}(f)\|_{\mathbb{L}^p} \lesssim_{\mu,p} |t-s|^{\frac{1}{2}} \|Cf\|_{\mathfrak{h}}$ for any $p \in [1, \infty]$.*

**Proof.** The action of the automorphism group on $\beta \equiv \beta_{\varrho = \mu\mathbb{1}_{\mathcal{H}}}$, see (3.4) and (3.8) is simply $\sigma_t(\beta(f)) = \mu^{-4it}\beta(f)$, hence its analytical continuation is $[\beta(f)]_\tau = \mu^{-4\tau}\beta(f)$. We can write $X_t(f) = \beta(f_t \oplus \tilde{f}_t)$, see Lemma 3.15, with $f_t := \mu\,\mathscr{C}_t f$ and $\tilde{f}_t := \mu^{-1}\tilde{\mathscr{C}}_t f$ so that $[X_t(f)]_\tau = \mu^{-4\tau} X_t(f)$. This fact together with Proposition 3.17 implies the first identity in (3.14). The second identity follows by noting that $\mathscr{C}_t f \propto_\mu \mathbf{1}_{[0,t]} \otimes UCf$ and $\tilde{\mathscr{C}}_t f \propto_\mu \mathbf{1}_{[0,t]} \otimes \Theta\, Cf$, so that

$$\omega(X_t^*(f)\,X_s(g)) = \mu^{-2}\langle\tilde{\mathscr{C}}_s g, \tilde{\mathscr{C}}_t f\rangle_{\mathcal{H}} + \mu^2\langle\mathscr{C}_t f, \mathscr{C}_s g\rangle_{\mathcal{H}} = \frac{\mu^2+\mu^{-2}}{\mu^{-2}-\mu^2}\,(s\wedge t)\,\langle f, C^2 g\rangle_{\mathfrak{h}}.$$

Regarding the bound, by Hölder's inequality and by (3.13) we have

$$\|T_\tau^{(p)}(X_{s,t}(f))\|_p \leqslant \|[X_{s,t}(f)]_\tau\| \lesssim_{\mu,\tau} |t-s|^{\frac{1}{2}} \|Cf\|_{\mathfrak{h}},$$

implying the claim. □

On the other hand, the following identity will be useful in Section 5.

**Lemma 3.19.** *Let $X_t$ be as in Definition 3.16, then*

$$\begin{aligned}&\langle T_\tau^{(2)}([\![X_t(f_1)\cdots X_t(f_n)]\!]), T_\tau^{(2)}([\![X_s(g_1)\cdots X_s(g_n)]\!])\rangle_{L^2}\\ &\quad= \mu^{-2(1+4\tau)n}\langle(f_{1,t}\oplus\tilde{f}_{1,t})\wedge\cdots\wedge(f_{n,t}\oplus\tilde{f}_{n,t}), (g_{1,s}\oplus\tilde{g}_{1,s})\wedge\cdots\wedge(g_{n,s}\oplus\tilde{g}_{n,s})\rangle_{\Gamma_a(\mathcal{K})},\end{aligned} \tag{3.15}$$

*where $f_{j,t} := \mu\,\mathscr{C}_t f_j$ and $\tilde{f}_{j,t} := \mu^{-1}\tilde{\mathscr{C}}_t f_j$ and where $g_{j,s} := \mu\,\mathscr{C}_s g_j$ and $\tilde{g}_{j,s} := \mu^{-1}\tilde{\mathscr{C}}_s g_j$.*

**Proof.** We note that

$$\begin{aligned} T_\tau^{(2)}([\![X_t(f_1)\cdots X_t(f_n)]\!]) &= [[\![\beta(f_{1,t}\oplus \tilde{f}_{1,t})\cdots\beta(f_{n,t}\oplus\tilde{f}_{n,t})]\!]]_{\frac{1}{4}+\tau} D^{\frac{1}{2}} \\ &= \mu^{-(1+4\tau)n}[\![\beta(f_{1,t}\oplus \tilde{f}_{1,t})\cdots\beta(f_{n,t}\oplus\tilde{f}_{n,t})]\!] D^{\frac{1}{2}} \end{aligned}$$

and thus, by the definition of Haagerup's trace and by (3.7) we obtain (3.15). □

## 3.3 Hypercontractivity

Let us now recall the hypercontractivity result in [LR11] for the special case of Araki–Wyss factors. Let $P_t$ be the Ornstein–Uhlenbeck (OU) semigroup on $\mathscr{A} \equiv \mathscr{A}(\mathcal{H},\varrho)$, defined by

$$P_t(x)\,\Omega := \mathrm{e}^{-t\mathcal{N}}(x\,\Omega), \qquad \forall x\in\mathscr{A}(\mathcal{H},\varrho),$$

where $\mathcal{N} = \mathrm{d}\Gamma(\mathbb{1})$ is the number operator on $\Gamma_a(\mathcal{K})$. By (3.6), we have

$$P_t([\![\beta_\varrho(f_1)\cdots\beta_\varrho(f_n)]\!]) = \mathrm{e}^{-tn}[\![\beta_\varrho(f_1)\cdots\beta_\varrho(f_n)]\!], \tag{3.16}$$

that is, the Wick's polynomials are eigenvectors of the OU semigroup. Furthermore, note that the modular operator $\Delta$ commutes with the number operator $\mathcal{N}$ and thus

$$\sigma_s(P_t(x)) = P_t(\sigma_s(x)), \qquad \forall s,t\in\mathbb{R}, \tag{3.17}$$

so that $P_t$ maps the subalgebra of analytic elements $\mathscr{A}_a$ into itself. To extend this map to the $L^p$ spaces, one introduces the densely-defined operators, for $1\leqslant p,q\leqslant\infty$,

$$P_t^{(p,q)}: \mathscr{A} D^{\frac{1}{p}} \to \mathscr{A} D^{\frac{1}{q}}, \qquad x D^{\frac{1}{p}} \mapsto P_t(x) D^{\frac{1}{q}}. \tag{3.18}$$

These operators are not of the form (2.8), and in fact there is no general extension theorem that holds. Still, for suitable $p,q$, Lee–Ricard prove the following result [LR11].

**Theorem 3.20.** *(Lee–Ricard, 2011). Let $1<p<2$. Then $P_t^{(p,2)}$ extends to a contraction from $L^p(\mathscr{A}(\mathcal{H},\varrho))$ to $L^2(\mathscr{A}(\mathcal{H},\varrho))$, that is, $\|P_t^{(p,2)}\|_{L^p\to L^2}\leqslant 1$ provided that*

$$\mathrm{e}^{-2t} \lesssim \|\varrho\|_{\mathcal{H}\to\mathcal{H}}^{\frac{8}{p}-4}(p-1).$$

**Remark 3.21.** Note that the map $P_t^{(p,q)}$ can be more generally defined on $D^{\frac{1}{2p}+\tau}\mathscr{A} D^{\frac{1}{2p}-\tau}$ by

$$P_t^{(p,q)}(T_\tau^{(p)}(x)) := T_\tau^{(q)}(P_t(x)), \tag{3.19}$$

compare with (2.9), and would not depend on $\tau$ because of Lemma 2.12 and the state preserving property (2.7). The statement of the theorem holds for any such map as well since in fact $T_t^{(p)}(x) = [x]_{\frac{1}{2p}+t} D^{\frac{1}{p}}$ and we by (3.17) we have $[P_t(x)]_s = P_t([x]_s)$ for $x\in\mathscr{A}_a$, $P_t$ mapping $\mathscr{A}_a$ into itself.

**Corollary 3.22.** *Let $2<p<\infty$. Then $P_t^{(2,p)}$ extends to a contraction from $L^2$ to $L^p$ provided that $\mathrm{e}^{-2t}\lesssim \|\varrho\|_{\mathcal{H}\to\mathcal{H}}^{\frac{8}{p'}-4}(p'-1)$, with $1/p'+1/p=1$.*

**Proof.** By duality and by the density of $\mathscr{M}_a D^{\frac{1}{p}}$, we can write, for $x\in\mathscr{M}_a$ and $t$ as in the claim

$$\left\|P_t^{(2,p)}\left(x D^{\frac{1}{2}}\right)\right\|_p = \sup_{\substack{y\in\mathscr{M}_a D^{1/p'},\\ \|y\|_{p'}=1}} \left|\mathrm{tr}_{\mathrm{H}}\left(y^* P_t(x) D^{\frac{1}{p}}\right)\right| = \sup_{\substack{y\in\mathscr{M}_a D^{1/p'},\\ \|y\|_{p'}=1}} \left|\mathrm{tr}_{\mathrm{H}}\left(P_t^{(p',2)}(y)^* x D^{\frac{1}{2}}\right)\right| \leqslant \left\|x D^{\frac{1}{2}}\right\|_2,$$

where the bound follows by Hölder's inequality and by Theorem 3.20. This implies the claim. □

The following hypercontractive estimates hold true.

**Corollary 3.23.** *There exists a constant $C_{p,\varrho} \propto \|\varrho\|_{\mathcal{H}\to\mathcal{H}}^{2-\frac{4}{p}} (p-1)^{-\frac{1}{2}}$ such that, for any polynomial $F \in \Gamma_a(\mathcal{K})$*

$$\|T_\tau^{(2)}([\![\beta_\varrho(F)]\!])\|_2 \leqslant C_{p,\varrho}^{\deg(F)} \|T_\tau^{(p)}([\![\beta_\varrho(F)]\!])\|_p \qquad p \in [1,2), \tag{3.20}$$

$$\|T_\tau^{(p)}([\![\beta_\varrho(F)]\!])\|_p \leqslant C_{\frac{p}{p-1},\varrho}^{\deg(F)} \|T_\tau^{(2)}([\![\beta_\varrho(F)]\!])\|_2 \qquad p \in (2,\infty]. \tag{3.21}$$

**Proof.** Take a monomial $F = f_1 \wedge \cdots \wedge f_n$. Then, by (3.16) and (3.19)

$$T_\tau^{(2)}([\![\beta_\varrho(F)]\!]) = \mathrm{e}^{tn} T_\tau^{(2)}(P_t([\![\beta_\varrho(F)]\!])) = \mathrm{e}^{tn} P_t^{(p,2)}(T_\tau^{(p)}([\![\beta_\varrho(F)]\!])) \tag{3.22}$$

and the estimate for $\|T_\tau^{(2)}([\![\beta_\varrho(F)]\!])\|_2$ follows by Theorem 3.20, by choosing $\mathrm{e}^t \sim \|\varrho\|_{\mathcal{H}\to\mathcal{H}}^{2-\frac{4}{p}} (p-1)^{-\frac{1}{2}}$. The estimate for $\|T_\tau^{(p)}([\![\beta_\varrho(F)]\!])\|_p$ follows in the same way. To extend it to general polynomial $F = \sum_{n \geqslant 0} F_n$ with $F_n$ monomials, since $T_\tau^{(2)}(x) = [x]_{\frac{1}{4}+\tau} D^{\frac{1}{2}}$ by (3.7) we have

$$\|T_\tau^{(2)}([\![\beta_\varrho(F)]\!])\|_2^2 = \sum_{n,n'} \omega\Big([[\![\beta_\varrho(F_{n'})]\!]]_{\frac{1}{4}+\tau}^* [[\![\beta_\varrho(F_n)]\!]]_{\frac{1}{4}+\tau}\Big) = \sum_n \Big\|[[\![\beta_\varrho(F_n)]\!]]_{\frac{1}{4}+\tau}\Big\|_2^2, \tag{3.23}$$

which allows one to obtain (3.20) from (3.22). The bound (3.21) follows by writing

$$\begin{aligned}
\|T_\tau^{(p)}([\![\beta_T(F)]\!])\|_p^2 &= \left\|\sum_n \mathrm{e}^{tn} P^{(2,p)}(T_\tau^{(2)}([\![\beta_\varrho(F_n)]\!]))\right\|_p^2 \\
&\leqslant \left\|\sum_n \mathrm{e}^{tn} T_\tau^{(2)}([\![\beta_\varrho(F_n)]\!])\right\|_2^2 \\
&= \sum_n \mathrm{e}^{2tn} \|T_\tau^{(2)}([\![\beta_\varrho(F_n)]\!])\|_2^2 \\
&\leqslant \left(C_{\frac{p}{p-1},\varrho}^{\deg(F)}\right)^2 \|T_\tau^{(2)}([\![\beta_\varrho(F)]\!])\|_2^2,
\end{aligned}$$

where in the first and second inequality we used Theorem 3.20 and (3.23). □

# 4 Stochastic calculus

## 4.1 Itô integration

We extend the construction of the Itô–Clifford integral carried out in [BSW82, BSW83b, BSW83c, PX97] to the Grassmann setting. We begin by describing a general theory of integration with respect to $\mathbb{L}^p$ martingales in Subsections 4.1.1–4.1.2 under some mild assumptions on the martingale and then restrict to the special case of integration in the Grassmann algebra in Subsection 4.1.3. Unlike the case of Clifford random variables, we deal with martingales whose square is not simply a multiple of the identity, property which was crucially used in the former works in order to compute the quadratic variation of the Clifford martingale.

### 4.1.1 Scalar Itô integral

We begin by considering the integration of an $\mathbb{L}^p$ simple adapted process against an $\mathbb{L}^q$ martingale.

**Definition 4.1.** *(Itô integral). Let $q,p\in[1,\infty]$ such that $1/p+1/q\leqslant 1$. Let $F_t$ be a simple adapted process in $\mathbb{L}^p$ as defined in (2.11) and let $X_t$ be an $\mathbb{L}^q$ martingale, letting $X_{s,t}:=X_t-X_s$ for $s\leqslant t$, the Itô integral of $F$ with respect to $X$ is the process $Y=(Y_t)_t$*

$$Y_t=\int_0^t F_s\cdot \mathrm{d}X_s:=\sum_{j=0}^{k} F_{t_j}\cdot X_{t_j,t_{j+1}\wedge t}, \tag{4.1}$$

*where $k$ is such that $t\in[t_k,t_{k+1})$ and where the product between elements of the twisted spaces was given in Definition 2.17.*

**Remark 4.2.** The Itô integral defines an $\mathbb{L}^r$ martingale, with $1/r=1/p+1/q$.

In the rest of this section, we will extend the Itô integral to a larger class of adapted process. Note that in Definition 4.1 we could consider the integration of an $L^p$ simple adapted process against an $L^q$ martingale. However, one cannot extend the definition to more general $L^p$ adapted processes unless the martingale is analytic, stressing once more the relevance of the twisted $\mathbb{L}^p$ spaces.

As a warm up, we consider $\mathbb{L}^2$ valued processes and integration against a martingale in $\mathscr{M}_a$. We introduce the following subspace of $\mathbb{H}^p([0,t])$, see Definition 2.34, for $p\in[2,\infty)$

$$\tilde{\mathbb{H}}^p([0,t]):=\left\{F\in\mathbb{H}^p([0,t])\mid \int_0^t \|F_s\|_{\mathbb{L}^p}^2\mathrm{d}s<\infty\right\}. \tag{4.2}$$

That this is indeed a subspace follows from the fact that for any $F\in\mathbb{S}_{\mathrm{ad}}$

$$\|F\|^2_{\mathbb{H}^p_\sharp([0,t])}=\sup_{|\tau|\leqslant 1-\frac{1}{2p}}\|T_\tau^{(p)}(F)\|^2_{\mathcal{H}^p_\sharp([0,t])}=\sup_{|\tau|\leqslant 1-\frac{1}{2p}}\left\|\int_0^t |T_\tau^{(p)}(F_s)^\#|^2\mathrm{d}s\right\|_{\frac{p}{2}}\leqslant\int_0^t\|F_s\|^2_{\mathbb{L}^p}\mathrm{d}s,$$

where we used the shorthand $F^\#$ for $F$ or $F^*$ depending on whether $\sharp=c$ or $\sharp=r$.

**Proposition 4.3.** *Let $F$ be an $\mathbb{L}^2$ simple adapted process, let $X$ be a $\mathscr{M}_a$ martingale with independent increments such that, for any $0\leqslant s\leqslant t$*

$$\omega_s([X_{s,t}^*]_r[X_{s,t}]_{r'})=(t-s)\,c'_{r,r'},\qquad \forall 0\leqslant s\leqslant t,\qquad \forall r,r'\in\mathbb{R}, \tag{4.3}$$

*for some constant $c'_{r,r'}$. Let $Y_t$ denote the Itô integral as in Definition 4.1. Then, for any $|\tau|\leqslant\frac{3}{4}$ there is a constant $c_\tau$ (independent of $F$) such that*

$$\|T_\tau^{(2)}(Y_t)\|_2^2=c_\tau\int_0^t\|T_\tau^{(2)}(F_r)\|_2^2\mathrm{d}r. \tag{4.4}$$

*In particular, the Itô integral extends to a map from $\tilde{\mathbb{H}}^2([0,t])$ to $C^0([0,t],\mathbb{L}^2)$ and*

$$\|Y_t\|^2_{\mathbb{L}^2}\lesssim\int_0^t\|F_s\|^2_{\mathbb{L}^2}\mathrm{d}s.$$

**Remark 4.4.** Note that by Lemma 3.18 the GBM introduced in Definition 3.16 satisfies (4.3), having chosen some direction $f\in\mathfrak{h}$, or simply taking $\mathfrak{h}$ one dimensional.

**Proof.** We fix $t\in\mathbb{R}_+$ and the set $(t_j)_{j\geqslant 0}\subset\mathbb{R}_+$, we let $k$ be such that $t\in[t_k,t_{k+1})$ and set $\tilde{t}_j=t_j\wedge t$ for $j=0,\dots,k$. We also abridge our notation by introducing $F_j\equiv F_{\tilde{t}_j}$ and $\delta X_j\equiv X_{\tilde{t}_{j+1}}-X_{\tilde{t}_j}$. For $F_j\in\mathscr{M}_a\cap\mathscr{M}_j$, we have by the definition of Haagerup's trace

$$\|T_\tau^{(2)}(Y_t)\|_2^2=\mathrm{tr}_{\mathrm{H}}(T_{-\tau}^{(2)}(Y_t^*)\,T_\tau^{(2)}(Y_t))=\omega(Y_t^*\,[Y_t]_s),$$

where we set $s=\frac{1}{2}+2\tau$. Then, by the modular property we can write

$$\|T_\tau^{(2)}(Y_t)\|_2^2=\sum_{j,j'}\omega(\delta X_j^*F_j^*\,[F_{j'}\,\delta X_{j'}]_s)=\sum_{j,j'}\omega(F_j^*\,[F_{j'}]_s[\delta X_{j'}]_s\,[\delta X_j^*]_1).$$

If $j' > j$, then we have

$$\omega(F_j^*[F_{j'}]_s[\delta X_{j'}]_s[\delta X_j^*]_1) = \omega(F_j^*[F_{j'}]_s\,\omega_{j'}([\delta X_{j'}]_s)\,[\delta X_j^*]_1) = 0,$$

since the conditional expectation is state preserving, see (2.7), and since $X$ is an $\mathcal{M}_a$ martingale. Similarly, is the expectation vanishing if $j < j'$. Therefore, by assumption (4.3)

$$\|T_\tau^{(2)}(Y_t)\|_2^2 = \sum_j \omega(F_j^*[F_j]_s)\,\omega_j(\delta X_j^*[\delta X_j]_s) = c'_{0,s}\sum_j \|T_\tau^{(2)}(F_j)\|_2^2\,\delta t_j, \tag{4.5}$$

so that, by continuity, this holds true for any $\mathbb{L}^2$ valued simple adapted $F$, proving the claim (4.4). Since $c'_{0,s}$ is locally bounded, we can take the sup over $\tau$ we obtain $\|Y_t\|_{\mathbb{L}^2}^2 \lesssim \int_0^t \|F_s\|_{\mathbb{L}^2}^2\,\mathrm{d}s$ for any $F \in \mathbb{S}_{\mathrm{ad}}$. The conclusion follows by the density of $\mathbb{S}_{\mathrm{ad}}$ in $\tilde{\mathbb{H}}^2([0,t])$, see Lemma 2.36. □

To extend this result to the general case of integration of an $\mathbb{L}^p$ adapted process against a suitable $\mathbb{L}^q$ martingale, we shall now resort to martingales inequalities, see Theorem 2.31. The following is the main result of this section.

**Theorem 4.5.** *Let $p, q, r \in (2, \infty]$ with $1/r \geqslant 1/p + 1/q$. Let $X$ be an $\mathbb{L}^q$ martingale with independent increments and such that*

$$\|X_{s,t}\|_{\mathbb{L}^q} \lesssim (t-s)^{\frac{1}{2}}, \qquad \forall 0 \leqslant s \leqslant t. \tag{4.6}$$

*Then, the Itô integral extends to a map $F \mapsto Y_\cdot = \int_0^\cdot F_s \cdot \mathrm{d}X_s$ from $\mathbb{H}^p_{\mathrm{loc}}(\mathbb{R}_+)$ to $C^0(\mathbb{R}_+, \mathbb{L}^r)$ and*

$$\|Y_t\|_{\mathbb{L}^r} \lesssim_r \|F\|_{\mathbb{H}^p([0,t])}. \tag{4.7}$$

**Remark 4.6.** Note once more that (4.6) is satisfied by the GBM introduced in Definition 3.16, see Lemma 3.18.

**Proof.** For the sake of brevity, we denote by $\mathcal{M}_j \equiv \mathcal{M}_{t_j}$ and write $F_j := F_{t_j}$, $X_j := X_{t_j}$, $\delta X_j := X_{j+1} - X_j$ and $\delta Y_j = F_{j-1} \cdot \delta X_{j-1}$ for $j \geqslant 1$ and $\delta Y_0 := 0$. The martingale property follows by direct inspection, since $F_j$ is adapted and $\mathcal{E}_{j'}^{(q)}(\delta X_j) = 0$, $\forall j' \leqslant j$, being $X_t$ a martingale.

To prove (4.7), we control the norm $\|T_\tau^{(r)}(Y_t)\|_r$ at fixed $\tau$ by means of the Burkholder's inequality (2.16). Thus, we shall bound the norms $\|\cdot\|_{h_c^r}$, $\|\cdot\|_{h_r^r}$ and $\|\cdot\|_{h_d^r}$. First of all, we consider $F_j, X_j \in \mathcal{M}_a \cap \mathcal{M}_j$ and note that

$$|T_\tau^{(r)}(\delta Y_{j+1})|^2 = D^{\frac{1}{2r}-\tau}\,\delta X_j^* F_j^*[F_j \delta X_j]_{\frac{1}{r}+2\tau} D^{\frac{3}{2r}+2\tau}.$$

Thus, by definition of $\omega_j^{(\frac{r}{2})}$ and by the properties of $\omega_j$

$$\begin{aligned}
\omega_j^{(\frac{r}{2})}(|T_\tau^{(r)}(\delta Y_{j+1})|^2) &= D^{\frac{1}{2r}-\tau}\omega_j\Big(\delta X_j^* F_j^*[F_j\delta X_j]_{\frac{1}{r}+2\tau}\Big)D^{\frac{3}{2r}+\tau}\\
&= D^{\frac{1}{2r}-\tau}\omega_j\Big(F_j^*[F_j\delta X_j]_{\frac{1}{r}+2\tau}[\delta X_j^*]_1\Big)D^{\frac{3}{2r}+\tau}\\
&= |T_\tau^{(r)}(F_j)|^2\,\omega_j\Big(\delta X_j^*[\delta X_j]_{\frac{1}{r}+2\tau}\Big),
\end{aligned}$$

where we used that $F_j, F_j^* \in \mathcal{M}_j$ and that $\omega_j\Big(\delta X_j^*[\delta X_j]_{\frac{1}{r}+2\tau}\Big) \in \mathbb{C}$, $X$ having independent increments. The latter property also implies

$$\Big\|\omega_j^{(\frac{r}{2})}(|T_\tau^{(r)}(\delta X_j)|^2)\Big\|_{\frac{r}{2}} = \Big\|D^{\frac{1}{2r}-\tau}\,\omega_j\Big(\delta X_j^*[\delta X_j]_{\frac{1}{r}+2\tau}\Big)\,D^{\frac{3}{2r}+\tau}\Big\|_{\frac{r}{2}} = \Big|\omega_j\Big(\delta X_j^*[\delta X_j]_{\frac{1}{r}+2\tau}\Big)\Big|.$$

Accordingly, for any $|\tau| \leqslant 1 - \frac{1}{2r}$, by Proposition 2.22 and by Hölder's inequality we have, recalling that $r \leqslant q$

$$\Big|\omega_j\Big(\delta X_j^*[\delta X_j]_{\frac{1}{r}+2\tau}\Big)\Big| \leqslant \||T_\tau^{(r)}(\delta X_j)|\|_r^2 \leqslant \|\delta X_j\|_{\mathbb{L}^q}^2.$$

Since by Proposition 2.22 $\omega_j^{(\frac{r}{2})}(|T_\tau^{(r)}(\delta Y_{j+1})|^2) \geqslant 0$, we obtain

$$\omega_j^{(\frac{r}{2})}(|T_\tau^{(r)}(\delta Y_{j+1})|^2) \leqslant \sum_{j\geqslant 0} |T_\tau^{(r)}(F_j)|^2 \|\delta X_j\|_{\mathbb{L}^q}^2, \tag{4.8}$$

which can be extended by continuity to $X_j \in \mathbb{L}^q$. All in all, by assumption (4.6) we have

$$\sum_{j\geqslant 0} \omega_j^{(\frac{r}{2})}(|T_\tau^{(r)}(\delta Y_{j+1})|^2) \lesssim \sum_{j\geqslant 0} |T_\tau^{(r)}(F_j)|^2 \delta t_j = D^{\frac{1}{q}} \left(\sum_{j\geqslant 0} \left|T^{(p)}_{\tau+\frac{1}{2q}}(F_j)\right|^2 \delta t_j\right) D^{\frac{1}{q}}$$

We observe that

$$\left(D^{\frac{1}{q}} \left(\sum_{j\geqslant 0} \left|T^{(p)}_{\tau+\frac{1}{2q}}(F_j)\right|^2 \delta t_j\right) D^{\frac{1}{q}}\right)^{\frac{1}{2}} = \left|\left(\sum_{j\geqslant 0} \left|T^{(p)}_{\tau+\frac{1}{2q}}(F_j)\right|^2 \delta t_j\right)^{\frac{1}{2}} D^{\frac{1}{q}}\right|,$$

and thus, since $\||x|\|_r = \|x\|_r$, by Hölder's inequality we have

$$\begin{aligned} \|T_\tau^{(r)}(Y_t)\|_{h_c^r} &\lesssim \left\|\left(D^{\frac{1}{q}} \left(\sum_{j\geqslant 0} \left|T^{(p)}_{\tau+\frac{1}{2q}}(F_j)\right|^2 \delta t_j\right) D^{\frac{1}{q}}\right)^{\frac{1}{2}}\right\|_r = \left\|\left|\left(\sum_{j\geqslant 0} \left|T^{(p)}_{\tau+\frac{1}{2q}}(F_j)\right|^2 \delta t_j\right)^{\frac{1}{2}} D^{\frac{1}{q}}\right|\right\|_r \\ &\leqslant \left\|\left(\sum_{j\geqslant 0} \left|T^{(p)}_{\tau+\frac{1}{2q}}(F_j)\right|^2 \delta t_j\right)^{\frac{1}{2}}\right\|_p. \end{aligned} \tag{4.9}$$

We observe that if $|\tau| \leqslant 1-\frac{1}{2r}$, clearly $\left|\tau+\frac{1}{2q}\right| \leqslant 1-\frac{1}{2r}+\frac{1}{2q} = 1-\frac{1}{2p}$, hence

$$\sup_{|\tau|\leqslant 1-\frac{1}{2r}} \|T_\tau^{(r)}(Y_t)\|_{h_c^r} \lesssim \sup_{|\tau|\leqslant 1-\frac{1}{2p}} \left\|\left(\int_0^t |T_\tau^{(p)}(F_s)|^2 \mathrm{d}s\right)^{\frac{1}{2}}\right\|_p \equiv \|F\|_{\mathbb{H}_c^p([0,t])},$$

and similarly for the $\|\cdot\|_{h_r^r}$ norm. We control the diagonal Hardy norm by using Hölder's inequality for the twisted spaces and the monotonicity of $x \mapsto x^{\frac{1}{r}}$ and $x \mapsto x^r$ for $x \in \mathbb{R}_+$,

$$\sup_\tau \left(\sum_j \|\delta Y_j\|_r^r\right)^{\frac{1}{r}} \leqslant \left(\sum_j \left(\sup_\tau \|\delta Y_j\|_r\right)^r\right)^{\frac{1}{r}} = \left(\sum_j \|\delta Y_j\|_{\mathbb{L}^r}^r\right)^{\frac{1}{r}} \lesssim \left(\sum_j \|F_j\|_{\mathbb{L}^p}^r \delta t_j^{\frac{r}{2}}\right)^{\frac{1}{r}},$$

which can be made as small as desirable for $r>2$. Thus, the $h^p$ norm is dominated by the contribution $h_\sharp^p$ for $\sharp = c, r$, and the claim follows by the Burkholder's inequality

$$\|Y_t\|_{\mathbb{L}^r} \lesssim_p \max_{\sharp=c,r} \sup_{|\tau|\leqslant 1-\frac{1}{2r}} \|T_\tau^{(r)}(Y_t)\|_{h_\sharp^r} \lesssim \max_{\sharp=c,r} \|F\|_{\mathbb{H}_\sharp^p([0,t])} = \|F\|_{\mathbb{H}_{\mathrm{loc}}^p([0,t])}.$$

By Lemma 2.36, $\mathbb{S}_{\mathrm{ad}}$ is dense in $\mathbb{H}^p$ hence (4.9) extends by continuity to any $F \in \mathbb{H}^p([0,t])$. □

**Remark 4.7.** It is possible to modify the hypotheses of Theorem 4.5, in such a way to consider processes $X$ satisfying inequality (4.6) without having independent increments. Indeed, if we consider a martingale $X$ such that $\|X_{s,t}\|_{\mathbb{L}^\infty} \lesssim (t-s)^{\frac{1}{2}}$ and $F \in \tilde{\mathbb{H}}_{\mathrm{loc}}^p$, for some $p \in (2,\infty]$, see (4.2) then $Y_t = \int_0^t F_s \cdot \mathrm{d}X_s$ is well-defined and, for every $r \in (2,p]$, $\|Y_t\|_{\mathbb{L}^r} \lesssim \|F\|_{\tilde{\mathbb{H}}^p([0,t])}$. Under the described hypotheses, for any simple process $F \in \tilde{\mathbb{H}}_{\mathrm{loc}}^p$, we get

$$\left\|\sum_{j\geqslant 0} \omega_j^{(\frac{r}{2})}(|T_\tau^{(r)}(\delta Y_{j+1})|^2)\right\|_{\mathbb{L}^r} \lesssim \sum_{j\geqslant 0} \|F_j\|_{\mathbb{L}^p} \delta t_j,$$

where we used the fact that $\left\|\omega_j\left(\delta X_j^*[\delta X_j]_{\frac{1}{r}+2\tau}\right)\right\|_{L^\infty} \leqslant \|\delta X_j\|_{\mathbb{L}^\infty}^2 \lesssim \delta t_j$.

### 4.1.2 Multidimensional Itô integral

We shall now consider martingales and adapted processes indexed in a Hilbert space $\mathfrak{h}$ with conjugation $\Theta$. We begin with the following definition.

**Definition 4.8.** *Let $p \in [1, +\infty]$. We say that $F$ is a simple adapted process in $\mathbb{L}^p$ indexed in the Hilbert space $\mathfrak{h}$ if $F: \mathbb{R}_+ \times \mathfrak{h} \to \mathbb{L}^p$ such that for any $t \in \mathbb{R}_+$, the map $v \mapsto F_t(v)$, $v \in \mathfrak{h}$, is linear, and for any $v \in \mathfrak{h}$, the process $t \mapsto F_t(v)$ is simple and adapted. Furthermore, we say that $F$ is an $\mathbb{L}^p$ finite-dimensional simple adapted process, if furthermore there exists a set $\{v_1, \ldots, v_k\} \subset \mathfrak{h}$ such that $F_t(v) = 0$ for any $v \in \operatorname{span}\{v_1, \ldots, v_k\}^\perp$.*

For finite-dimensional simple adapted processes the multidimensional Itô integral is the following generalisation of Definition 4.1.

**Definition 4.9.** *Let $q, p \in [1, \infty]$ such that $1/p + 1/q \leqslant 1$. Let $F$ be an $\mathbb{L}^p$ finite-dimensional simple adapted process and $X_t$ an $\mathbb{L}^q$ martingale, both indexed in the Hilbert space $\mathfrak{h}$. The multidimensional Itô integral of $F$ with respect to $X$ is the process $Y = (Y_t)_t$*

$$Y_t = \int_0^t \langle F_s, \mathrm{d}X_s \rangle := \sum_\alpha \int_0^t F_s(v_\alpha) \cdot \mathrm{d}X_s(v_\alpha),$$

*where $\{v_\alpha\}_{\alpha \in \mathbb{N}}$ is any orthonormal real basis with respect the conjugation $\Theta$ and where the r.h.s. is the scalar Itô integral of Definition 4.1.*

**Remark 4.10.** Note that $Y$ does not depend on the choice of the orthonormal real basis $\{v_\alpha\}_{\alpha \in \mathbb{N}}$.

**Notation 4.11.** *Let $A: \mathfrak{h} \to \mathfrak{h}$ be a bounded linear operator and $F: \mathfrak{h} \to \mathbb{L}^p$, $H: \mathfrak{h} \to \mathbb{L}^q$, for some $q$, $p \geqslant 2$, be some linear (or anti-linear) maps. We set*

$$\operatorname{Tr}_A(F \otimes H) := \sum_{\alpha, \beta} \langle A v_\alpha, v_\beta \rangle_{\mathfrak{h}} F(v_\alpha) \cdot H(v_\beta)$$

*where $\{v_\alpha\}_{\alpha \in \mathbb{N}}$ is some orthonormal real basis of $\mathfrak{h}$ with respect to the conjugation $\Theta$, whenever the series is absolutely convergent in $\mathbb{L}^r$, $1/r = 1/p + 1/q$, in which case it does not depend on the choice of the real orthonormal basis $\{v_\alpha\}$.*

**Remark 4.12.** When $A \geqslant 0$ commutes with $\Theta$, there is $\tilde{A}$ such that $\tilde{A}^2 = A$ and for any orthonormal (real) basis $\{v_\alpha\}$ of $\mathfrak{h}$ we have $\langle \tilde{A} v_\alpha, v_\beta \rangle_{\mathfrak{h}} = \langle \tilde{A} v_\beta, v_\alpha \rangle_{\mathfrak{h}} \in \mathbb{R}$ and as a series of positive operators

$$\operatorname{Tr}_A(F^* \otimes F) = \sum_\alpha |F(\tilde{A}(v_\alpha))|^2. \tag{4.10}$$

For this reason it is appealing to introduce the notation

$$|F|_A^2 := \operatorname{Tr}_A(F^* \otimes F).$$

As long as $F$ is finite-dimensional, the Itô integral can be extended as was done in the previous section to processes such that $F(v_\alpha) \in \mathbb{H}^p([0, t])$ for $\alpha$ in a finite set. Otherwise, we introduce a finer topology in the following way. We let $\mathbb{S}_{\mathrm{ad}}(\mathfrak{h})$ denote the linear space of finite-dimensional simple adapted processes indexed in the Hilbert space $\mathfrak{h}$ and taking value in $\mathscr{M}$_a. On $\mathbb{S}_{\mathrm{ad}}(\mathfrak{h})$ we introduce the norms, compare with (2.17)

$$\|F\|_{\mathcal{H}_c^p([0,t]);A} := \left\| \left( \int_0^t |F_s|_A^2 \, \mathrm{d}s \right)^{\frac{1}{2}} \right\|_p, \qquad \|F\|_{\mathcal{H}_r^p([0,t]);A} := \left\| \left( \int_0^t |F_s^*|_A^2 \, \mathrm{d}s \right)^{\frac{1}{2}} \right\|_p. \tag{4.11}$$

We then introduce the following extension of the twisted Hardy spaces of Definition 2.34.

**Definition 4.13.** *Let $\mathfrak{h}$ be a separable Hilbert space with conjugation $\Theta$. Let $A \geqslant 0$ be a bounded linear operator on $\mathfrak{h}$ commuting with $\Theta$. For $\sharp = c, r$ let $\mathbb{H}^p_{A,\sharp}([0,t])$ be the completion of $\mathbb{S}_{\mathrm{ad}}(\mathfrak{h})$ with respect to the twisted norms*

$$\|F\|_{\mathbb{H}^p_{\sharp,A}([0,t])} := \sup_{|\tau| \leqslant 1-\frac{1}{2p}} \|T^{(p)}_\tau(F)\|_{\mathcal{H}^p_\sharp([0,t]);A}.$$

*For $p \in [2,\infty)$ set*

$$\mathbb{H}^p_A([0,t]) := \mathbb{H}^p_{A,c}([0,t]) \cap \mathbb{H}^p_{A;r}([0,t]), \qquad \|\cdot\|_{\mathbb{H}^p_A([0,t])} := \max_{\sharp=c,r} \|\cdot\|_{\mathbb{H}^p_{A,\sharp}([0,t])},$$

*and let $\mathbb{H}^p_{A,\mathrm{loc}}(\mathbb{R}_+)$ denote the space of processes on $\mathbb{R}_+$ whose restriction to $[0,t]$ is in $\mathbb{H}^p([0,t])$.*

**Remark 4.14.** It is important to note that $\mathbb{S}_{\mathrm{ad}}$ is dense also in $\mathbb{H}^p_A$ and this is a simple generalisation of Corollary 2.36. In fact, it suffices to observe that (2.19), which was used in the proof of the Stein's inequality, holds for $|\cdot|_A$ as well: by (4.10), we have

$$\left|\sum_j r_j F_j\right|^2_A = \sum_\alpha \left|\sum_j r_j F_j(\tilde{A} v_\alpha)\right|^2 \leqslant \sum_\alpha \sum_j r_j |F_j(\tilde{A} v_\alpha)|^2 = \sum_j r_j |F_j|^2_A.$$

We also introduce the following subspace of $\mathbb{H}^2_A([0,t])$, compare with (4.2)

$$\tilde{\mathbb{H}}^2_A([0,t]) := \left\{F \in \mathbb{H}^2_A([0,t]) | \int_0^t \||F_s|_A\|^2_{\mathbb{L}^2} \mathrm{d}s < \infty\right\}.$$

We consider integration of an $\mathbb{L}^p$ adapted process for $p \geqslant 2$ against a $\mathscr{M}_a$ martingale. In comparison with Proposition 4.3, note that we consider more general adapted processes, as was done in Theorem 4.5, albeit we still have some strong assumption on the martingale.

**Theorem 4.15.** *Let $\mathfrak{h}$ be a separable Hilbert space with conjugation $\Theta$. Let $X$ be an $\mathscr{M}_a$ martingale indexed in $\mathfrak{h}$ such that*

$$\omega([X^*_{s,t}(f)]_r [X_{s,t}(g)]_{r'}) = (t-s)\, c'_{r,r'} \langle f, A g\rangle_{\mathfrak{h}}, \qquad \forall 0 \leqslant s \leqslant t, \quad \forall f, g \in \mathfrak{h}, \quad \forall r, r' \in \mathbb{R}, \tag{4.12}$$

*for some constant $c'_{r,r'}$ and some $A \geqslant 0$ commuting with $\Theta$. Then, the Itô integral extends to a map $F \mapsto Y_\cdot = \int_0^\cdot \langle F_s, \mathrm{d}X_s\rangle$ from $\tilde{\mathbb{H}}^2_A([0,t])$ to $C^0([0,t], \mathbb{L}^2)$ and*

$$\|Y_t\|^2_{\mathbb{L}^2} \lesssim \int_0^t \||F_s|_A\|^2_{\mathbb{L}^2} \mathrm{d}s.$$

*Assume furthermore that $\|[X_{s,t}(f)]_r\| \lesssim_{r,f} |t-s|^{1/2}$. Then, the Itô integral extends to a map $F \mapsto Y_\cdot = \int_0^\cdot \langle F_s, \mathrm{d}X_s\rangle$ from $\mathbb{H}^p_{A;\mathrm{loc}}(\mathbb{R}_+)$ to $C^0(\mathbb{R}_+, \mathbb{L}^p)$ for $p \in (2,\infty]$ and*

$$\|Y_t\|_{\mathbb{L}^p} \lesssim_p \|F\|_{\mathbb{H}^p_A([0,t])}.$$

**Proof.** In the case of $p = 2$ the proof of Proposition 4.3 carries over and one obtains

$$\begin{aligned}
\|T^{(2)}_\tau(Y_t)\|^2_2 &= \sum_j \sum_{\alpha,\beta} \omega(F^*_j(v_\alpha)\,[F_j(v_\beta)]_s)\,\omega_j(\delta X^*_j(v_\alpha)\,[\delta X_j(v_\beta)]_s) \\
&= c'_{0,s} \sum_j \sum_{\alpha,\beta} \mathrm{tr}_{\mathrm{H}}(T^{(2)}_\tau(F_j(v_\alpha))^* T^{(2)}_\tau(F_j(v_\beta)))\,\delta t_j \langle v_\alpha, A v_\beta\rangle_{\mathfrak{h}} \\
&= c'_{0,s} \sum_j \||T^{(2)}_\tau(F_j)|^2_A\|\, \delta t_j,
\end{aligned}$$

in place of (4.5), hence the claim. The case $p \in (2, \infty]$ is an adaptation of the proof of Theorem 4.5. Following the said proof, we obtain the identity

$$\begin{aligned}\omega_j^{(\frac{p}{2})}(|T_\tau^{(p)}(\delta Y_{j+1})|^2) &= \sum_{\alpha,\beta} T_\tau^{(p)}(F_j(v_\alpha))^* T_\tau^{(p)}(F_j(v_\beta))\, \omega_j\Big(\delta X_j(v_\alpha)^* [\delta X_j(v_\beta)]_{\frac{1}{p}+2\tau}\Big) \\ &= c'_{0,\frac{1}{p}+2\tau} |T_\tau^{(p)}(F_j)|_A^2\, \delta t_j\end{aligned}$$

compare with (4.8). Thus, instead of (4.9) we obtain the estimate

$$\|T_\tau^{(p)}(Y_t)\|_{h_c^p} \lesssim \left\| \left( \sum_{j \geqslant 0} |T_\tau^{(p)}(F_j)|_A^2\, \delta t_j \right)^{\frac{1}{2}} \right\|_p,$$

and likewise for the $\|\cdot\|_{h_r^p}$ norm. To control the diagonal Hardy norm, we use the additional assumption $\|[X_{s,t}(f)]_r\| \lesssim_{r,f} |t-s|^{1/2}$ and obtain

$$\sup_\tau \left( \sum_j \|\delta Y_j\|_p^p \right)^{\frac{1}{p}} \lesssim \left( \sum_j \left( \sum_\alpha C_\alpha \|F_j(v_\alpha)\|_{\mathbb{L}^p} \right)^p \delta t_j^{\frac{p}{2}} \right)^{\frac{1}{p}},$$

which for $p > 2$ and finite-dimensional processes can be made arbitrary small by making the partition of $[0, t]$ finer. The conclusion holds by the density of $\mathbb{S}_{\mathrm{ad}}$ in $\mathbb{H}_A^p$, see Remark 4.14. □

**Remark 4.16.** If one requires weaker assumptions, e.g.

$$\|X_{s,t}(f)\|_{\mathbb{L}^q} \lesssim (t-s)^{\frac{1}{2}} \|A^{1/2} f\|_{\mathfrak{h}}, \qquad \forall 0 \leqslant s \leqslant t, \quad \forall f \in \mathfrak{h}, \tag{4.13}$$

then one cannot obtain the result of Theorem 4.15 without modifying the topology on $\mathbb{S}_{\mathrm{ad}}^p$ in a base dependent way. Indeed, if we suppose that $F \in \mathbb{H}_A^p$ is such that there is a basis $\{v_\alpha\}$ of $\mathfrak{h}$ for which, for any $t \in \mathbb{R}_+$ and for $F^\# = F, F^*$, we have

$$\left\| \left( \int_0^t \left| \sum_\alpha F_s^\#(v_\alpha)\, \|A^{1/2} v_\alpha\|_{\mathfrak{h}} \right|^2 \mathrm{d}s \right)^{\frac{1}{2}} \right\|_{\mathbb{L}^r} < \infty, \tag{4.14}$$

then the integral $Y_t = \int_0^t \langle F_s, \mathrm{d}X_s \rangle$ is well-defined. Indeed, for any $F \in \mathbb{S}_{\mathrm{ad}}$, we get

$$\begin{aligned}\omega_j^{(\frac{r}{2})}(|T_\tau^{(r)}(\delta Y_{j+1})|^2) &= \sum_{\alpha,\beta} T_\tau^{(r)}(F_j(v_\alpha))^* T_\tau^{(r)}(F_j(v_\beta))\, \omega_j\Big(\delta X_j(v_\alpha)^* [\delta X_j(v_\beta)]_{\frac{1}{r}+2\tau}\Big) \\ &\leqslant \delta t_j \left| \sum_\alpha T_\tau^{(r)}(F_j(v_\alpha))\, \|A^{1/2} v_\alpha\|_{\mathfrak{h}} \right|^2,\end{aligned}$$

from which one can conclude along the lines of the proof of Theorem 4.15.

### 4.1.3 Itô–Grassmann integral

Now we want to focus on processes belonging to a Grassmann algebra (see Definition 3.7) generated by some analytic Grassmann Brownian motion $X$ and an independent and compatible Grassmann algebra $\tilde{\mathscr{M}}_0$, see Definition 3.9 and Lemma 3.10.

**Definition 4.17.** *Consider the filtered modular space $(\mathscr{M}, \omega, (\mathscr{M}_t)_t)$, and let $X: \mathbb{R}_+ \times \mathfrak{h} \to \mathscr{M}_a$ be an analytic GBM, and suppose that $\mathscr{M}_0$ contains a Grassmann algebra $\tilde{\mathscr{M}}_0$ of analytic elements which are independent of and compatible with the Grassmann algebra generated by X. Then we call the space*

$$\mathcal{G}_X := \mathrm{span}\{a_0 X_{t_1}(v_1) \cdots X_{t_n}(v_n),\ \textit{where } a_0 \in \tilde{\mathscr{M}}_0,\ n \in \mathbb{N}_0,\ v_1, \ldots, v_n \in \mathfrak{h},\ t_1, \ldots, t_n \in \mathbb{R}_+\},$$

*the Grassmann algebra generated by ($\tilde{\mathscr{M}}_0$ and) $X$. For every $p\in[1,\infty]$ we let $\mathcal{G}_X^p$ be the closure of $\mathcal{G}_X$ with respect to the $\mathbb{L}^p$ topology.*

**Remark 4.18.** The algebra $\tilde{\mathscr{M}}_0$ will be needed to study SDE with non-trivial initial conditions, see Sections 4.4 and 5.2.

The space $\mathcal{G}_X$ is a Grassmann algebra, and thus we can define the set of *even elements* $\mathcal{G}_{X,+}$ and *odd elements* $\mathcal{G}_{X,-}$ and we have $\mathcal{G}_X=\mathcal{G}_{X,+}\oplus\mathcal{G}_{X,-}$ , see Remark 3.8. Likewise, we can write $\mathcal{G}_X^p=\mathcal{G}_{X,+}^p\oplus\mathcal{G}_{X,-}^p$ where $\mathcal{G}_{X,+}^p$and $\mathcal{G}_{X,-}^p$ are the closure of $\mathcal{G}_{X,+}$ and $\mathcal{G}_{X,-}$ with respect the topology of $\mathbb{L}^p$. If an adapted process takes values in $\mathcal{G}_{X,+}^p$ or $\mathcal{G}_{X,-}^p$ for some $p\in[1,\infty]$, we will henceforth say that it is an even or odd process respectively. In particular we have that, for any $v\in\mathfrak{h}$, $X(v)$ is an odd process. We have the following fundamental result.

**Proposition 4.19.** *Let $X$ be an analytic GBM with covariance $G_t=tG$ and let $F\in\mathbb{H}^p_{|G|;\mathrm{loc}}$ be an adapted process taking values in $\mathcal{G}_X^p$. Then, for any $t\in\mathbb{R}_+$, $\int_0^t\langle F_s,\mathrm{d}X_s\rangle\in\mathcal{G}_X^p$. Furthermore if $F_s$ takes values in $\mathcal{G}_{X,+}^p$ then $\int_0^t\langle F_s,\mathrm{d}X_s\rangle\in\mathcal{G}_{X,-}^p$ and if $F_s$ takes values in $\mathcal{G}_{X,-}^p$ then $\int_0^t\langle F_s,\mathrm{d}X_s\rangle\in\mathcal{G}_{X,+}^p$.*

**Proof.** The thesis is obvious in the case where $F_s$ is a finite-dimensional simple process. The generic case can be obtained by the density of finite-dimensional simple processes in $\mathbb{H}^p_{|G|;\mathrm{loc}}$ and the continuity of the Itô integral in $\mathbb{L}^p$ (an thus in $\mathcal{G}_X^p$). □

We conclude this section by introducing the notion of Itô processes with respect to a GBM, which, as shown in Lemma 3.18 satisfy the assumption of Theorem 4.15.

**Definition 4.20.** *(Itô process). Let $X$ be an analytic GBM with covariance $G_t=tG$. We say that $Y\in C^0(\mathbb{R}_+,\mathcal{G}_X^p)$ is a Itô process (of integrability $p\geqslant 2$) if $Y_0\in\mathcal{G}_X^p$, $H\in\mathbb{H}^p_{|G|;\mathrm{loc}}\cap\mathcal{G}_X^p$, and $K\colon\mathbb{R}_+\to\mathcal{G}_X^p$, $\|K\|_{\mathbb{L}^p}\in L^1_{\mathrm{loc}}(\mathbb{R}_+)$ such that*

$$Y_t=Y_0+\int_0^t\langle H_s,\mathrm{d}X_s\rangle+\int_0^t K_s\,\mathrm{d}s \tag{4.15}$$

*where the second integral is in the Bochner sense. If $H_s\in\mathcal{G}_{X,+}^p$, $K_s\in\mathcal{G}_{X,-}^p$ (resp. $H_s\in\mathcal{G}_{X,-}^p$, $K_s\in\mathcal{G}_{X,+}^p$) we call $Y$ an odd (resp. even) Itô process.*

**Remark 4.21.** The above definition can be extended to the case of a multidimensional process $X$ indexed in (pre-)Hilbert space $W$. Indeed, we call the map $Y_\cdot\colon\mathbb{R}_+\times W\to\mathcal{G}_X^p$ an Itô process indexed in $W$ and of integrability $p\geqslant 2$), a map that for any $t\in\mathbb{R}_+$, $Y_t(\cdot)$ is a linear operator, for any $w\in W$, $Y_\cdot(w)\in C^0(\mathbb{R}_+,\mathcal{G}_X^p)$, and there are $H\colon\mathbb{R}_+\times W\times K\to\mathcal{G}_X^p$ and $K\colon\mathbb{R}_+\times W\to\mathcal{G}_X^p$, which are linear in the $W$ variable, and for any $w\in W$ we have $H_\cdot(w,\cdot)\in\mathbb{H}^p_{|G|;\mathrm{loc}}\cap\mathcal{G}_X^p$, $\|K_\cdot(w)\|_{\mathbb{L}^p}\in L^1_{\mathrm{loc}}(\mathbb{R}_+)$, for which for any $t\in\mathbb{R}_+$ and $w\in W$ we have

$$Y_t(w)=Y_0(w)+\int_0^t\langle H_s(w,\cdot),\mathrm{d}X_s\rangle+\int_0^t K_s(w)\mathrm{d}s. \tag{4.16}$$

If, for any $s\in\mathbb{R}_+$, $v\in\mathfrak{h}$ and $w\in W$, we have $F_s(w,v)\in\mathcal{G}_{X,+}^p$, $K_s(w)\in\mathcal{G}_{X,-}^p$ (resp. $F_s\in\mathcal{G}_{X,-}^p$, $K_s\in\mathcal{G}_{X,+}^p$) we call $Y$ an odd (resp. even) Itô random field.

**Remark 4.22.** The fact that $Y\in C^0(\mathbb{R}_+,\mathcal{G}_X^p)$ or $Y(w)\in C^0(\mathbb{R}_+,\mathcal{G}_X^p)$ follows directly from the representations (4.15) and (4.16). Furthermore, by Proposition 4.19, an odd/even Itô process is made (for fixed $s\in\mathbb{R}_+$ and $w\in W$) by odd/even elements of $\mathcal{G}_X^p$.

**Remark 4.23.** If $Y_t$ is an Itô process in $\mathbb{L}^p$ then, as consequence of Theorem 4.15, for every $T>0$, we have

$$\lim_{\pi\in\{\text{partitions of }[0,T]\},|\pi|\to 0}\left(\sup_{t_i\in\pi}\sup_{t\in[t_i,t_{i-1}]}\|Y_{t_{i-1}}-Y_s\|_{\mathbb{L}^{2p}}\right)=0$$

where the limit is taken with respect to the diameter $|\pi|$ of partitions $\pi$ going to 0. This property is analogous to the uniform continuity in $L^p$ for commutative stochastic processes.

## 4.2 Itô formula

In this section we want to prove an Itô formula for Itô processes in the sense of Definition 4.20. We remark that the following discussion could be extended to the anticommutative analogous of continuous commutative martingales - see, e.g., [RW00, Chapter IV, Section 5] for the latter notion.

**Definition 4.24.** *Let $p,q\in[1,\infty]$ and let $B,B'$ be two martingales in $C^0(\mathbb{R}_+,\mathbb{L}^p)$ and $C^0(\mathbb{R}_+,\mathbb{L}^q)$ respectively. The quadratic variation $[B,B']\in C^0(\mathbb{R}_+,\mathbb{L}^r)$ with $1/r=1/p+1/q$ is a process of bounded variation (in the $\mathbb{L}^r$ norm) such that the process*

$$t\mapsto B_t\cdot B_t'-[B,B']_t$$

*is an martingale.*

**Lemma 4.25.** *Let $X$ be an analytic GBM with covariance $G_t=tG$. Then, for any $v,v'\in\mathfrak{h}$ we have*

$$[X(v),X(v')]_t=\langle\Theta v,Gv'\rangle_{\mathfrak{h}}t.$$

**Proof.** The proof is a consequence of the fact that the processes $X_t(v)$ and $X_t(v')$ have independent increments and so $\omega_s((X_{s,t})(v)(X_{s,t})(v'))=\omega_s((X_{s,t})(v)(X_{s,t})(v'))=\langle\Theta\, v,Gv'\rangle(t-s)$, for any $s\leqslant t\in\mathbb{R}_+$. □

**Theorem 4.26.** *Let $Y_t,Y_t'$ be two Itô processes of integrability $p\geqslant 2$ of the form*

$$Y_t=Y_0+\int_0^t\langle H_s,\mathrm{d}X_s\rangle,\quad Y_t'=Y_0'+\int_0^t\langle H_s',\mathrm{d}X_s\rangle,$$

*where $X$ is an analytic GBM indexed in $\mathfrak{h}$ with covariance $G_t=tG$. Then, we have*

$$[Y,Y']_t=\int_0^t\mathrm{Tr}_{G^*\Theta}(H_s\otimes H_s')\mathrm{d}s, \tag{4.17}$$

*where $\mathrm{Tr}_{G^*\Theta}(\cdot)$ is defined in Notation 4.11, $\Theta$ being the conjugation in $\mathfrak{h}$.*

**Proof.** We prove the theorem in the special case of $Y,Y'$ being simple finite-dimensional processes. The general case follows from the linearity of the quadratic variation and by density. On the other hand, formula (4.17) follows directly from the definition of simple processes, the properties of conditional expectation and Lemma 4.25. □

It will be convenient to work in some subspaces of the spaces $\mathbb{H}_A^p$, see Definition 4.13.

**Definition 4.27.** *Let $F$ an $\mathbb{L}^p$ adapted process indexed in a Hilbert space $\mathfrak{h}$ with conjugation $\Theta$ and let $A\geqslant 0$ be a linear bounded operator on $\mathfrak{h}$ commuting with $\Theta$. We let*

$$\|F_s\|_{\mathbb{L}^p,A}:=\||F_s|_A\|_{\mathbb{L}^p}, \tag{4.18}$$

*and introduce the following subspace of $\mathbb{H}_{A,\mathrm{loc}}^p$*

$$\tilde{\mathbb{H}}_{A,\mathrm{loc}}^p:=\{F\in\mathbb{H}_{A,\mathrm{loc}}^p|\|F_s\|_{\mathbb{L}^p,A}\in L_{\mathrm{loc}}^2(\mathbb{R}_+)\}.$$

**Remark 4.28.** One can check that $\tilde{\mathbb{H}}_{A,\mathrm{loc}}^p$ is indeed a subspace of $\mathbb{H}_{A,\mathrm{loc}}^p$, in fact

$$\|F\|_{\mathbb{H}_A^p([0,t])}\lesssim\left(\int_0^t\|F_s\|_{\mathbb{L}^p,A}^{2p}\mathrm{d}s\right)^{\frac{1}{2}}.$$

We also note that $\tilde{\mathbb{H}}^p_{A,\mathrm{loc}}$ are semimetric spaces with the set of semidistances

$$d_t(F,G) \mapsto \Big(\int_0^t \|F_s - G_s\|^{2p}_{\mathbb{L}^p,A}\,\mathrm{d}s\Big)^{\frac{1}{2}}.$$

Furthermore the finite-dimensional processes taking values in $\mathbb{L}^p$ defined in Definition 4.9 are dense in $\tilde{\mathbb{H}}^p_{A,\mathrm{loc}}$ with respect to its natural topology, compare with Remark 4.14.

**Corollary 4.29.** *Under the hypothesis of Theorem 4.26, the following inequality holds*

$$\int_0^t \|\mathrm{Tr}_{G^*\Theta}(H_s \otimes H'_s)\|_{\mathbb{L}^r}\,\mathrm{d}s \lesssim \|H\|^2_{\tilde{\mathbb{H}}^p_{|G|}([0,t])} + \|H'\|^2_{\tilde{\mathbb{H}}^q_{|G|}([0,t])}. \tag{4.19}$$

**Proof.** We will use the representation (4.17). Consider first $H \in \mathbb{S}_{\mathrm{ad}}$ and $H \in \mathbb{S}_{\mathrm{ad}}$ and let $[s,t] \subset \mathbb{R}_+$ an interval where both $H$ and $H'$ are constant so that $H_r \cdot H_{r'} \in \mathscr{M}_a \cap \mathscr{M}_s$ for $r \in [s,t]$. Then, by Theorem 4.26 and by using that $t \mapsto Y_t \cdot Y'_t - [Y,Y']_t$ is a martingale

$$\begin{aligned}\int_s^t \|\mathrm{Tr}_{G^*\Theta}(H_r \otimes H'_r)\|_{\mathbb{L}^r}\,\mathrm{d}r &= \Big\|\omega_s\Big(\int_s^t \mathrm{Tr}_{G^*\Theta}(H_r \otimes H'_r)\,\mathrm{d}r\Big)\Big\|_{\mathbb{L}^r}\\ &= \|\omega_s([Y,Y']_t - [Y,Y']_s)\|_{\mathbb{L}^r}\\ &= \Big\|\omega_s\Big(\int_s^t \langle H_r, \mathrm{d}X_r\rangle\Big)\Big(\int_s^t \langle H'_{r'}, \mathrm{d}X_{r'}\rangle\Big)\Big\|_{\mathbb{L}^r}.\end{aligned}$$

By Hölder's inequality, Theorem 4.15 and the bound $\|\cdot\|^2_{\mathbb{H}^p_A([0,t])} \leqslant \|\cdot\|^2_{\tilde{\mathbb{H}}^p_A([0,t])}$

$$\begin{aligned}\int_s^t \|\mathrm{Tr}_{G^*\Theta}(H_r \otimes H'_r)\|_{\mathbb{L}^r}\,\mathrm{d}r &\leqslant \Big\|\Big(\int_s^t \langle H_r, \mathrm{d}X_r\rangle\Big)\Big(\int_s^t \langle H'_{r'}, \mathrm{d}X_{r'}\rangle\Big)\Big\|_{\mathbb{L}^r}\\ &\lesssim \Big(\int_s^t \|H_r\|^2_{\mathbb{L}^p,|G|}\,\mathrm{d}r\Big) + \Big(\int_s^t \|H'_s\|^2_{\mathbb{L}^q,|G|}\,\mathrm{d}r\Big),\end{aligned}$$

Taking the sum over a partition where both $H_s$ and $H'_s$ are constant, we obtain (4.19). □

We introduce the following standard notation.

**Notation 4.30.** *If $Y, Y'$ are Itô processes of the form $Y_t = Y_0 + \int_0^t \langle H_s, \mathrm{d}X_s\rangle + \int_0^t K_s\,\mathrm{d}s$, $Y'_t = Y'_0 + \int_0^t \langle H'_s, \mathrm{d}X_s\rangle + \int_0^t K'_s\,\mathrm{d}s$, and $\tilde{H}$ and $\tilde{K}$ are adapted processes, we use the notation*

$$\int_0^t \tilde{H}_s\,\mathrm{d}Y_s := \int_0^t \langle \tilde{H}_s H_s, \mathrm{d}X_s\rangle + \int_0^t \tilde{H}_s K_s\,\mathrm{d}s, \qquad \int_0^t \tilde{K}_s\,\mathrm{d}[Y,Y']_s := \int_0^t \tilde{K}_s \mathrm{Tr}_{G^*\Theta}(H_s \otimes H'_s)\mathrm{d}s.$$

We want to prove an Itô formula for Itô processes indexed in some finite-dimensional linear space $W$. To this end, we introduce the space of polynomials of degree at most $k \in \mathbb{N}$ of commutative and anticommutative random variables indexed in $W$:

$$\mathcal{F}^k(W) = \bigoplus_{n=0}^k \bigoplus_{\ell=0}^n \Lambda^{n-\ell} W \otimes \odot^\ell W,$$

i.e. the tensor product of the antisymmetric and symmetric vector spaces generated by $W$. We also use the notation $\mathcal{F}^\infty(W) = \bigcup_{k\in\mathbb{N}} \mathcal{F}^k(W)$. The space $\mathcal{F}^\infty(W)$ is an algebra with the following multiplication: if $F \in \mathcal{F}^k(W)$, that is,

$$F = \sum_{n=0}^k \sum_{\ell=0}^n F^{(n,\ell)} \otimes \hat{F}^{(n,\ell)}, \tag{4.20}$$

where $F^{(n,\ell)} \in \Lambda^{n-\ell} W$ and $\hat{F}^{(n,\ell)} \in \bigodot^\ell W$, for any $G \in \Lambda^m W$ and $G' \in \bigodot^{m'} W$ we set

$$(G \otimes G') \cdot F := \sum_{n=0}^k \sum_{\ell=0}^n (G \wedge F^{(n,\ell)}) \otimes (G' \odot \hat{F}^{(n,\ell)}). \tag{4.21}$$

If $Z_{-}: W \to \mathcal{G}^p_{X,-}$ is an odd Grassmann random field and $F = \sum_i w_{\alpha_i^1} \wedge \cdots \wedge w_{\alpha_i^k} \in \Lambda^k W$, $n \leqslant p$, where $w_{\alpha_i^n} \in W$, we can define

$$F(Z_-) = \sum_i Z_-(w_{\alpha_i^1}) \cdots Z_-(w_{\alpha_i^k})$$

which is a well-defined object (not depending on the explicit representation of $F$) since $Z$ is odd. In the same way if $Z_+$ is even and $F = \sum_i w_{\alpha_i^1} \odot \cdots \odot w_{\alpha_i^k} \in \bigodot^k W$ we can define $F(Z_-) = \sum_i Z_+(w_{\alpha_i^1}) \cdots Z_+(w_{\alpha_i^k})$. In this way if $Z: W \to \mathcal{G}^p_X$ is a generic random field, which can be decomposed in a unique way $Z = Z_- + Z_+$, where $Z_-, Z_+$ are even and odd respectively, and considering a generic element $F \in \mathcal{F}^k(W)$, see (4.20), we define

$$F(Z) := \sum_{n=0}^{k} \sum_{\ell=0}^{n} F^{(n,\ell)}(Z_-) \cdot \hat{F}^{(n,\ell)}(Z_+). \tag{4.22}$$

With the previous notation we have that if $F, G \in \mathcal{F}^\infty(W)$ and $Z$ has enough integrability we have $(F \cdot G)(Z) = F(Z) \cdot G(Z)$ (where the first product is the product of polynomials in the sense of definition (4.21), and the second is the product in the twisted space $\mathbb{L}^p$).

We introduce commutative and anticommutative derivatives on $\mathcal{F}^k(W)$.

**Definition 4.31.** *For any $w \in W$, we define two linear maps $\partial^a_w: \mathcal{F}^k(W) \to \mathcal{F}^{k-1}(W)$ and $\partial^c_w: \mathcal{F}^k(W) \to \mathcal{F}^{k-1}(W)$ by the following recursive relation, for any $v \in W$,*

$$\partial^a_w(1 \otimes v) = \partial^c_w(v \otimes 1) = 0,$$

$$\partial^a_w((v \otimes 1) \cdot F) = \langle w, v \rangle_W F - (w \otimes 1) \cdot \partial^a_w F, \qquad \partial^c_w((1 \otimes v) \cdot F) = \langle w, v \rangle_h F + (1 \otimes w) \cdot \partial^a_w F.$$

This allows us to provide a general Taylor formula for polynomials in $\mathcal{F}^k(W)$.

**Lemma 4.32.** *Consider $F \in \mathcal{F}^k(W)$, let $Z = Z_- + Z_+, R = R_- + R_+: W \to \mathcal{G}^p_X$ be random fields indexed in $h$ (where $Z_-, R_-$ and $Z_+, R_+$ are respectively even and odd ), and let $\{w_j\}_{j=1}^N$ be an orthonormal basis of $W$, then*

$$F(Z) - F(R) = \sum_{n=0}^{k} \sum_{\ell=0}^{n} \sum_{j_1, \ldots, j_{n-\ell}, j'_1, \ldots, j'_\ell = 1}^{N} \frac{1}{(n-\ell)!\,\ell!} \partial^a_{w_{j_1}} \partial^a_{w_{j_2}} \cdots \partial^a_{w_{j_{n-\ell}}} \partial^c_{w_{j'_1}} \cdots \partial^c_{w_{j'_\ell}}(F)(R) \cdot$$
$$\cdot \left( \prod_{o=1}^{n-\ell} (Z_-(w_{j_o}) - R_-(w_{j_o})) \right) \cdot \left( \prod_{o=1}^{\ell} (Z_+(w_{j'_o}) - R_+(w_{j'_o})) \right).$$

We omit the proof, it being a simple generalisation of the usual Taylor formula. We are now ready to prove an Itô formula for polynomials of Itô processes.

**Theorem 4.33.** *Let $\mathfrak{h}$ be a Hilbert space with conjugation $\Theta$. Let $X$ is an analytic GBM indexed in $\mathfrak{h}$ with covariance $G_t = tG$. Consider $Y_t = Y_{t,-} + Y_{t,+}$ be a generic Itô process of integrability $p \geqslant 2$, indexed in a finite-dimensional vector space $W$, where $Y_{t,-}$ and $Y_{t,+}$ are an odd and even Itô random fields such that*

$$Y_\rho = Y_{0,\rho} + \int_0^t K_{s,\rho}(\cdot) \mathrm{d}s + \int_0^t \langle H_{s,\rho}(\cdot,\cdot), \mathrm{d}X_s \rangle, \qquad \rho = \pm$$

*where $K_-, K_+: W \to \mathcal{G}^p_X$ are odd and even processes respectively and $H_{s,-}(w,\cdot), H_{s,+}(w,\cdot) \in \mathbb{H}^p_{|G|,\mathrm{loc}}$ are even and odd processes respectively. Consider $F \in \mathcal{F}^n(W)$, where $n \leqslant \frac{p}{2}$, then, for any orthonormal basis $\{v_j\}_{\alpha=1}^N$ of $W$, we have*

$$F(Y_t) - F(Y_0) = \sum_\alpha \int_0^t \partial^a_{v_\alpha} F(Y_s) \, \mathrm{d}Y^\alpha_{s,-} + \int_0^t \partial^c_{v_\alpha} F(Y_s) \, \mathrm{d}Y^\alpha_{s,+} + \frac{1}{2} \sum_{\alpha,\beta} \int_0^t \partial^a_{v_\alpha} \partial^a_{v_\beta} F(Y_s) \, \mathrm{d}[Y^\alpha_-, Y^\beta_-]_s$$
$$+ \frac{1}{2} \sum_{\alpha,\beta} \int_0^t \partial^c_{v_\alpha} \partial^c_{v_\beta} F(Y_s) \, \mathrm{d}[Y^\alpha_+, Y^\beta_+]_s + \sum_{\alpha,\beta} \int_0^t \partial^c_{v_\alpha} \partial^a_{v_\beta} F(Y_s) \, \mathrm{d}[Y^\alpha_+, Y^\beta_-]_s \tag{4.23}$$

*where* $Y_\rho^\alpha = Y_\rho(v_\alpha)$, $\rho = \pm$, *and the integrals are defined in function of* $X$, $K_-$, $K_+$, $H_-$ *and* $H_+$ *as in Notation 4.30.*

**Remark 4.34.** If $Y$ is odd and $F \in \bigoplus_{n=0}^k \Lambda^n W$ the composition $F(Y_t)$ is well-defined (as a particular case of equation (4.22)) and Theorem 4.33 implies

$$F(Y_t) - F(Y_0) = \sum_\alpha \int_0^t \partial_{v_\alpha}^a F(Y_s) \mathrm{d}Y_s^\alpha + \frac{1}{2}\sum_{\alpha,\beta} \int_0^t \partial_{v_\alpha}^a \partial_{v_\beta}^a F(Y_s) \mathrm{d}[Y^\alpha, Y^\beta]_s, \tag{4.24}$$

which is the perfect analogous of the Itô formula in standard stochastic calculus. If $F$ is even and $F \in \bigoplus_{n=0}^k \bigodot^n W$ Itô formula (4.24) holds with $\partial_{v_\alpha}^c$ in place of $\partial_{v_\beta}^a$.

**Proof of Theorem 4.33.** The proof is along the lines of the standard proof of the Itô formula in the commutative setting (see, e.g., [Le 16, Section 5.2]). For this reason, we sketch it in the special case $K_- = K_+ = H_+ = 0$. Furthermore, we consider the convergence in $\mathbb{L}^2$ since $\mathbb{L}^p \subset \mathbb{L}^2$ for any $p > 2$; in fact, by proving the equality in $\mathbb{L}^2$, because both members belong to $\mathbb{L}^p$ we have that the equality holds in $\mathbb{L}^p$ as well.

We write $Y_t^\alpha = Y_{t,-}(v_\alpha)$, and $H_s^\alpha(\cdot) = H_{s,-}(v_\alpha, \cdot)$. Furthermore, since in the proof only anticommutative derivatives $\partial^a$ are needed, we use the notation by $\partial_{v_\alpha} \equiv \partial_{v_\alpha}^a$, $\partial_{v_{\alpha_1} \dots v_{\alpha_i}}^i \equiv \partial_{v_{\alpha_1}}^a \cdots \partial_{v_{\alpha_i}}^a$. Being $F$ a polynomial of degree at most $n$, for any $\alpha_1, \dots, \alpha_i$ we have

$$\|\partial_{v_{\alpha_1} \dots v_{\alpha_i}}^i F(Y_t)\|_{\mathbb{L}^{2n/(n-i)}} \lesssim (1 + \|Y_t\|_{\mathbb{L}^{2n}})^{n-i}, \tag{4.25}$$

where the constant in the symbol $\lesssim$ depends only on $F$. Furthermore, by Lemma 4.32, we get

$$\begin{aligned} F(Y_t) - F(Y_0) &= \sum_{t_i \in \pi} F(Y_{t_i}) - F(Y_{t_{i-1}}) \\ &= \sum_{t_i \in \pi} \sum_\alpha \partial_{v_\alpha} F(Y_{t_{i-1}})(Y_{t_i}^\alpha - Y_{t_{i-1}}^\alpha) + \frac{1}{2} \sum_{t_i \in \pi} \sum_{\alpha,\beta} \partial_{v_\alpha v_\beta}^2 F(Y_{t_{i-1}})(Y_{t_i}^\alpha - Y_{t_{i-1}}^\alpha)(Y_{t_i}^\beta - Y_{t_{i-1}}^\beta) \\ &\quad + \sum_{t_i \in \pi} \sum_{i=3}^n \frac{1}{i!} \sum_{\alpha_1, \dots, \alpha_i = 1}^k \partial_{v_{\alpha_1} \dots v_{\alpha_i}}^i F(Y_{t_{i-1}})(Y_{t_i}^{\alpha_1} - Y_{t_{i-1}}^{\alpha_1}) \cdots (Y_{t_i}^{\alpha_i} - Y_{t_{i-1}}^{\alpha_i}). \end{aligned} \tag{4.26}$$

Furthermore, combining (4.25) and (4.26), we get

$$\|\partial_{v_{\alpha_1} \dots v_{\alpha_i}}^i F(Y_t) - \partial_{v_{\alpha_1} \dots v_{\alpha_i}}^i F(Y_s)\|_{\mathbb{L}^{2n/(n-i)}} \lesssim (1 + \|Y_t\|_{\mathbb{L}^{2n}} + \|Y_s\|_{\mathbb{L}^{2n}})^{n-i} \|Y_t - Y_s\|_{\mathbb{L}^{2n}}. \tag{4.27}$$

Let us focus on each term in the sum in (4.26). For any $|\tau| \leqslant \frac{3}{4}$, we have that

$$\begin{aligned} &\left\| T_\tau^{(2)} \left( \sum_{t_i \in \pi} \partial_{v_\alpha} F(Y_{t_{i-1}})(Y_{t_i}^\alpha - Y_{t_{i-1}}^\alpha) - \int_0^t \langle \partial_{v_\alpha} F(Y_s) H_{s,-}, \mathrm{d}X_s \rangle \right) \right\|_{L^2}^2 \\ &= \left\| T_\tau^{(2)} \left( \sum_{t_i \in \pi} \left( \int_{t_{i-1}}^{t_i} \langle (\partial_{v_\alpha} F(Y_{t_{i-1}}) - \partial_{v_\alpha} F(Y_s)) H_{s,-}, \mathrm{d}X_s \rangle \right) \right) \right\|_{L^2}^2 \\ &\lesssim \sum_{t_i \in \pi} \int_{t_{i-1}}^{t_i} \left\| T_{\tau_1}^{\left(\frac{2n}{n-1}\right)} (\partial_{v_\alpha} F(Y_{t_{i-1}}) - \partial_{v_\alpha} F(Y_s)) \cdot \right. \\ &\qquad \left. \cdot \mathrm{Tr}_{G^* \Theta}(T_{\tau_2}^{(2n)}(H_s) \otimes T_{\tau_2}^{(2n)}(H_s^*)) T_{\tau_1}^{\left(\frac{2n}{n-1}\right)} (\partial_{v_\alpha} F(Y_{t_{i-1}}) - \partial_{v_\alpha} F(Y_s)) \right\|_{\mathbb{L}^1} \mathrm{d}s \\ &\lesssim \sum_{t_i \in \pi} \sup_{s \in [t_i, t_{i-1}]} \|\partial_{v_\alpha} F(Y_{t_{i-1}}) - \partial_{v_\alpha} F(Y_s)\|_{\mathbb{L}^{2n/(n-1)}}^2 \left( \sum_\alpha \int_{t_{i-1}}^{t_i} \|\mathrm{Tr}_{G^* \Theta}(H_s^\alpha \otimes H_s^{\alpha *})\|_{\mathbb{L}^{2n}} \mathrm{d}s \right) \\ &\lesssim \left( \sum_\alpha \int_0^t \|H_s^\alpha\|_{\mathbb{L}^p, |G|} \mathrm{d}s \right) \left( \sup_{s \in [0,t]} (1 + \|Y_t\|_{\mathbb{L}^p}) \right)^{2n-4} \left( \max_{t_i \in \pi} \sup_{s \in [t_i, t_{i-1}]} \|Y_{t_{i-1}} - Y_s\|_{\mathbb{L}^p}^2 \right), \end{aligned} \tag{4.28}$$

where $\tau_1, \tau_2 \geqslant 0$ are suitable real numbers, and in the last step we used the inequalities (4.19) and (4.27).

By the inequality (4.27) and by Remark 4.23 applied to the Itô process $Y_t$, the last term goes to zero. Taking the sup over $|\tau| \leqslant \frac{3}{4}$ in (4.28), we obtain that $\sum_{t_i \in \pi} \partial_{v_\alpha} F(Y_{t_{i-1}})(Y_{t_i}^\alpha - Y_{t_{i-1}}^\alpha)$ converges to $\int_0^t \langle \partial_{v_\alpha} F(Y_s) H_s^\alpha(\cdot), \mathrm{d}X_s \rangle$ as $|\pi| \to 0$ in $\mathbb{L}^2$.

By using the fact that $Y_s^\alpha Y_s^\beta - [Y^\alpha, Y^\beta]_s$ is a martingale, see Definition 4.24, we get that, for any $|\tau| \leqslant \frac{3}{4}$,

$$\left\| T_\tau^{(2)} \left( \sum_{t_i \in \pi} \partial^2_{v_\alpha v_\beta} F(Y_{t_{i-1}})((Y_{t_i}^\alpha - Y_{t_{i-1}}^\alpha)(Y_{t_i}^\beta - Y_{t_{i-1}}^\beta) - [Y^\alpha, Y^\beta]_{t_i} + [Y^\alpha, Y^\beta]_{t_{i-1}}) \right) \right\|_{L^2}^2 =$$

$$= \sum_{t_i \in \pi} \left\| T_{\tau_1}^{(n)}((Y_{t_i}^\alpha - Y_{t_{i-1}}^\alpha)(Y_{t_i}^\beta - Y_{t_{i-1}}^\beta) - [Y^\alpha, Y^\beta]_{t_i} + [Y^\alpha, Y^\beta]_{t_{i-1}}) T_{\tau_2}^{\left(\frac{2n}{n-2}\right)}(\partial^2_{v_\alpha v_\beta} F(Y_{t_{i-1}})) \cdot \right.$$

$$\left. \cdot T_{\tau_2}^{\left(\frac{2n}{n-1}\right)}(\partial^2_{v_\alpha v_\beta} F(Y_{t_{i-1}}))^* T_{\tau_2}^{(n)}((Y_{t_i}^\alpha - Y_{t_{i-1}}^\alpha)(Y_{t_i}^\beta - Y_{t_{i-1}}^\beta) - [Y^\alpha, Y^\beta]_{t_i} + [Y^\alpha, Y^\beta]_{t_{i-1}})^* \right\|_{L^1}$$

$$\lesssim \sum_{t_i \in \pi} \|(Y_{t_i}^\alpha - Y_{t_{i-1}}^\alpha)(Y_{t_i}^\beta - Y_{t_{i-1}}^\beta) - [Y^\alpha, Y^\beta]_{t_i} + [Y^\alpha, Y^\beta]_{t_{i-1}}\|_{\mathbb{L}^n}^2 \|\partial^2_{v_\alpha v_\beta} F(Y_{t_{i-1}})\|_{\mathbb{L}^{2n/(n-2)}}^2$$

$$\lesssim \left( \sup_{s \in [0,t]} (1 + \|Y_t\|_{\mathbb{L}^{2n}}) \right)^{2n-4} \left( \sum_{t_i \in \pi} \|(Y_{t_i}^\alpha - Y_{t_{i-1}}^\alpha)(Y_{t_i}^\beta - Y_{t_{i-1}}^\beta)\|_{\mathbb{L}^n}^2 + \|[Y^\alpha, Y^\beta]_{t_i} - [Y^\alpha, Y^\beta]_{t_{i-1}}\|_{\mathbb{L}^n}^2 \right)$$

$$\lesssim \left( \sup_{s \in [0,t]} (1 + \|Y_t\|_{\mathbb{L}^p}) \right)^{2n-4} \left( \sum_{t_i \in \pi} \left( \|H^\alpha\|^2_{\mathbb{H}^n_{|G|}([t_{i-1}, t_i])} + \|H^\beta\|^2_{\mathbb{H}^n_{|G|}([t_{i-1}, t_i])} \right)^2 \right) \quad (4.29)$$

where we used that $n \leqslant p/2$, which goes to zero as $|\pi| \to 0$ since the map $t \mapsto \|H^\alpha\|^2_{\mathbb{H}^p_{|G|}([0,t])}$ (and thus the map $t \mapsto \|H^\alpha\|^2_{\mathbb{H}^n_{|G|}([0,t])}$) is an absolutely continuous function. On the other hand, we have that

$$\left\| T_\tau^{(2)} \left( \sum_{t_i \in \pi} \partial^2_{v_\alpha v_\beta} F(Y_{t_{i-1}})([Y^\alpha, Y^\beta]_{t_i} - [Y^\alpha, Y^\beta]_{t_{i-1}}) - \int_0^t \partial^2_{v_\alpha v_\beta} F(Y_s) \mathrm{Tr}_{G^*\Theta}(H_s^\alpha(\cdot) \otimes H_s^\beta(\cdot)) \mathrm{d}s \right) \right\|_{L^2}$$

$$= \left\| T_\tau^{(2)} \left( \sum_{t_i \in \pi} \int_{t_{i-1}}^{t_i} (\partial^2_{v_\alpha v_\beta} F(Y_{t_{i-1}}) - \partial^2_{v_\alpha v_\beta} F(Y_s)) \mathrm{Tr}_{G^*\Theta}(H_s^\alpha(\cdot) \otimes H_s^\beta(\cdot)) \mathrm{d}s \right) \right\|_{L^2}$$

$$\lesssim \sum_{t_i \in \pi} \int_{t_{i-1}}^{t_i} \|(\partial^2_{v_\alpha v_\beta} F(Y_{t_{i-1}}) - \partial^2_{v_\alpha v_\beta} F(Y_s))\|_{\mathbb{L}^{2n/(n-2)}} \|\mathrm{Tr}_{G^*\Theta}(H_s^\alpha(\cdot) \otimes H_s^\beta(\cdot))\|_{\mathbb{L}^n} \mathrm{d}s$$

$$\lesssim \left( \sup_{s \in [0,t]} (1 + \|Y_t\|_{\mathbb{L}^{2p}}) \right)^{n-3} \left( \sup_{t_i \in \pi} \sup_{t \in [t_i, t_{i-1}]} \|Y_{t_{i-1}} - Y_s\|_{\mathbb{L}^{2p}} \right) \left( \|H^\alpha(\cdot)\|^2_{\mathbb{H}^p_{|G|}([0,t])} + \|H^\beta(\cdot)\|^2_{\mathbb{H}^p_{|G|}([0,t])} \right) \quad (4.30)$$

which converges to 0 by Remark 4.23 applied to the Itô process $Y_t$. Thus, taking the supremum over $|\tau| \leqslant \frac{3}{4}$, the inequalities (4.29) and (4.29) imply that $\frac{1}{2} \sum_{t_i \in \pi} \partial^2_{v_\alpha v_\beta} F(Y_{t_{i-1}})(Y_{t_i}^\alpha - Y_{t_{i-1}}^\alpha)(Y_{t_i}^\beta - Y_{t_{i-1}}^\beta)$ converges to $\int_0^t \partial^2_{v_\alpha, v_\beta} F(Y_s) \mathrm{Tr}_{G^*\Theta}(H_s^\alpha(\cdot) \otimes H_s^\beta(\cdot)) \mathrm{d}s$ in $\mathbb{L}^2$. Finally we have

$$\begin{aligned}
&\left\| \sum_{t_i \in \pi} \partial^k_{v_{\alpha_1} \ldots v_{\alpha_k}} F(Y_{t_{i-1}})(Y_{t_i}^{\alpha_1} - Y_{t_{i-1}}^{\alpha_1}) \cdots (Y_{t_i}^{\alpha_i} - Y_{t_{i-1}}^{\alpha_i}) \right\|_{\mathbb{L}^2} \\
&\quad \leqslant \sum_{t_i \in \pi} \|\partial^k_{v_{\alpha_1} \ldots v_{\alpha_k}} F(Y_{t_{i-1}})(Y_{t_i}^{\alpha_1} - Y_{t_{i-1}}^{\alpha_1}) \cdots (Y_{t_i}^{\alpha_i} - Y_{t_{i-1}}^{\alpha_i})\|_{\mathbb{L}^2} \\
&\quad \lesssim \sup_{s \in [0,t]} (1 + \|Y_s\|_{\mathbb{L}^p})^{n-k} \sum_{t_i \in \pi} \|Y_{t_i} - Y_{t_{i-1}}\|_{\mathbb{L}^p}^k \\
&\quad \lesssim \sup_{s \in [0,t]} (1 + \|Y_s\|_{\mathbb{L}^p})^{n-k} \sum_\alpha \left( \sum_{t_i \in \pi} \left( \int_{t_{i-1}}^{t_i} \mathrm{Tr}_{|G|}(H_s^\alpha(\cdot) \otimes H_s^{\alpha *}(\cdot)) \mathrm{d}s \right)^{\frac{k}{2}} \right).
\end{aligned}$$

The last term converges to 0 since $k>2$. This means that

$$\begin{aligned}F(Y_t)-F(Y_0) &= \lim_{|\pi|\to 0}\sum_{t_i\in\pi} F(Y_{t_i})-F(Y_{t_{i-1}})\\ &= \sum_{\alpha}\int_0^t \langle \partial_{v_\alpha}F(Y_s)H_s^\alpha(\cdot),\mathrm{d}X_s\rangle+\frac{1}{2}\sum_{\alpha,\beta}\int_0^t \partial^2_{v_\alpha,v_\beta}F(Y_s)\mathrm{Tr}_{G^*\Theta}(H_s^\alpha(\cdot)\otimes H_s^\beta(\cdot))\mathrm{d}s.\end{aligned}$$

□

**Remark 4.35.** The fact that the Itô process $Y$ takes values in a Grassmann algebra (namely $\mathcal{G}_X^p$) is used in a essential way in the proof for obtaining the last three terms in the Itô formula (4.23), namely the terms involving the quadratic variation. Indeed in order to get the quadratic variation, we need to pull through the term $(Y_{t_i}^\alpha-Y_{t_{i-1}}^\alpha)(Y_{t_i}^\beta-Y_{t_{i-1}}^\beta)$ in the Taylor expansion which is possible thanks to the anticommutativity or commutativity of the $Y^\alpha$'s.

**Remark 4.36.** An important consequence of Theorem 4.33 is that sums and products of Itô processes taking values in a Grassmann algebra are again an Itô processes. In other words, Itô processes are closed with respect to algebraic operations.

## 4.3 Girsanov's formula

In this section we will prove a form of Girsanov's theorem involving Itô processes in the algebra $\bar{\mathcal{G}}_X := \bigcap_{p\geqslant 1}\mathcal{G}_X^p$ (see Remark 4.40 below for more details). Since $\bar{\mathcal{G}}_X$ is only a Fréchet algebra, but the notions of GBM and conditional expectation were so far discussed for modular spaces ($\mathscr{M}$, $\omega$, $(\mathscr{M}_t)_t$), we shall now extend them. First of all, we introduce the notion of signed expectation (or state) for a generic Fréchet algebra $\mathcal{A}$.

**Definition 4.37.** *Let $\mathcal{A}$ be a Fréchet algebra. We call a continuous linear functional $\bar{\mathcal{E}}:\mathcal{A}\to\mathbb{C}$ a signed expectation (or state) on $\mathcal{A}$ . Let $\{\mathcal{A}_t\}_{t\in\mathbb{R}_+}$ be a filtration of Fréchet subalgebra of $\mathcal{A}$ and let $\bar{\mathcal{E}}_t:\mathcal{A}\to\mathcal{A}_t$ be a family of linear continuous operators, we say $\bar{\mathcal{E}}_t$ is a conditional expectation associated with $\bar{\mathcal{E}}$ if for any $a\in\mathcal{A}$, $t,s\in\mathbb{R}_+$, $t\geqslant s\geqslant 0$, $b,c\in\mathcal{A}_s$ we have*

1. $\bar{\mathcal{E}}(\bar{\mathcal{E}}_t(a))=\bar{\mathcal{E}}(a)$,
2. $\bar{\mathcal{E}}_s(\bar{\mathcal{E}}_t(a))=\bar{\mathcal{E}}_s(a)$,
3. $\bar{\mathcal{E}}_s(bac)=b\bar{\mathcal{E}}_s(a)c$.

**Remark 4.38.** Note that in a generic Fréchet algebra (without a * operation) we do not have the notion of positive operator and thus of positive linear functional. Yet, in the particular case of a Fréchet *-algebra the expectation $\bar{\mathcal{E}}$ is generally not a positive linear functional.

A process $X\in\mathcal{A}$ is a martingale with respect to the signed expectation $\bar{\mathcal{E}}$ if, as usual,

$$\bar{\mathcal{E}}_s(X_t)=X_s,\qquad \forall s\leqslant t.$$

We can also define a concept of law for random fields with respect to a signed expectation $\bar{\mathcal{E}}$: we say that two random fields $B,B':\mathfrak{h}\to\mathcal{A}$, defined on the same vector space $\mathfrak{h}$, has the same law with respect to the signed state $\bar{\mathcal{E}}$ if and only if the $n$ points functions of the two fields are the same; i.e. for any $v_1,\ldots,v_n\in\mathfrak{h}$ we have

$$\bar{\mathcal{E}}(B(v_1)\cdots B(v_n))=\bar{\mathcal{E}}(B'(v_1)\cdots B'(v_n)).$$

We can extend the same notion for processes, and thus we are able to define the notion of a Brownian motion with covariance $G$ with respect to $\bar{\mathcal{E}}$ and of martingale.

**Definition 4.39.** *Let $\mathcal{A}$ be a Fréchet algebra with signed expectation $\bar{\mathcal{E}}$ and let $\mathcal{B}=\mathcal{B}_-\oplus\mathcal{B}_+\subset\mathcal{A}$ be a Grassmann subalgebra of $\mathcal{A}$. Let $\mathfrak{h}$ be a separable Hilbert space with conjugation $\Theta$. We say that a process $X\colon\mathbb{R}_+\times\mathfrak{h}\to\mathcal{B}_-$ is a weak Grassmann Brownian motion with respect to $\bar{\mathcal{E}}$ and with covariance $G$, if it is a martingale and if, for any $v_1,\ldots,v_{2n}\in\mathfrak{h}$ and $t_1,\ldots,t_{2n}\in\mathbb{R}_+$, we have*

$$\bar{\mathcal{E}}\left(\prod_{i=1}^{2n} X_{t_i}(v_i)\right)=\sum_{\mathcal{M}\in\text{Perfect matches of }\{1,\ldots,2n\}}(-1)^{\mathcal{M}}\prod_{(i,j)\in\mathcal{M}}\langle\Theta v_i,Gv_j\rangle(t_i\wedge t_j).$$

Obviously the $n$-point functions of Brownian motion defined with respect a signed state $\bar{\mathcal{E}}$ are equal to the ones of Brownian motion with respect to a positive state, as defined in Definition 3.14.

**Remark 4.40.** Hereafter we consider the following setting: we fix a filtered modular space $(\mathscr{M},\omega,(\mathscr{M}_t)_t)$, for which there is an analytic Brownian motion $X_t$, defined on $\mathfrak{h}$, and the Grassmann algebra $\mathcal{G}_X$ and $\mathcal{G}_X^p$ generated by $X$ and some compatible independent Grassmann algebra $\tilde{\mathscr{M}}_0\subset\mathscr{M}_0$ (see Definition 4.17). We then consider the space

$$\bar{\mathcal{G}}_X:=\bigcap_{p\geqslant 1}\mathcal{G}_X^p, \tag{4.31}$$

which is, thanks to Hölder's inequality, a Fréchet algebra with the seminorms $\|\cdot\|_{\mathbb{L}^p}$. Furthermore $\bar{\mathcal{G}}_X$ is Grassmann and the even and odds elements are $\bar{\mathcal{G}}_{X,\pm}=\bigcap_{p\geqslant 1}\mathcal{G}_{X,\pm}^p$. We call the case where $\mathcal{A}=\bar{\mathcal{G}}_X$ the *standard setting*.

In Section 2.4 we extended the conditional expectation $\omega_t$ to the $\mathbb{L}^p$ spaces, see Remark 2.23, so that the the state $\omega$ is also a signed expectation with conditional expectation $\omega_t$ on $\bar{\mathcal{G}}_X$ in the sense of Definition 4.37. We shall now describe a standard way of building from $\omega$ and $\omega_t$ a signed expectation $\bar{\mathcal{E}}$ with a conditional expectation $\bar{\mathcal{E}}_t$ on the Fréchet algebra $\bar{\mathcal{G}}_X$.

**Lemma 4.41.** *In the standard setting of Remark 4.40, consider a martingale $\{Z_t\}_{t\in\mathbb{R}_+}\subset\bar{\mathcal{G}}_{X,+}$ such that $Z_0=1$, $Z_\infty=\lim_{t\to+\infty}Z_t$ exists, $\omega_t(Z_\infty)=Z_t$, and, for any $t\in\mathbb{R}_+$ there is $Z_t^{-1}\in\bar{\mathcal{G}}_{X,+}$ such that $Z_t\cdot Z_t^{-1}=1$. Then, if for every $a\in\bar{\mathcal{G}}_X$ and $t\in\mathbb{R}_+$, we define*

$$\bar{\mathcal{E}}^Z(a)=\omega(aZ_\infty),\qquad \bar{\mathcal{E}}_t^Z(a)=\omega_t(aZ_\infty)Z_t^{-1}. \tag{4.32}$$

*Then $\bar{\mathcal{E}}^Z$ is an signed expectation on $\bar{\mathcal{G}}_X$ with associated conditional expectation $\bar{\mathcal{E}}_t^Z$ with respect to the filtration $\{\bar{\mathcal{G}}_{X,t}\}_{t\in\mathbb{R}_+}:=\{\bar{\mathcal{G}}_X\cap\mathbb{L}^1(\mathscr{M}_t)\}_{t\in\mathbb{R}_+}$ of $\bar{\mathcal{G}}_X$, in the sense of Definition 4.37.*

**Proof.** The continuity of $\bar{\mathcal{E}}^Z$ and $\bar{\mathcal{E}}_t^Z$ is given by Hölder's inequality for the twisted $\mathbb{L}^p$ spaces. Furthermore properties 1, 2, 3 of Definition 4.37 hold. Indeed, we have that, for any $a\in\bar{\mathcal{G}}_X$ and $s\leqslant t$,

$$\bar{\mathcal{E}}_s^Z(\bar{\mathcal{E}}_t^Z(a))=\omega_s(\omega_t(aZ_\infty)Z_t^{-1}Z_\infty)Z_s^{-1}=\omega_s(\omega_t(aZ_\infty)Z_t^{-1}\omega_t(Z_\infty))Z_s^{-1}=\omega_s(\omega_t(aZ_\infty))Z_s^{-1}=\bar{\mathcal{E}}_s^Z(a),$$

and similarly in the case where $\bar{\mathcal{E}}_s^Z$ is replaced by $\bar{\mathcal{E}}^Z$. Furthermore for any $b,c\in\bar{\mathcal{G}}_{X,t}$ we get

$$\bar{\mathcal{E}}_t^Z(bac)=\omega_t(bacZ_\infty)Z_t^{-1}=\omega_t(baZ_\infty c)Z_t^{-1}=b\omega_t(aZ_\infty)cZ_t^{-1}=b\bar{\mathcal{E}}_t^Z(a)c$$

where we used that $Z_\infty,Z_t^{-1}$ commute with any element of $\bar{\mathcal{G}}_X$ being both elements of $\bar{\mathcal{G}}_{X,+}$. □

Below, it will be important to characterize weak Grassmann Brownian motion. To this end, we will use the following theorem, which is a non-commutative version of Lévy's characterization theorem.

**Theorem 4.42.** *(Non-commutative Lévy's characterization). Suppose we are in the standard setting, see Remark 4.40, and let $\bar{\mathcal{E}}$ be a signed state on $\bar{\mathcal{G}}_X$ and admitting the conditional expectation $\bar{\mathcal{E}}_t$. Consider the Itô process $B_t = \int_0^t \langle H_s', \mathrm{d}X_s\rangle + \int_0^t K_s\,\mathrm{d}s$, with $\|H'\|_{\mathbb{L}^p,|G|} \in L^\infty_{\mathrm{loc}}(\mathbb{R}_+)$, $K \in L^\infty_{\mathrm{loc}}(\mathbb{R}_+, \mathcal{G}_X^p)$, for every $p \geqslant 2$, which is a martingale with respect to $\bar{\mathcal{E}}_t$. Then, $B$ is a weak Brownian motion with covariance $G$ and with respect to the signed state $\bar{\mathcal{E}}$ if and only if $B_0 = 0$ and the quadratic variation $[B(v), B(v')]_t$ with respect to $\bar{\mathcal{E}}_t$ is $[B(v), B(v')]_t = \langle \Theta v, G v'\rangle\, t$.*

The proof of Theorem 4.42 requires the following technical lemma.

**Lemma 4.43.** *Let $Y_t, Y_t'$ be two Itô processes taking values in $\bar{\mathcal{G}}_X$ such that $Y \in L^2_{\mathrm{loc}}(\mathbb{R}_+, \mathcal{G}_X^{2p})$ and $Y_t' = \int_0^t \langle H_s', \mathrm{d}X_s\rangle + \int_0^t K_s' \mathrm{d}s$ with $\|H_s'\|_{\mathbb{L}^{2q},|G|} \in L^2_{\mathrm{loc}}(\mathbb{R}_+)$, $K \in L^2(\mathbb{R}_+, \mathcal{G}_X^{2q})$ and $\frac{1}{p} + \frac{1}{q} = 1$. Then, for each $t \in \mathbb{R}_+$, the following identity holds in $\mathbb{L}^2$*

$$\lim_{|\pi| \to 0} \sum_{t_i \in \pi} Y_{t_{i-1}\wedge t}(Y'_{t_i \wedge t} - Y'_{t_{i-1}\wedge t}) = \int_0^t Y_s\,\mathrm{d}Y_s'.$$

**Proof.** The proof is similar to the one of Theorem 4.33 (see in particular equation (4.33)). □

An important consequence of Lemma 4.43 is the following.

**Corollary 4.44.** *Let $Y_t, Y_t'$ be two Itô processes satisfying the hypotheses of Lemma 4.43. Then, if $Y_t'$ is a martingale with respect to the signed expectation $\bar{\mathcal{E}}$, also $\int_0^t Y_s \mathrm{d}Y_s'$ is a martingale with respect to $\bar{\mathcal{E}}$. Furthermore, if $U, U'$ is an other pair of Itô processes satisfying the hypotheses of Lemma 4.43, such that $U_t' = \int_0^t \langle T_s', \mathrm{d}X_s\rangle + \int_0^t L_s\,\mathrm{d}s$ is a martingale with respect to the signed expectation $\bar{\mathcal{E}}$, then the quadratic variation of $\int_0^t Y_s \mathrm{d}Y_s'$ and $\int_0^t U_s \mathrm{d}U_s'$ with respect to $\bar{\mathcal{E}}$ is*

$$\left[\int_0^\cdot Y_s \mathrm{d}Y_s', \int_0^\cdot W_s \mathrm{d}W_s'\right]_t = \int_0^t \mathrm{Tr}_{G^*\Theta}((Y_s \cdot H_s') \otimes (U_s \cdot T_s'))\,\mathrm{d}s.$$

**Proof.** The proof of the first part is a consequence of the fact that $\sum_{t_i \in \pi} Y_{t_{i-1}\wedge t}(Y'_{t_i\wedge t} - Y'_{t_{i-1}\wedge t})$ is a $\bar{\mathcal{E}}$ martingale whenever $Y'$ is a $\bar{\mathcal{E}}$ martingale, and the fact that the martingale property is preserved by taking limits in $\mathbb{L}^2$, $\bar{\mathcal{E}}_t$ being continuous. The second part of the statement is a consequence of the first part and the Itô formula. □

**Proof of Theorem 4.42.** Let $\theta_1, \ldots, \theta_n \in \bar{\mathcal{G}}_{X,-,0}$ be independent and compatible random variables which are (non-zero) elements of $\tilde{\mathscr{M}}_0$ (by Lemma 3.10 it is always possible to assume that $\tilde{\mathscr{M}}_0$ contains a numerable set of independent and compatible Grassmann random variables). Consider $f_1, \ldots, f_n \colon \mathbb{R}_+ \to \mathfrak{h}$ smooth compactly supported functions. Let $\mathcal{A} \subset \mathcal{N} = \{1, \ldots, n\}$ and define

$$\mathcal{S}_{\mathcal{A}}(t) = \bar{\mathcal{E}}_0\left(\prod_{i \in \mathcal{A}} \left(1 + \theta_i \int_0^t \langle f_i(s), \mathrm{d}B_s\rangle\right)\right).$$

The functions $\mathcal{S}_{\mathcal{A}}(t)$ take value in the finite-dimensional Grassmann algebra $\mathfrak{G} = \mathrm{span}\{\prod_{j \in \mathcal{A}} \theta_j, \mathcal{A} \subset \mathcal{N}\}$. Furthermore, the coefficients of the function $\mathcal{S}_{\mathcal{A}}$ contain all the possible $n$-point functions of the Grassmann random variables $\int_0^t \langle f_i(s), \mathrm{d}B_s\rangle$, indeed

$$\mathcal{S}_{\mathcal{A}}(t) = \sum_{\mathcal{B} \subset \mathcal{A}} \mathbb{S}_{\mathcal{B}}(t)\,\theta_{\mathcal{B}}$$

where $\theta_{\mathcal{B}} = \prod_{\ell, \{i_1 < \ldots < i_k\} = \mathcal{B}} \theta_{i_\ell}$ and

$$\mathbb{S}_{\mathcal{B}}(t) = (-1)^{\mathcal{B}} \bar{\mathcal{E}}\left(\prod_{\ell, \{i_1 < \ldots < i_k\} = \mathcal{B}} \int_0^t \langle f_{i_\ell}(s), \mathrm{d}B_s\rangle\right).$$

This means that, if for any $\mathcal{A} \subset \mathcal{N}$, the function $\mathcal{S}_{\mathcal{A}}(t)$ is equal to

$$\bar{\mathcal{S}}_{\mathcal{A}}(t) = \mathcal{E}_0\left(\prod_{i\in\mathcal{A}} \left(1 + \theta_i \int_0^t \langle f_i(s), \mathrm{d}X_s\rangle\right)\right)$$

where $X$ is some Brownian motion with covariance $G$ with respect to $\omega$, then $B$ has the same law of a Brownian motion.

Suppose that $B$ is a $\bar{\mathcal{E}}$ martingale and that it has quadratic variation with respect to $\bar{\mathcal{E}}$ equal to $\langle G^*\Theta\cdot,\cdot\rangle t$. By Corollary 4.44 also $\int_0^t \langle f_j(s)\mathrm{d}B_s\rangle$ are martingales with quadratic variation given by $\int_0^t \langle G^*\Theta f_i(s), f_j(s)\rangle \mathrm{d}s$. This means that, if we apply the Itô formula to the product $\prod_{i\in\mathcal{A}} (1+\theta_i\int_0^t \langle f_i(s), \mathrm{d}B_s\rangle)$, we get that $\mathcal{S}_{\mathcal{A}}$ satisfies the following system of ODE

$$\mathcal{S}_{\mathcal{A}}(t) = 1 + \sum_{\{i,j\}\subset\mathcal{A}} \int_0^t \langle G^*\Theta f_i(s), f_j(s)\rangle \theta_j\theta_i \mathcal{S}_{\mathcal{A}\setminus\{i,j\}}(s)\,\mathrm{d}s.$$

Since the previous system of ODE for $\{\mathcal{S}_{\mathcal{A}}\}_{\mathcal{A}\subset\mathcal{N}}$ is takes value in the finite-dimensional (Grassmann) algebra $\mathfrak{G}$ it has a unique solution. On the other hand, using again the Itô formula, the previous system of ODEs is satisfied by the $n$-point functions of the Brownian motion $X$ with covariance $G$, namely

$$\bar{\mathcal{S}}_{\mathcal{A}}(t) = 1 + \sum_{\{i,j\}\subset\mathcal{A}} \int_0^t \langle G^*\Theta f_i(s), f_j(s)\rangle \theta_j\theta_i \bar{\mathcal{S}}_{\mathcal{A}\setminus\{i,j\}}(s)\,\mathrm{d}s.$$

This implies that $\{\mathcal{S}_{\mathcal{A}}(t)\}_{\mathcal{A}\subset\mathcal{N}}$ (where $T\in\mathbb{R}_+$ is such that $\mathrm{supp}(f_j)\subset[0,T]$) are equal to $\{\bar{\mathcal{S}}_{\mathcal{A}}(t)\}_{\mathcal{A}\subset\mathcal{N}}$, which is equivalent to saying that $B_t$ has, with respect to $\bar{\omega}$, the law of the Brownian motion $X$ with covariance $G$. The opposite direction of the implication is trivial. □

Theorem 4.42 is crucial for proving the following non-commutative version of Girsanov's theorem.

**Theorem 4.45.** *(Non-commutative Girsanov's formula). Let $(\mathscr{M}, \omega, (\mathscr{M}_t)_t)$ be a filtered modular space and let $X$ be an analytic Brownian motion of covariance $G$. Consider a process $H$: $\mathbb{R}_+\times\mathfrak{h}\to\bar{\mathcal{G}}_{X,-}$, where $\bar{\mathcal{G}}_X$ is the Fréchet algebra introduced in Remark 4.40, and let $Z_t$ be an even adapted process in $\mathbb{L}^2$ satisfying*

$$Z_t = 1 + \int_0^t \langle Z_s\cdot H_s, \mathrm{d}X_s\rangle. \tag{4.33}$$

*Let $\bar{\mathcal{E}}^Z$ be the expectation associated with it. Then, the random field, such that*

$$B_t(v) = X_t(v) - \int_0^t H_s(G^*\Theta(v))\mathrm{d}s, \qquad \forall v\in\mathfrak{h} \quad t\in\mathbb{R}_+,$$

*is a Brownian motion of covariance $G$ with respect to the signed expectation $\bar{\mathcal{E}}^Z$.*

**Remark 4.46.** By requiring that $Z\in\mathbb{L}^2$ satisfies equation (4.33), we automatically get that $Z$ is a martingale by definition of the Itô integral, see Section 4.1.

**Proof Theorem 4.45.** The proof follows the same strategy of the commutative case (see, e.g., [Le 16, Section 5.5]). First of all, by (4.32) for the conditional expectation $\bar{\mathcal{E}}_t^Z$, we have that the process $B_t$ is a martingale with respect to $\bar{\mathcal{E}}_t^Z$ if and only if the process $B_tZ_t$ is a martingale with respect to the conditional expectation $\omega_t$. On the other hand, by the Itô formula, the fact that $Z$ commutes with $B$ and (4.33), we get

$$\begin{aligned}
B_t(v)\,Z_t &= \int_0^t Z_s\mathrm{d}B_s(v) + \int_0^t B_s(v)\,\mathrm{d}Z_s + [B(v), Z]_t \\
&= \int_0^t Z_s\mathrm{d}X_s(v) - \int_0^t Z_sH_s\,(G^*\Theta\,v)\,\mathrm{d}s + \int_0^t \langle B_s(v)\,Z_s\cdot H_s, \mathrm{d}X_s\rangle + \int_0^t Z_sH_s(G^*\Theta\,G\,v)\,\mathrm{d}s \\
&= \int_0^t Z_s\mathrm{d}X_s + \int_0^t \langle B_s(v)\,Z_sH_s, \mathrm{d}X_s\rangle
\end{aligned}$$

which is a martingale with respect to $\omega_t$, since it is an Itô integral with respect to $X$. On the other hand, by Corollary 4.44, the quadratic variation of $B_t$ with respect to $\bar{\mathcal{E}}_t^Z$ is

$$[B(v),B(v')]_t=[X(v),X(v')]_t=\langle\Theta\, v,G v'\rangle.$$

Thus since $B$ is a martingale with respect to $\bar{\mathcal{E}}_t^Z$ with quadratic variation $\langle\Theta\, v,G v'\rangle$, by Theorem 4.42, $B$ is a weak Brownian motion with respect to $\bar{\mathcal{E}}_t^Z$ with covariance $G$. □

Unlike the commutative case, we do not know if (4.33) has a solution because the exponential of an Itô process in $\bar{\mathcal{G}}_X$ is generally not well-defined. We will now consider some sufficient conditions on $H$ such that the solution $Z$ to (4.33) exists.

**Definition 4.47.** *Let $X$ be an analytic Brownian motion of covariance $G$. An adapted process $H_.\colon \mathbb{R}_+\times\mathfrak{h}\to\mathcal{G}_-^X$ satisfies the Novikov condition if for every $p\geqslant 2$ we have $\|H_.\|_{\mathbb{L}^p,|G|}\in L^\infty_{\mathrm{loc}}(\mathbb{R}_+)$, and for every $t>0$, we get*

$$\sum_{n=0}^{+\infty}\frac{1}{n!}\left\|\left(\int_0^t\langle H_s,\mathrm{d}X_s\rangle-\frac{1}{2}\int_0^t\mathrm{Tr}_{G^*\Theta}(H_s\otimes H_s)\mathrm{d}s\right)^n\right\|_{\mathbb{L}^p}<\infty.$$

Let $H_.\colon\mathbb{R}_+\times\mathfrak{h}\to\mathcal{G}_{X,-}$ be an adapted process satisfying the Novikov condition. We set

$$Z_t^{(n)}:=\sum_{j=0}^{n}\frac{1}{j!}\left(\int_0^t\langle H_s,\mathrm{d}X_s\rangle-\frac{1}{2}\int_0^t\mathrm{Tr}_{G^*\Theta}(H_s\otimes H_s)\mathrm{d}s\right)^j.$$

For every $p\geqslant 1$, $Z_t^{(n)}$ converges to

$$Z_t=\sum_{j=0}^{\infty}\frac{1}{j!}\left(\int_0^t\langle H_s,\mathrm{d}X_s\rangle-\frac{1}{2}\int_0^t\mathrm{Tr}_{G^*\Theta}(H_s\otimes H_s)\mathrm{d}s\right)^j \tag{4.34}$$

in $C^0(\mathbb{R}_+,\mathbb{L}^p)$, and $Z_t$ is an invertible element of $\mathcal{G}_+^X$ with $Z_t^{-1}$ given by $Z_t^{-1}=\sum_{j=0}^\infty\frac{(-1)^j}{j!}\big(\int_0^t\langle H_s,\mathrm{d}B_s\rangle-\frac{1}{2}\int_0^t\mathrm{Tr}_{G^*\Theta}(H_s\otimes H_s)\mathrm{d}s\big)^j$.

**Lemma 4.48.** *Let $X$ be an analytic Brownian motion of covariance $G$ and let $H_.\colon\mathbb{R}_+\times\mathfrak{h}\to\mathcal{G}_-^X$ be an adapted process satisfying the Novikov condition. Then, $Z_t$ defined in (4.34) satisfies (4.33).*

**Proof.** By the Itô formula we get that

$$Z_t^{(n)}=1+\int_0^t\langle Z_s^{(n-1)}\cdot H_s,\mathrm{d}X_s\rangle+\frac{1}{2}\int_0^t(Z_s^{(n-2)}-Z_s^{(n-1)})\cdot\mathrm{Tr}_{G^*\Theta}(H_s\otimes H_s)\mathrm{d}s.$$

Taking the limit $n\to\infty$, which is well-defined by the Novikov condition, we obtain the claim. □

When $H$ takes values in $\mathcal{G}_X^\infty$ the Novikov condition is automatically satisfied. More precisely, we have the following proposition.

**Proposition 4.49.** *Let $X$ be an analytic Brownian motion of covariance $G$ and let $H_.\colon\mathbb{R}_+\times\mathfrak{h}\to\mathcal{G}_-^X$ be an adapted process such that $\mathrm{Tr}_{G^*\Theta}(H_.\otimes H_.)\in L^\infty(\mathbb{R}_+,\mathcal{G}_X^\infty)$. Then, $H$ satisfies the Novikov condition.*

**Proof.** We note that for any $n\in\mathbb{N}$

$$\frac{1}{n!}\left(\int_0^t\langle H_s,\mathrm{d}X_s\rangle-\frac{1}{2}\int_0^t\mathrm{Tr}_{G^*\Theta}(H_s\otimes H_s)\mathrm{d}s\right)^n=Z_t^{(n)}-Z_t^{(n-1)}. \tag{4.35}$$

By Theorem 4.5, equation (4.48) and Hölder's inequality for the twisted $\mathbb{L}^p$ spaces we get

$$\begin{aligned}&\|Z_t^{(n)}-Z_t^{(n+1)}\|_{\mathbb{L}^p}\\&\leqslant C_p\int_0^t \|Z_s^{(n-1)}-Z^{(n)}\|_{\mathbb{L}^p}(\|\mathrm{Tr}_{G^*\Theta}(H_s\otimes H_s)\|_{\mathbb{L}^\infty})^{\frac{1}{2}}\mathrm{d}s\\&\quad+\frac{C_p}{2}\int_0^t (\|(Z_s^{(n-2)}-Z^{(n-1)})\|_{\mathbb{L}^p}+\|Z_s^{(n-1)}-Z^{(n)}\|_{\mathbb{L}^p})\|\mathrm{Tr}_{G^*\Theta}(H_s\otimes H_s)\|_{\mathbb{L}^\infty}\mathrm{d}s\Big)\end{aligned} \tag{4.36}$$

for a suitable constant $C_p>0$ (depending on $p\geqslant 2$). By induction, it is simple to see that

$$\sup_{t\in[0,T]}\|Z_t^{(n)}-Z_t^{(n+1)}\|_{\mathbb{L}^p}\leqslant\frac{(3C_p(1+\|\mathrm{Tr}_{G^*\Theta}(H_s\otimes H_s)\|_{\mathbb{L}^\infty})\,(T+1))^n}{n!},$$

hence the claim. □

## 4.4 Weak solutions of finite-dimensional Grassmann SDE

In this section we define the notion of weak and strong solution to an Grassmann stochastic differential equation (SDE), when the vector space $\mathfrak{h}$ is finite dimensional. First of all we introduce some notation. We note that the space $\mathcal{G}_X^\infty=\bar{\mathcal{G}}_X\cap\mathbb{L}^\infty$ is a Banach algebra, thus if an adapted process $B$ belongs to $C^0(\mathbb{R}_+,\mathcal{G}_X^\infty)$ we can consider the Banach algebra $\mathcal{G}_B^\infty$ as the natural Banach subalgebra generated by $B$ in $\mathcal{G}_X^\infty$.

We now consider two linear functions

$$\mu:\mathfrak{h}\to\bigoplus_{n=0}^{+\infty}\Lambda^{2n+1}\mathfrak{h},\quad \sigma:\mathfrak{h}\times\mathfrak{h}\to\left(\bigoplus_{n=0}^{+\infty}\Lambda^{2n}\mathfrak{h}\right). \tag{4.37}$$

**Definition 4.50.** *In a standard setting, let $\bar{\mathcal{E}}$ a signed expectation on $\bar{\mathcal{G}}_X$ and consider two linear functions $\mu,\sigma$ as in (4.37), we say that the couple of adapted Itô odd random fields $(\Psi_\cdot,B_\cdot)$ is a weak solution to the SDE $(\mu,\sigma)$ (with respect the signed expectation $\bar{\mathcal{E}}$) on $[0,T]$ and with initial condition $\Psi_0$ if $B$ is a weak Brownian motion with covariance $G$ with respect to $\bar{\mathcal{E}}$ and, for any $v\in\mathfrak{h}$ and $t\leqslant T$, the following equation holds:*

$$\Psi_t(v)-\Psi_0(v)=\int_0^t\mu(v)(\Psi_s)\mathrm{d}s+\int_0^t\langle\sigma(v,\cdot)(\Psi_s),\mathrm{d}B_s\rangle.$$

*Furthermore, if $B\in C^0(\mathbb{R}_+,\mathcal{G}_X^\infty)$ we say that $\Psi$ is a strong solution to the SDE $(\mu,\sigma)$ driven by $B$ (till the time $T$) if $(\Psi,B)$ is a weak solution to the SDE $(\mu,\sigma)$ and $\Psi_t\in\mathcal{G}_B^\infty$ for any $0\leqslant t\leqslant T$.*

In the special case where $\sigma(v,v')=\langle v,v'\rangle_{\mathfrak{h}}$ (i.e. *additive noise SDEs*) and $\mathfrak{h}$ is finite-dimensional, we can prove the existence of strong solutions.

**Theorem 4.51.** *Suppose that $\mathfrak{h}$ is finite-dimensional and let $\sigma(v,v')=\langle v,v'\rangle_{\mathfrak{h}}$ then for any $\Psi_0\in\mathcal{G}_{X,-}^\infty=\mathcal{G}_X^\infty\cap\bar{\mathcal{G}}_{X,-}$, $B\in C^0(\mathbb{R}_+,\mathcal{G}_X^\infty)$ Brownian motion with respect to $\bar{\mathcal{E}}$, and $T\geqslant 0$, there is a unique strong solution $\Psi_t\in C^0(\mathbb{R}_+,\mathcal{G}_B^\infty)$ to the SDE $(\mu,\sigma)$, driven by $B$.*

**Proof.** This theorem, using a slightly difference language and notations, has been proved in [ABDG22, Theorem 30 and 31], in the case where $\|B_t-B_s\|_{\mathbb{L}^\infty}\lesssim|t-s|^{\frac{1}{2}}$. The case where the map $t\mapsto B_t$ is continuous in $\mathbb{L}^\infty$ can be proved in a similar way. □

Here we introduce the notation

$$\mu^A(v):=\mu(v)-A\,v,$$

where $A:\mathfrak{h}\to\mathfrak{h}$ is a linear map.

**Theorem 4.52.** *Suppose that $\mathfrak{h}$ is finite-dimensional, suppose that $\bar{G} = G^* \Theta$ is an invertible anti-linear map, and consider $\sigma(v, v') = \langle v, v' \rangle_{\mathfrak{h}}$. Write*

$$X_t^A(v) = \tilde{X}_0(\mathrm{e}^{At} v) + \int_0^t \langle \mathrm{e}^{A(t-s)}(v), \mathrm{d}X_s \rangle.$$

*and consider the processes*

$$\begin{aligned} Z_t^{\mu^A} &= \exp\Big( \int_0^t \mathbf{1}_{[0,T]}(s) \langle \mu^A(\bar{G}^{-1} \cdot, (X_s^A + \mathrm{e}^{As} h_0)), \mathrm{d}X_s \rangle + \\ &\quad -\frac{1}{2} \int_0^t \mathbf{1}_{[0,T]}(s) \mathrm{Tr}_{\bar{G}}(\mu^A(\bar{G}^{-1} \cdot, (X_s^A + \mathrm{e}^{As} h_0)) \otimes \mu^A(\bar{G}^{-1} \cdot, (X_s^A + \mathrm{e}^{As} h_0))) \mathrm{d}s \Big) \\ B_t(\cdot) &= \int_0^t \mathbf{1}_{[0,T]}(s) \mu^A(\cdot, (X_s^A + h_0(\mathrm{e}^{As} \cdot))) \mathrm{d}s + X_t. \end{aligned}$$

*The the pair of processes $(X^A + h_0(\mathrm{e}^{A \cdot}), B)$ is a weak solution of the SDE $(\mu, \sigma)$ on $[0, T]$ with initial condition $\tilde{X}_0 + h_0$ and with respect to the expectation $\bar{\mathcal{E}}^{Z^{\mu^A}}$. Furthermore $X_t^A + h_0(\mathrm{e}^{At})$ coincides with the strong solution of Theorem 4.51.*

**Proof.** First we note that the process $Z_t^{\mu^A}$ is well-defined. Indeed, $\mu^A(\bar{G}^{-1} \cdot, (X_s^A + \mathrm{e}^{As} h_0)) \in L^\infty_{\mathrm{loc}}(\mathbb{R}_+, \mathcal{G}_X^\infty)$ and by Proposition 4.49 it satisfies the Novikov condition. Furthermore, the process $X_t^A + h_0(\mathrm{e}^{At})$ solves the SDE $X_t^A(v) + h_0(\mathrm{e}^{At} v) - \tilde{X}_0(v) - h_0(v) = \int_0^t (X_s^A(Av) + h_0(\mathrm{e}^{At} A v)) \mathrm{d}s + X_t(v)$. Which means that we have the following relation, for any $t \leqslant T$

$$X_t^A(v) + h_0(\mathrm{e}^{At} v) - \tilde{X}_0(v) - h_0(v) = \int_0^t \mu(v, X_s^A + h_0(\mathrm{e}^{As} \cdot)) \mathrm{d}s + B_t(v).$$

Since, by Theorem 4.45, $B_t$ is a Brownian motion with respect to the expectation $\bar{\mathcal{E}}^{Z_t^{\mu^A}}$, equation (4.52) implies that $(X^A + h_0(\mathrm{e}^{A \cdot}), B)$ is a weak solution to the SDE $(\mu, \sigma)$. Finally, since by Theorem 4.51 the strong solution exists and is unique, we have $\mathcal{G}_X^\infty = \mathcal{G}_B^\infty$ so that any weak solution is also a strong solution. □

# 5 Applications

## 5.1 Twisted $L^p$ bound of the exponential

In Subsection 5.1 we prove the existence of a family of weak Gibbsian quartic perturbation of $\omega$, a special case being the Grassmann $\Psi_2^4$ on the torus. The brevity of the arguments presented should convince the reader about the usefulness of the twisted $\mathbb{L}^p$ spaces introduced in Section 2.3. Our strategy is based on the control of the exponential of a Wick's polynomial in the spirit of the well-known bound by Nelson [Nel66, Sim74].

In order to consider spin-$\frac{1}{2}$ fermionic fields on $\mathbb{T}^2$, we let $\mathfrak{h} = L^2(\mathbb{T}^2; \mathbb{C}^4)$ and let $(\mathscr{M}^{(\mu)}, \omega, (\mathscr{M}_t^{(\mu)})_{t \geqslant 0})$ be the filtered modular space associated with $\mathfrak{h}$, see Section 3.2. For the sake of brevity, we henceforth fix $\mu = \frac{1}{2}$ and remove it from our notation. We let $(\Psi_t)_t$ be a GBM with covariance

$$G_t = (1 - \Delta)^{\theta - 1} \chi(|-\Delta| / t) (\mathbb{1} \oplus -\mathbb{1}), \tag{5.1}$$

where $-\Delta$ is the Laplacian on $h$, where $\theta \geqslant 0$ and where $\chi \in C_c^\infty(\mathbb{R}^+)$ is such that $\chi([0, 1/2]) = 1$ and $\mathrm{supp}(\chi) = [0, 1]$. We let $\delta_{x,t}(y) = \sum_{k \in \mathbb{Z}^2} \mathrm{e}^{ik(y-x)} \chi(|k| / t)$ together with $\delta_{x,t,+} = (\delta_{x,t}, 0)$ and $\delta_{x,t,-} = (0, \delta_{x,t})$, and introduce

$$\bar{\Psi}_{t,\sigma}(x) := \Psi_t(\delta_{x,t,\sigma} \oplus 0), \qquad \Psi_{t,\sigma}(x) := \Psi_t(0 \oplus \delta_{x,t,\sigma}) \qquad \sigma \in \{\pm\}.$$

Note that

$$\omega(\Psi_t(x)\bar{\Psi}_t(y)) = G_t(x;y) = \sum_{k\in\mathbb{Z}^2} \frac{\mathrm{e}^{\mathrm{i}k(x-y)}\chi(|k|/t)}{(1+k^2)^{1-\theta}} 1, \tag{5.2}$$

and that by Proposition 3.17, if $\theta > 0$

$$\sup_{x,\sigma} \|\Psi_{t,\sigma}(x)\|^2 \lesssim \sum_{k\in\mathbb{Z}^2} \frac{\chi(|k|/t)}{(1+k^2)^{1-\theta}} \sim t^{2\theta} \tag{5.3}$$

whereas the divergence is as $\sim \log t$ if $\theta = 0$. We prove the existence of weak Gibbsian quartic perturbation of $\omega$, for $\theta$ small enough.

**Theorem 5.1.** *For any $\lambda \in \mathbb{C}$ and $t \geqslant 0$ define $Z_t^{(\lambda)} = \exp(\lambda V_t)$ with*

$$V_t := \int_{\mathbb{T}^2} [\![(\bar{\Psi}_t(x), \Psi_t(x))^2]\!] \mathrm{d}x,$$

*having set $(\bar{\Psi}_t(x), \Psi_t(x)) = \sum_{\sigma=\pm} \bar{\Psi}_{t,\sigma}(x)\, \Psi_{t,\sigma}(x)$. Then, if $7\theta < 1$ and if $|\lambda|$ is small enough depending on $\mu$, the weak limit*

$$\omega^{(\lambda)}(\cdot) := \lim_{t\to\infty} \frac{\omega(\cdot Z_t^{(\lambda)})}{\omega(Z_t^{(\lambda)})}, \tag{5.4}$$

*is a well-defined continuous normalized linear functional on $\mathbb{L}^p$ for any $1 < p \leqslant \infty$.*

**Remark 5.2.** In particular, $\omega^{(\lambda)}$ is a signed expectation on the Fréchet algebra $\bar{\mathcal{G}}_\Psi$, see Definition 4.37 and Remark 4.40. The fermionic $\Psi_2^4$ theory corresponds to the choice $\theta = 0$. Note that the analysis we present can be easily generalized to the case of a general polynomial interaction and of any dimension, provided that the covariance is suitably adjusted.

The claim in Theorem 5.1 is a straightforward consequence of Theorem 5.3 and Corollary 5.5 below. The crucial point here is to show that even though $\lim_{t\to\infty} Z_t^{(\lambda)}$ is not a bounded random variable, it belongs to the twisted $\mathbb{L}^p$ spaces for any $1 \leqslant p < \infty$.

**Theorem 5.3.** *Let $7\theta < 1$. Then, for any $\lambda \in \mathbb{C}$ we have $Z_t^{(\lambda)} \in \bigcap_{p<\infty} \mathbb{L}^p$ uniformly in $t \geqslant 0$ and $Z_t^{(\lambda)} \to Z^{(\lambda)}$ as $t \to \infty$ in the $\mathbb{L}^p$ topology, for any $1 \leqslant p < \infty$. Furthermore, for any $x \in \mathbb{L}^p$, $1 < p \leqslant \infty$, $\omega(x \cdot Z_t^{(\lambda)})$ has a limit as $t \to \infty$.*

The following technical lemma will be used in the proof of Theorem 5.3.

**Lemma 5.4.** *The following bounds hold true for any $\nu < 1 - 3\theta$*

$$\|V_t\|_{\mathbb{L}^\infty} \lesssim t^{4\theta}, \qquad \|T_\tau^{(2)}(V_t - V_s)\|_2^2 \lesssim_\tau s^{-2\nu},$$

*where the constant in $\lesssim_\tau$ is locally bounded in $\tau$.*

**Proof.** To prove the first bound, we write the Wick's polynomial explicitly

$$[\![(\bar{\Psi}_t(x), \Psi_t(x))^2]\!] = (\bar{\Psi}_t(x), \Psi_t(x))^2 + 2C_t(\bar{\Psi}_t(x), \Psi_t(x)) + 2C_t^2$$

with $C_t := \omega(\Psi_{t,\sigma}(x)\bar{\Psi}_{t,\sigma}(x)) \sim t^{2\theta}$, see (5.2), and therefore, by (5.3) and by Lemma 3.19 we obtain the bound

$$\|V_t\|_{\mathbb{L}^\infty} \lesssim \sup_{x\in\mathbb{T}^2} \|[\![(\bar{\Psi}_t(x), \Psi_t(x))^2]\!]\| \lesssim t^{4\theta}.$$

To prove the other bound, we follow the strategy of [Seg86], see also [Sim74], and switch to the Fourier components:

$$V_t = \sum_{k_1,\ldots,k_4\in\mathbb{Z}^2} \mathbf{1}_{\sum_i k_i=0}\left(\prod_i \chi(|k_i|/t)\right) [\![ (\hat{\hat{\Psi}}_t(k_1), \hat{\Psi}_t(k_2))(\hat{\hat{\Psi}}_t(k_3), \hat{\Psi}_t(k_4)) ]\!],$$

where now, denoting spinorial plane waves by $e^+_{k,\sigma} := e_{k,\sigma}\oplus 0$, $e^-_{k,\sigma} := 0\oplus e_{k,\sigma}$, with $e_{k,+}(y) = (\mathrm{e}^{\mathrm{i}ky}, 0)$ and $e_{k,-}(y) = (0, \mathrm{e}^{\mathrm{i}ky})$, we set $\hat{\Psi}_{t,\sigma}(k) := \Psi_t(e^-_{k,\sigma})$ and $\hat{\hat{\Psi}}_{t,\sigma}(k) = \Psi_t(e^+_{k,\sigma})$. With the notation introduced in Lemma 3.15, we also set $e^\ell_{k,\sigma,t} := \mathscr{C}_t e^\ell_{k,\sigma}$ together with $\tilde{e}^\ell_{k,\sigma,t} := \tilde{\mathscr{C}}_t e^\ell_{k,\sigma}$ and note that, for some constant $C$

$$\langle e^\ell_{k,\sigma,t}, e^{\ell'}_{k',\sigma',s}\rangle_{\mathcal{H}} = \langle \tilde{e}^\ell_{k,\sigma,t}, \tilde{e}^{\ell'}_{k',\sigma',s}\rangle_{\mathcal{H}} = C\delta_{\ell,\ell'}\delta_{\sigma,\sigma'}\delta_{k,k'}\frac{\chi(|k|/(t\wedge s))}{(1+k^2)^{1-\theta}}. \tag{5.5}$$

Then, by Lemma 3.19, noting that the r.h.s. of (5.5) depends only on $r = t\wedge s$, we have

$$\begin{aligned}
&\langle T^{(2)}_\tau([\![\hat{\hat{\Psi}}_{t,\sigma}(k_1)\,\hat{\Psi}_{t,\sigma}(k_2)\,\hat{\hat{\Psi}}_{t,\rho}(k_3)\,\hat{\Psi}_{t,\rho}(k_4)]\!]),\\
&\qquad T^{(2)}_\tau([\![\hat{\hat{\Psi}}_{s,\sigma'}(k_1')\,\hat{\Psi}_{s,\sigma'}(k_2')\,\hat{\hat{\Psi}}_{s,\rho'}(k_3')\,\hat{\Psi}_{s,\rho'}(k_4')]\!])\rangle_{L^2}\\
&= C_\tau\Big\langle (e^+_{k_1,\sigma,r}\oplus \tilde{e}^+_{k_1,\sigma,r})\wedge\cdots\wedge(e^-_{k_4,\rho,r}\oplus\tilde{e}^-_{k_4,\rho,r}),\\
&\qquad (e^+_{k_1',\sigma',r}\oplus\tilde{e}^+_{k_1',\sigma',r})\wedge\cdots\wedge(e^-_{k_4',\rho',r}\oplus\tilde{e}^-_{k_4',\rho',r})\Big\rangle_{\Gamma_a(\mathcal{H}\oplus\mathcal{H})}
\end{aligned} \tag{5.6}$$

for some constant $C_\tau$ locally bounded in $\tau$. Therefore, combining (5.5) and (5.6), we have

$$\begin{aligned}
&\|T^{(2)}_\tau(V_t - V_s)\|_2^2 = \|T^{(2)}_\tau(V_t)\|_2^2 - \|T^{(2)}_\tau(V_s)\|_2^2\\
&= C_\tau \sum_{k_1,\ldots,k_4\in\mathbb{Z}^2}\left(\prod_i \chi(|k_i|/t) - \prod_i \chi(|k_i|/s)\right)\frac{\mathbf{1}_{\sum_i k_i=0}}{\prod_i (1+k_i^2)^{1-\theta}},
\end{aligned}$$

for some other constant $C_\tau$ locally bounded in $\tau$. By the Young's inequality for convolutions we obtain, for any $\varepsilon > 0$

$$\|T^{(2)}_\tau(V_t - V_s)\|_2^2 \lesssim_\tau \|f_{s,t}\|_{L^p(\mathbb{R}^2)}\|f_{0,t}\|^3_{L^q(\mathbb{R}^2)} \lesssim s^{-2\nu} \tag{5.7}$$

where $\lesssim_\tau$ is up to a constant locally bounded in $\tau$, where

$$q = \frac{1+\varepsilon}{1-\theta}, \qquad p^{-1} = 3\left(\frac{\varepsilon+\theta}{1+\varepsilon}\right), \qquad \nu = \frac{1-2\varepsilon-3\theta}{1+\varepsilon},$$

and where $f_{s,t}(k) = (\chi(|k_i|/t) - \chi(|k_i|/s))(1+k^2)^{\theta-1}$ is such that $\|f_{s,t}\|_{L^q(\mathbb{R}^2)} \lesssim 1$ uniformly in $s$, $t\in\mathbb{R}$ and $\varepsilon > 0$. This concludes the proof. □

**Proof of Theorem 5.3.** By Hölder's inequality for the twisted $\mathbb{L}^p$ spaces, we have

$$\|Z^{(\lambda)}_t\|_{\mathbb{L}^p} \leqslant \sum_{n\geqslant 0}\frac{|\lambda|^n\|V_t\|^n_{\mathbb{L}^{pn}}}{n!}. \tag{5.8}$$

By Lemma 5.4 and by the hypercontractivity bounds, see Lemma 3.23, noting that $C_{\frac{p}{p-1},\varrho=\frac{1}{2}\mathbb{1}_{\mathcal{H}}} \lesssim p^{\frac{1}{2}}$, we have for $p\geqslant 2$

$$\|T^{(p)}_\tau(V_t - V_s)\|_p \lesssim_\tau p^2 s^{-\nu}. \tag{5.9}$$

Therefore, by Proposition 2.19, by Lemma 5.4 and by (5.9) we obtain for any $1\leqslant p<\infty$

$$\|V_t\|_{\mathbb{L}^p} \leqslant \|V_s\|_{\mathbb{L}^\infty} + \|V_t - V_s\|_{\mathbb{L}^{p\vee 2}} \lesssim s^{4\theta} + p^2 s^{-\nu}.$$

In particular, choosing $s=p^{\frac{2}{4\theta+\nu}}$ we have $\|V_t\|_{\mathbb{L}^p}\lesssim p^{\frac{8\theta}{4\theta+\nu}}$, which by (5.8) and by the fact that $7\theta<1$ implies

$$\|Z_t^{(\lambda)}\|_{\mathbb{L}^p}\leqslant\sum_{n\geqslant 0}\frac{|\lambda|^n c^n (pn)^{\frac{8\theta n}{4\theta+\nu}}}{n!}<\infty,$$

for some constant $c$ and for any $1\leqslant p<\infty$, hence $Z_t^{(\lambda)}\in\bigcap_{p<\infty}\mathbb{L}^p$ uniformly in $t\geqslant 0$. To prove that the limit exists, we use that $V_t$ are commuting random variables and by the fundamental theorem of calculus write

$$\mathrm{e}^{\lambda V_t}-\mathrm{e}^{\lambda V_s}=\lambda\int_0^1 \mathrm{e}^{\lambda V_t(1-r)}(V_t-V_s)\,\mathrm{e}^{\lambda V_s r}\mathrm{d}r.$$

Therefore, for any $1\leqslant p<\infty$ and any $s\leqslant t$

$$\|Z_t^{(\lambda)}-Z_s^{(\lambda)}\|_{\mathbb{L}^p}\leqslant|\lambda|\left(\sup_{r=s,t}\|Z_r^{(\lambda)}\|_{\mathbb{L}^{2p}}\right)\|V_t-V_s\|_{\mathbb{L}^{2p}}\lesssim_{\lambda,p}s^{-\nu},$$

proving the continuity in the $\mathbb{L}^p$ topology. The continuity of $\omega(x\cdot Z_t^{(\lambda)})$ for $x\in\mathbb{L}^p$, $1<p\leqslant\infty$ follows by the fact that $x\cdot Z_t^{(\lambda)}\in\mathbb{L}^1$ for any $t$ and that

$$|\omega(x\cdot Z_t^{(\lambda)})-\omega(x\cdot Z_s^{(\lambda)})|=|\omega(x\cdot(Z_t^{(\lambda)}-Z_s^{(\lambda)}))|\leqslant\|x\|_{\mathbb{L}^p}\|Z_t^{(\lambda)}-Z_t^{(\lambda)}\|_{\mathbb{L}^{p'}},$$

where $p'^{-1}=1-p^{-1}<\infty$. □

On the other hand, the normalization factor can be controlled for small $|\lambda|$.

**Corollary 5.5.** *If $|\lambda|$ is small enough, then $|1-\omega(Z_t^{(\lambda)})|\lesssim|\lambda|$ for any $t\geqslant 0$. Furthermore, the limit $\lim_{t\to\infty}\omega(Z_t^{(\lambda)})$ exists.*

**Proof.** Along the lines of the proof of Theorem 5.3, by using the properties of the twisted $\mathbb{L}^p$ spaces we have

$$|1-\omega(Z_t^{(\lambda)})|=|\omega(1-Z_t^{(\lambda)})|\leqslant\|1-Z_t^{(\lambda)}\|_{\mathbb{L}^1}\leqslant\sum_{n\geqslant 1}\frac{|\lambda|^n c^n n^{\frac{8\theta n}{4\theta+\nu}}}{n!}\lesssim_\theta|\lambda|.$$

Furthermore, for any $s\leqslant t$

$$|\omega(Z_t^{(\lambda)})-\omega(Z_s^{(\lambda)})|=|\omega(Z_t^{(\lambda)}-Z_s^{(\lambda)})|\leqslant\|Z_t^{(\lambda)}-Z_s^{(\lambda)}\|_{\mathbb{L}^1}\lesssim s^{-\nu},$$

hence the existence of the limit. □

## 5.2 Weak solution to $\Psi_2^4$ stochastic quantization SPDE

In this section we provide the solution to the stochastic quantization equation for a fermionic $\Psi_2^4$ on the two-dimensional torus $\mathbb{T}^2$ (namely the fermionic SPDE having the signed expectation built in Section 5.1 as invariant solution), in the same way in which Jona-Lasinio and Mitter [JM85] proved the existence of a weak solution in the bosonic case. More precisely, our strategy is as follows. First, we give a finite dimensional approximation of the $\Psi_2^4$ stochastic quantization equation (5.10). Thanks to this finite-dimensional approximation, see (5.13), we can apply the tools developed in Sections 4.2-4.3, namely the Itô and Girsanov's formulas, to provide a weak solution to (5.13), compare with Section 4.4. Finally, we prove that the process $Z^N$, i.e. the density of the signed expectation $\bar{\mathcal{E}}^{Z^N}$ that appears in Girsanov's formula, converges to a well-defined random variable $Z$ as the $N\to\infty$. As a result, we obtain a weak solution to the original equation (5.10) under the signed expectation $\bar{\mathcal{E}}^Z$.

Let us now delve into the details of the construction. We consider $\mathfrak{h} = L^2(\mathbb{T}^2; \mathbb{C}^4)$ and an analytic Brownian motion $X_t$ with covariance $U = (\mathbb{1} \oplus -\mathbb{1})$. Hereafter if $\chi: \mathbb{R}_+ \times \mathfrak{h} \to \mathcal{G}_X^p$ we write $\chi_t(x) = (\chi_t^1(x), \bar{\chi}_t^1(x), \chi_t^2(x), \bar{\chi}_t^2(x))$ to be the limit, if exists,

$$\chi_t^j(x) := \lim_{\tau \to +\infty} \chi_t(\delta_{x,\tau,j} \oplus 0), \quad \bar{\chi}_t^j(x) := \lim_{\tau \to +\infty} \chi_t(0 \oplus \delta_{x,\tau,j}),$$

which has to be understood in $\mathcal{S}'(\mathbb{T}^2, \mathcal{G}_X^p)$ (see Appendix A for the definitions of $\mathcal{S}(\mathbb{T}^2, \mathcal{G}_X^p)$ and $\mathcal{S}'(\mathbb{T}^2, \mathcal{G}_X^p)$). In this way, if $v(x) = (v^1(x), \bar{v}^1(x), v^2(x), \bar{v}^2(x)) \in C^\infty(\mathbb{T}^2, \mathbb{C}^4) \subset \mathfrak{h}$, we have

$$\chi_t(v) = \sum_{j=1,2} (\langle \chi_t^j(\cdot), v^j(\cdot) \rangle_{\mathcal{S}', \mathcal{S}} + \langle \bar{\chi}_t^j(\cdot), \bar{v}^j(\cdot) \rangle_{\mathcal{S}', \mathcal{S}}).$$

Using the described identification, when possible, of $\chi: \mathbb{R}_+ \times \mathfrak{h} \to \mathcal{G}_X^p$ with the function $\chi: \mathbb{R}_+ \to \mathcal{S}'(\mathbb{T}^2, \mathcal{G}_X^p)^4$ allows us to speak about the Besov regularity of a stochastic process. In particular, we say that a stochastic process $\chi: \mathbb{R}_+ \to \mathcal{S}'(\mathbb{T}^2, \mathcal{G}_X^\ell)^k$ has Besov regularity $B_{p,q}^s$ if $\chi_\cdot(\cdot) \in C^0(\mathbb{R}_+, B_{p,q}^s(\mathbb{T}^2, \mathcal{G}_X^\ell))$ - see Appendix A for the definitions of space of distributions and Besov spaces taking values in a Banach space.

**Remark 5.6.** In the particular case where $\chi_\cdot(\cdot) \in C^0(\mathbb{R}_+, B_{p,p}^s(\mathbb{T}^2, \mathcal{G}_X^p))$, we can consider the following sequence of stronger seminorms, for any $T \geqslant 0$

$$\|\chi\|_{C^0([0,T], B_{p,p}^s(\mathbb{T}^2, \mathcal{G}_X^p))}^p \lesssim \sup_{0 \leqslant t \leqslant T} \sum_{j \geqslant -1} 2^{psj} \int_{\mathbb{T}^2} \|K_j * \chi(x)\|_{\mathbb{L}^p}^p \mathrm{d}x,$$

where $K_j = \mathcal{F}_{\mathbb{T}^2}^{-1}(\varphi_j)$, $j \geqslant -1$, and $\{\varphi_j\}_{j \geqslant -1}$ is a dyadic partition of unity of $\mathbb{R}^2$, see Appendix A.

Using the notation introduced above, the $\Psi_2^4$ stochastic quantization SPDE reads, for any $0 \leqslant \theta < \frac{1}{2}$, $\lambda \in \mathbb{R}_+$,

$$\chi_t(x) = \chi_0(x) - \int_0^t (A^{1-2\theta} \chi_s(x) + \lambda A^{-2\theta} (\chi_s |\chi_s|^2)(x) - \lambda \infty A^{-2\theta} \chi_s) \mathrm{d}s + A^{-\theta} X_t(x), \quad x \in \mathbb{T}^2, \tag{5.10}$$

where $A = (-\Delta + m^2)$, where

$$|\chi|^2 = \sum_{j=1,2} \chi^1 \bar{\chi}^1, \tag{5.11}$$

$\chi = (\chi^1, \bar{\chi}^1, \chi^2, \bar{\chi}^2) \in C^0(\mathbb{R}_+, \mathcal{S}'(\mathbb{T}^2, \mathcal{G}_X^\ell)^4)$ (for some $\ell \geqslant 3$) and $X_t$ is some Brownian motion of covariance $U$, and the subtraction of $-\infty A^{-2\theta} \chi_s$ stands for a renormalisation procedure which will be explained below. As usual, the necessity of a renormalisation procedure is due to the low expected regularity of the solution to the SPDE (5.10).

**Lemma 5.7.** *For any $0 \leqslant \theta < \frac{1}{2}$, consider the process*

$$X_t^A = \mathrm{e}^{-A^{1-2\theta} t} \tilde{X}_0 + \int_0^t \mathrm{e}^{-A^{1-2\theta}(t-s)} A^{-\theta} \mathrm{d}X_s, \tag{5.12}$$

*where $X_0$ is an analytic (odd) Gaussian with covariance $A^{-1}U$ independent of the process $(X_t)_{t \geqslant 0}$, then for any, $\varepsilon > 0$, and $2 \leqslant p \leqslant \infty$ we have $X^A \in C^0(\mathbb{R}_+, B_{p,p}^{-\varepsilon}(\mathbb{T}^2, \mathcal{G}_X^p))$.*

**Proof.** The proof is given in the case $p = \infty$ in [ABDG22, Lemma 62]. The generic $p \geqslant 2$ can be proved using a similar method and hypercontractivity of Gaussian random variables of Theorem 3.20 (see also [DD03] for an analogous proof in the commutative case). □

**Notation 5.8.** *Let $\mathfrak{Q}(\chi,\psi)$ be a (local) antisymmetric polynomial of the (regular) random fields, we denote by $[\![\mathfrak{Q}(\chi,\psi)]\!]$ the same polynomial where every product between the components of $\chi$ and $\psi$ is replaced by the Wick product, where $\chi,\psi$ are Gaussian random field, more precisely we suppose $\chi = X_t^A$ and $\psi = A^{1-2\theta}X_t^A$.*

*For example if $\mathfrak{Q}(\chi,\psi) = P_N(\chi^1)(x)P_N(\bar{\chi}^1)(x)P_N(\chi^2)(x)P_N(\bar{\chi}^2)(x)$, where $P_N$ is the $(L^2(\mathbb{T}^2)$-)orthogonal projection on the Fourier modes less or equal then $N\geqslant 0$, and recalling that $\mathcal{E}(P_N(X_t^{A,1})(x)P_N(\bar{X}_t^{A,1})(x)) = \mathcal{E}(P_N(X_t^{A,1})(x)P_N(\bar{X}_t^{A,1})(x)) = C_N$ where $C_N = \sum_{|k|\leqslant N}\frac{1}{|k|^2+m^2}$ and $\mathcal{E}(X^{A,i}(x)\bar{X}^{A,j}(x)) = 0$, for any $i,j = 1,2$ and $i\neq j$, we have*

$$[\![\mathfrak{Q}(\chi,\psi)]\!] = \chi^1(x)\,\bar{\chi}^1(x)\,\chi^2(x)\,\bar{\chi}^2(x) - C_N\chi^1(x)\,\bar{\chi}^1(x) - C_N\chi^2(x)\,\bar{\chi}^2(x).$$

Using the symbols introduced in Notation 5.8, we consider an approximate stochastic quantization equation for $x\in\mathbb{T}^2$

$$\begin{cases} \chi_t^{(N)}(x) = \chi_0^{(N)}(x) - \int_0^t [A^{1-2\theta}\chi_s^{(N)}(x) + \lambda P_N(A^{-2\theta}[\![\chi_s^{(N)}|\chi_s^{(N)}|^2]\!])(x)]\mathrm{d}s + A^{-\theta}X_t(x), \\ \chi_0^{(N)}(x) = \tilde{X}_0(x) + h_0(x). \end{cases} \tag{5.13}$$

**Remark 5.9.** We take an initial condition of the form Gaussian free field plus some regular random variable, since we want to cover the case where the initial condition is distributed as the non-commutative measure $\omega^{(\lambda)}$ defined in Theorem 5.1. The fact that a non-commutative random variable distributed as $\omega^{(\lambda)}$ can be realized as a regular shift of a Gaussian free field is proved in [DFG22].

Let us now consider the process

$$\begin{aligned} Z_t^{N,h_0} = \exp\Bigg( &\lambda\int_0^t\int_{\mathbb{T}^2} P_N(A^{-\theta}[\![\mathfrak{P}_3(P_N(X_s^A + \mathrm{e}^{-A^{1-2\theta}s}h_0))]\!])(x)\mathrm{d}X_s(x) \\ &-\frac{\lambda^2}{2}\int_0^t\int_{\mathbb{T}^2}\langle P_N(A^{-\theta}[\![\mathfrak{P}_3(P_N(X_s^A + \mathrm{e}^{-A^{1-2\theta}s}h_0))]\!]), UA^{-\theta}[\![\mathfrak{P}_3(P_N(X_s^A + \mathrm{e}^{-A^{1-2\theta}s}h_0))]\!]\rangle_{\mathbb{R}^4}\mathrm{d}x\mathrm{d}s\Bigg), \end{aligned}$$

where

$$\mathfrak{P}_3(P_N(X_s^A + \mathrm{e}^{-A^{1-2\theta}s}h_0)) = P_N(X_s^A + \mathrm{e}^{-A^{1-2\theta}s}h_0)|P_N(X_s^A + \mathrm{e}^{-A^{1-2\theta}s}h_0)|^2,$$

and the Itô random field

$$B_t^{N,h_0}(x) = X_t(x) + \int_0^t P_N(A^{-2\theta}(P_N(X_s^A + \mathrm{e}^{-A^{1-2\theta}s}h_0)|P_N(X_s^A + \mathrm{e}^{-A^{1-2\theta}s}h_0)|^2))(x)\mathrm{d}s.$$

The reason of the introduction of the previous processes is the following weak representation of the solution to the approximating SPDE (5.13).

**Proposition 5.10.** *For any $h_0\in C^1(\mathbb{T}^2,\mathcal{G}_X^\infty)$ and for any $\varepsilon>0$, there is a unique (global in time) strong solution to $\chi_t^{(N)}\in C^0(\mathbb{R}_+, C^{-\varepsilon}(\mathbb{T}^2,\mathcal{G}_X^\infty))$ to equation (5.13) driven by the Brownian motion $X_t$. Furthermore the couple of processes $(X^A + \mathrm{e}^{A^{1-2\theta}\cdot}h_0, B^{N,h_0})$ is a weak solution to the SPDE (5.13) with respect to the expectation $\bar{\mathcal{E}}^{Z^{N,h_0}}(\cdot) = \mathcal{E}_0(\cdot Z_T^{N,h_0})$, namely we have that, for any polynomial $F\in\bigoplus_{n=0}^k\Lambda^n(\mathcal{S}(\mathbb{T}^2)^r)$ and for any $t_1<\ldots<t_r\in\mathbb{R}_+$,*

$$\begin{aligned} \omega_0(F(\chi_{t_1}^{(N)},\ldots,\chi_{t_r}^{(N)})) &= \bar{\mathcal{E}}_0^{Z^{N,h_0}}(F(X_{t_1}^A + \mathrm{e}^{-A^{1-2\theta}t_1}h_0,\ldots,X_{t_r}^A + \mathrm{e}^{-A^{1-2\theta}t_r}h_0)) \\ &= \mathcal{E}_0\big(F(X_{t_1}^A + \mathrm{e}^{-A^{1-2\theta}t_1}h_0,\ldots,X_{t_r}^A + \mathrm{e}^{-A^{1-2\theta}t_r}h_0)Z_{t_r}^{N,h_0}\big). \end{aligned} \tag{5.14}$$

**Proof.** The result follows from Theorems 4.51 and 4.52 by noting that it is possible to split equation (5.13) into two independent equations, by projecting the solution on the image of the projection $P_N$ and $I-P_N$. The first equation is a non-linear finite dimensional equation, and thus Theorems 4.51 and 4.52 apply directly. The second equation is the linear equation (5.12) projected on the image of $I-P_N$. A linear equation of the form (5.12) has always global in time solution and, by the independence from the processes $P_N(\chi_t^{(N)})$ and $Z^{N,h_0}$, equality (5.14) can be checked directly. □

In order to take the limit $N\to\infty$, in the weak solution $(X^A, B^{N,h_0})$ and obtaining a weak solution to the limit equation (5.10), we first need a result on the regularity of Wick polynomials of $X_t^A$.

**Lemma 5.11.** *Let $\mathfrak{Q}(\chi,\psi)$ be any antisymmetric local polynomial of the random fields $\chi,\psi$ of degree $n$, which is at most of first degree in $\psi$. Then, for any $\frac{2n(n-2)+3}{4n(n-2)+8}<\theta<\frac{1}{2}$, $p\geqslant 2$ and $0\leqslant s<\frac{2n(n-2)+4}{n}\theta-\frac{2n(n-2)+3}{2n}$, the sequence of random fields $[\![\mathfrak{Q}(P_N(X_t^A),P_N(A^{1-2\theta}X_t^A))]\!]$ is a Cauchy sequence in $C^0\left(\mathbb{R}_+,B_{p,p}^{s-\frac{1}{2}}(\mathbb{T}^2,\mathcal{G}_X^p)\right)$. We denote the limit by $[\![\mathfrak{Q}(X_t^A,A^{1-2\theta}X_t^A)]\!]$. If $\mathfrak{Q}$ does not depend on $\psi$, then $[\![\mathfrak{Q}(P_N(X_t^A))]\!]\to[\![\mathfrak{Q}(X_t^A)]\!]$ in $C^0(\mathbb{R}_+,B_{p,p}^{-\varepsilon}(\mathbb{T}^2,\mathcal{G}_X^p))$ for any $\varepsilon>0$. Finally, we have that, for any $T\geqslant 0$,*

$$\sup_{t\in[0,T]}\|\mathfrak{Q}(X_t^A,A^{1-2\theta}X_t^A)-\mathfrak{Q}(P_N(X_t^A),P_N(A^{1-2\theta}X_t^A))\|_{B_{p,p}^{s-\frac{1}{2}}(\mathbb{T}^2,\mathcal{G}_X^p)}\lesssim p^{\nu(\theta,s)} \tag{5.15}$$

$$\sup_{t\in[0,T]}\|\mathfrak{Q}(X_t^A)-\mathfrak{Q}(P_N(X_t^A))\|_{B_{p,p}^{-\varepsilon}(\mathbb{T}^2,\mathcal{G}_X^p)}\lesssim p^{\tilde{\nu}(\theta,\varepsilon)} \tag{5.16}$$

*for some $\nu(\theta,s),\tilde{\nu}(\theta,\varepsilon)<1$, not depending on $T$.*

**Proof.** We give the proof for the case where $\mathfrak{Q}(\chi,\psi)$ is a monomial of the form $\mathfrak{Q}(\chi,\psi)=\psi^r\prod_{k_1=1}^{n_2}\chi^{j_{k_1}}\prod_{k_2=1}^{n_2}\bar{\chi}^{j_{k_2}}$ or $\mathfrak{Q}(\chi,\psi)=\prod_{k_1=1}^{n_2}\chi^{j_{k_1}}\prod_{k_2=1}^{n_2}\bar{\chi}^{j_{k_2}}$, for some $n_1,n_2,j_{k_1},j_{k_2},r\in\mathbb{N}$. By linearity the general result follows.

Let us call

$$V_N^{j,r}=\int_{\mathbb{T}^2}\left[\!\!\left[P_N(A^{1-2\theta}X_t^{A,r})\prod_{k_1=1}^{n_2}P_N\big(\bar{X}_t^{A,j_{k_1}}\big)\prod_{k_2=1}^{n_2}P_N\big(\bar{X}_t^{A,j_{k_2}}\big)\right]\!\!\right](x)K_j(x)\mathrm{d}x,$$

$$V_N^{j}=\int_{\mathbb{T}^2}\left[\!\!\left[\prod_{k_1=1}^{n_2}P_N\big(\bar{X}_t^{A,j_{k_1}}\big)\prod_{k_2=1}^{n_2}P_N\big(\bar{X}_t^{A,j_{k_2}}\big)\right]\!\!\right](x)K_j(x)\mathrm{d}x,$$

where $K_j=\mathcal{F}^{-1}(\varphi_j)$ is the function corresponding to the $j$-th Littlewood-Paley block.

By the invariance of the law, and of the norm of $P_N(X_t^A)$ and $P_N(A^{1-\theta}\bar{X}_t^A)$ with respect to spatial and temporal translation, we have that, for any $N,N'\in\mathbb{N}$,

$$\begin{aligned}
&\sup_{t\in[0,T]}\|[\![\mathfrak{Q}(P_N(X_t^A),P_N(A^{1-2\theta}X_t^A))]\!]-[\![\mathfrak{Q}(P_{N'}(X_t^A),P_{N'}(A^{1-2\theta}X_t^A))]\!]\|^p_{B_{p,p}^{s-\frac{1}{2}}}\\
&\lesssim\sum_{r=1}^{4}\sum_{j\geqslant -1}2^{\left(s-\frac{1}{2}\right)jp}\|V_N^{j,r}-V_{N'}^{j,r}\|^p_{\mathbb{L}^p}\sup_{t\in[0,T]}\|[\![\mathfrak{Q}(P_N(X_t^A),P_N(A^{1-2\theta}X_t^A))]\!]\|^p_{B_{p,p}^{s-\frac{1}{2}}} &(5.17)\\
&\lesssim\sum_{r=1}^{4}\sum_{j\geqslant -1}2^{-\left(s-\frac{1}{2}\right)jp}\|V_N^{j,r}\|^p_{\mathbb{L}^p}\mathrm{d}x, &(5.18)
\end{aligned}$$

and similar expression with $V_N^j$. Since $[\![\mathfrak{Q}(P_N(X_t^A),P_N(A^{1-2\theta}X_t^A))]\!]$, is a polynomial of $P_N(X_t^A)$ and $P_N(A^{1-\theta}\bar{X}_t^A)$ with coefficient bounded by $\log(N)$ and $N^{2-4\theta}$, by following the proof of Lemma 5.4 we get

$$\sup_{t\in[0,T],x\in\mathbb{T}^2}\|[\![\mathfrak{Q}(P_N(X_t^A),P_N(A^{1-2\theta}X_t^A))]\!](x)\|^2_{\mathbb{L}^\infty}\lesssim(\log N)^{n_1+n_2}N^{4-8\theta},$$

and thus, since $\|K_j\|_{L^1(\mathbb{T}^2)}=1$, $\|V_N^{j,r}\|_{\mathbb{L}^\infty}\lesssim \|K_j\|_{L^1(\mathbb{T}^2)}(\log N)^{n_1+n_2}N^{4-8\theta}\lesssim(\log N)^{n_1+n_2}N^{4-8\theta}$ (when $V_N^{j,r}$ is replaced by $V_N^j$ we obtain a bound proportional to $(\log(N))^{n_1+n_2}$).

Following again the proof of Lemma 5.4, we obtain that

$$\begin{aligned}
&\|T_\tau^{(2)}(V_{N'}^{j,r})-T_\tau^{(2)}(V_N^{j,r})\|_{L^2}^2\\
&=C_\tau\sum_{h,\ell_1,\ldots,\ell_{n_1+n_2}\in\mathbb{Z}^2}\varphi_j^2\left(h+\sum_{i=1}^{n_1+n_2}\ell_i\right)\left(\mathbf{1}_{B_1}\left(\frac{|h|}{N'}\right)\prod_{i=1}^{n_1+n_2}\mathbf{1}_{B_1}\left(\frac{|\ell_i|}{N'}\right)-\mathbf{1}_{B_1}\left(\frac{|g|}{N}\right)\prod_{i=1}^{n_1+n_2}\mathbb{I}_{B_1}\left(\frac{|\ell_i|}{N}\right)\right)\cdot\\
&\qquad\cdot\frac{1}{(|h|^2+m^2)^{4\theta-1}}\prod_{i=1}^{n_1+n_2}\frac{1}{(|\ell_i|^2+m^2)}\\
&=C_\tau\sum_{h,\ell_1,\ldots,\ell_{n_1+n_2}\in\mathbb{Z}^2}\varphi_j^2\left(h+\sum_{i=1}^{n_1+n_2}\ell_i\right)\left(\mathbf{1}_{B_{N'}\setminus B_N}(|h|)\prod_{i=1}^{n_1+n_2}\mathbf{1}_{B_1}\left(\frac{|\ell_i|}{N'}\right)\right)\frac{1}{(|h|^2+m^2)^{4\theta-1}}\prod_{i=1}^{n_1+n_2}\frac{1}{(|\ell_i|^2+m^2)}\\
&\quad+\sum_{h,\ell_1,\ldots,\ell_{n_1+n_2}\in\mathbb{Z}^2}\sum_{k=1}^{n_1+n_2}\varphi_j^2\left(h+\sum_{i=1}^{n_1+n_2}\ell_i\right)\left(\mathbf{1}_{B_{N'}\setminus B_N}(|\ell_k|)\mathbf{1}_{B_1}\left(\frac{|h|}{N}\right)\prod_{i=k+1}^{n_1+n_2}\mathbf{1}_{B_1}\left(\frac{|\ell_i|}{N}\right)\prod_{i=1}^{k-1}\mathbf{1}_{B_1}\left(\frac{|\ell_i|}{N'}\right)\right)\cdot\\
&\qquad\cdot\frac{1}{(|h|^2+m^2)^{4\theta-1}}\prod_{i=1}^{n_1+n_2}\frac{1}{(|\ell_i|^2+m^2)}. \qquad (5.19)
\end{aligned}$$

Following the proof of [DD03, Lemma 3.2], we get

$$(5.19)\lesssim 2^{-j\tilde{s}}\left(\|\tilde{\gamma}_{N,N'}\gamma_{N'}^{n_1+n_2}\|_{H^{1+\tilde{s}}}+\sum_{k=1}^{n_1+n_2}\|\tilde{\gamma}_N\gamma_{N'}^{k-1}\gamma_N^{n_1+n_2-k}\gamma_{N,N'}\|_{H^{1+\tilde{s}}}\right)$$

where

$$\begin{aligned}
\gamma_L(x)&=\sum_{|k|\leqslant L,k\in\mathbb{Z}^2}\frac{\mathrm{e}^{ik\cdot x}}{(m^2+|k|^2)}, &\gamma_{L,'L}(x)&=\sum_{L<|k|\leqslant L',k\in\mathbb{Z}^2}\frac{\mathrm{e}^{ik\cdot x}}{(m^2+|k|^2)},\\
\tilde{\gamma}_L(x)&=\sum_{|k|\leqslant L,k\in\mathbb{Z}^2}\frac{\mathrm{e}^{ik\cdot x}}{(m^2+|k|^2)^{2\theta-1}}, &\tilde{\gamma}_{L,'L}(x)&=\sum_{L<|k|\leqslant L',k\in\mathbb{Z}^2}\frac{\mathrm{e}^{ik\cdot x}}{(m^2+|k|^2)^{2\theta-1}}.
\end{aligned}$$

If we choose, $\tilde{s}<0$, $\delta,\delta'>0$ and $p',q'\geqslant 1$ such that

$$\tilde{s}>2s-1,\quad \delta,\delta'>(n-2)(2-4\theta),\quad \delta>(\tilde{s}+1)-\frac{1}{q'(n-1)},\quad ,8\theta-4-\delta'>(\tilde{s}+1)-\frac{1}{p'}\quad \frac{1}{p'}+\frac{1}{q'}=1$$

(the existence of the previous constant is ensured by the conditions $\frac{2n(n-2)+3}{4n(n-2)+8}<\theta<\frac{1}{2}$, and $0\leqslant s<\frac{2n(n-2)+4}{n}\theta-\frac{2n(n-2)+3}{2n}$), using the Besov embedding Theorem A.3 for the computation of $\|\tilde{\gamma}_{N,N'}\|_{B^{1+\tilde{s}}_{p,p}}$, $\|\gamma_{N,N'}\|_{B^{1+\tilde{s}}_{p,p}}$, and so on we get that

$$\|V_{N'}^{j,r}-V_N^{j,r}\|_{\mathbb{L}^2}^2=\sup_{|\tau|\leqslant\frac{3}{4}}\|T_\tau^{(2)}(V_{N'}^{j,r})-T_\tau^{(2)}(V_N^{j,r})\|_{L^2}^2\lesssim 2^{-j\tilde{s}}(N^{-\delta'}+N^{-\delta}). \qquad (5.20)$$

Since $[\![\Omega(P_N(X_t^A),P_N(A^{1-2\theta}X_t^A))]\!]$ is a polynomial of degree $n=n_1+n_2+1$, by the hypercontractivity inequalities (5.17) and (5.20) we have that $[\![\Omega(P_N(X_t^A),P_N(A^{1-2\theta}X_t^A))]\!]$ is a Cauchy sequence in $C^0\left(\mathbb{R}_+,B_{p,p}^{s-\frac{1}{2}}(\mathbb{T}^2,\mathcal{G}_X^p)\right)$. In order to get inequality (5.15), we proceed as in the proof of Theorem 5.3. Indeed, for any $p\geqslant 2$, we have, for any $N'>0$,

$$\begin{aligned}
\|V_{N'}^{j,r}\|_{\mathbb{L}^p}&\leqslant\|V_N^{j,r}\|_{\mathbb{L}^\infty}+\|V_{N'}^{j,r}-V_N^{j,r}\|_{\mathbb{L}^p}\lesssim(\log N)^{n_1+n_2}N^{2-4\theta}+p^{\frac{n}{2}}2^{-j\tilde{s}}(N^{-\delta}+N^{-\delta'}).\\
&\lesssim 2^{j\tilde{s}}\left((\log N)^{n_1+n_2}N^{2-4\theta}+p^{\frac{n}{2}}(N^{-\delta}+N^{-\delta'})\right) \qquad (5.21)
\end{aligned}$$

By our assumptions on $\delta, \delta'$ there is $\alpha < \frac{1}{2-4\theta}$ such that, $(2-4\theta)\alpha < \nu(\theta,s)$, $\frac{n}{2} - \alpha\delta' \leqslant \nu(\theta,s)$ and $\frac{n}{2} - \alpha\delta \leqslant \nu(\theta,s)$, for a suitable $0 < \nu(\theta,s) < 1$. This mean that by choosing $N = p^{\alpha}$ we get

$$\|V_{N'}^{j,r}\|_{\mathbb{L}^p} \lesssim 2^{j\tilde{s}}\Big((\log p)^{n_1+n_2} p^{\alpha(2-4\theta)} + p^{\frac{n}{2}-\alpha\delta} + p^{\frac{n}{2}-\alpha\delta'}\Big) \lesssim 2^{j\tilde{s}} p^{\nu(\theta,s)}.$$

By using inequality (5.18) the claim (5.15) follows.

In order to prove the convergence $[\![\mathfrak{Q}(P_N(X_t^A))]\!] \to [\![\mathfrak{Q}(X_t^A)]\!]$ and inequality (5.16), we can repeat the previous reasoning by replacing $V_N^{j,r}$ by $V^j$. The reason for the arbitrary (negative) Besov regularity $-\varepsilon$ in this second case is due to the fact that $\|V_N^j\|_{\mathbb{L}^\infty} \lesssim (\log(N))^{n_1+n_2}$. □

We observe now that, by the Itô formula,

$$\begin{aligned}
&Z_t^{N,h_0}\\
&= \exp\Big(\lambda \int_{\mathbb{T}^2} [\![P_N(|P_N(X_t^A + \mathrm{e}^{-A^{1-2\theta}s}h_0)|^4)]\!]\mathrm{d}x - \lambda \int_{\mathbb{T}^2} [\![P_N(|P_N(\tilde{X}_0 + h_0)|^4)]\!]\mathrm{d}x\\
&+ \lambda \int_0^t \int_{\mathbb{T}^2} [\![[P_N(X_s^A + \mathrm{e}^{-A^{1-2\theta}s}h_0) \cdot A^{1-2\theta}(P_N(X_s^A + \mathrm{e}^{-A^{1-2\theta}s}h_0))]|P_N(X_t^A + \mathrm{e}^{-A^{1-2\theta}s}h_0)|^2]\!]\mathrm{d}x\mathrm{d}s\\
&- \frac{\lambda^2}{2}\int_0^t \int_{\mathbb{T}^2} \langle P_N(A^{-\theta}[\![\mathfrak{P}_3(P_N(X_s^A + \mathrm{e}^{-A^{1-2\theta}s}h_0))]\!]),\\
&\qquad UA^{-\theta}[\![\mathfrak{P}_3(P_N(X_s^A + \mathrm{e}^{-A^{1-2\theta}s}h_0))]\!]\rangle_{\mathbb{R}^4}\mathrm{d}x\mathrm{d}s\Big),
\end{aligned} \tag{5.22}$$

where we use the notation

$$\chi \cdot \psi = \sum_{j=1,2} (\chi^i \bar{\psi}^i + \bar{\chi}^i \psi^i).$$

**Lemma 5.12.** *Consider $\frac{19}{40} < \theta < \frac{1}{2}$, under the hypotheses of Proposition 5.10, then the addends in the sum defining the exponential (5.22) (namely $[\![|(X_t^A + \mathrm{e}^{-A^{1-2\theta}t}h_0)|^4]\!]$, $[\![P_N(|P_N(\tilde{X}_0 + h_0)|^4)]\!]$ etc.) converges, as $N \to +\infty$, to some well defined random processes which we denote by*

$$[\![|(X_t^A + \mathrm{e}^{-A^{1-2\theta}t}h_0)|^4]\!], \quad [\![[(X_s^A + \mathrm{e}^{-A^{1-2\theta}s}h_0) \cdot A^{1-2\theta}(X_s^A + \mathrm{e}^{-A^{1-2\theta}s}h_0)]|X_t^A + \mathrm{e}^{-A^{1-2\theta}s}h_0|^2]\!], \tag{5.23}$$

$$\langle A^{-\theta}[\![(X_s^A + \mathrm{e}^{-A^{1-2\theta}s}h_0)|(X_s^A + \mathrm{e}^{-A^{1-2\theta}s}h_0)|^2]\!], UA^{-\theta}[\![(X_s^A + \mathrm{e}^{-A^{1-2\theta}s}h_0)|(X_s^A + \mathrm{e}^{-A^{1-2\theta}s}h_0)|^2]\!]\rangle_{\mathbb{R}^4}. \tag{5.24}$$

**Proof.** Thanks to Lemma 5.11 and Theorem A.4 (on the multiplication of distributions in Besov spaces), the expressions in equation (5.23) and equation (5.24) can be defined in the following way: taking for example $[\![|(X_t^A + \mathrm{e}^{-A^{1-2\theta}t}h_0)|^4]\!]$, if we (formally) expand the fourth power and using the properties of Wick products with respect to addition we can define $[\![|(X_t^A + e^{-A^{1-2\theta}t}h_0)|^4]\!]$

$$\begin{aligned}
[\![|(X_t^A + \mathrm{e}^{-A^{1-2\theta}t}h_0)|^4]\!] := {} & [\![|X_t^A|^4]\!] + 2[\![|X_t^A|^2 X_t^A]\!] \cdot \mathrm{e}^{-A^{1-2\theta}t}h_0 + 2[\![|X_t^A|^2]\!]|\mathrm{e}^{-A^{1-2\theta}t}h_0|^2\\
&+ 2[\![X_t^{A,1} X_t^{A,2}]\!](\mathrm{e}^{-A^{1-2\theta}t}\bar{h}_0^1 \mathrm{e}^{-A^{1-2\theta}t}\bar{h}_0^2) + [\![\bar{X}_t^{A,1}\bar{X}_t^{A,2}]\!](\mathrm{e}^{-A^{1-2\theta}t}h_0^1 \mathrm{e}^{-A^{1-2\theta}t}h_0^2)\\
&+ 2(X_t^A \cdot \mathrm{e}^{-A^{1-2\theta}t}h_0)|\mathrm{e}^{-A^{1-2\theta}t}h_0|^2 + |\mathrm{e}^{-A^{1-2\theta}t}h_0|^4.
\end{aligned} \tag{5.25}$$

The right hand side of the expression (5.25) is well-defined since $[\![|X_t^A|^4]\!]$, $[\![|X_t^A|^2 X_t^A]\!]$, etc. are, by Lemma 5.11, random distributions in $B_{p,p}^{-\varepsilon}(\mathbb{T}^2, \mathcal{G}_X^p)$, meanwhile $\mathrm{e}^{-A^{1-2\theta}t}h_0$, $|\mathrm{e}^{-A^{1-2\theta}t}h_0|^2$, etc. are random fields in $C^1(\mathbb{T}^2, \mathcal{G}_X^\infty)$. Since $-\varepsilon + 1 > 0$ (for $\varepsilon$ small enough) and $\frac{1}{p} + \frac{1}{\infty} \leqslant 1$ (for any $p \geqslant 1$), by Theorem A.4, the products $[\![|X_t^A|^2 X_t^A]\!] \cdot \mathrm{e}^{-A^{1-2\theta}t}h_0$, $[\![|X_t^A|^2]\!]|\mathrm{e}^{-A^{1-2\theta}t}h_0|^2$, etc. are well-defined as random distribution in $B_{p,p}^{-\varepsilon}(\mathbb{T}^2, \mathcal{G}_X^\infty)$. In a similar way we can define all the other random distribution in expression (5.23) and (5.24).

By observing that

$$
\begin{aligned}
\llbracket |P_N(X_t^A+\mathrm{e}^{-A^{1-2\theta}t}h_0)|^4\rrbracket &= \llbracket |P_NX_t^A|^4\rrbracket + 2\llbracket |P_NX_t^A|^2P_NX_t^A\rrbracket\cdot \mathrm{e}^{-A^{1-2\theta}t}P_Nh_0 + 2\llbracket |P_NX_t^A|^2\rrbracket |\mathrm{e}^{-A^{1-2\theta}t}P_Nh_0|^2 + \\
&\quad + \llbracket P_NX_t^A\otimes P_NX_t^A\rrbracket\cdot(\mathrm{e}^{-A^{1-2\theta}t}P_Nh_0\otimes \mathrm{e}^{-A^{1-2\theta}t}P_Nh_0) + |\mathrm{e}^{-A^{1-2\theta}t}P_N(h_0)|^4,
\end{aligned} \tag{5.26}
$$

it is easy to see that, by Lemma 5.11, Theorem A.4 and the convergence of $\mathrm{e}^{-A^{1-2\theta}t}P_N(h_0)$ to $\mathrm{e}^{-A^{1-2\theta}}h_0$ in $C^{1-\varepsilon}(\mathbb{T}^2,\mathcal{G}_X^{\infty})$, as $N\to+\infty$, the expression at the right hand side of equation (5.26) converges to the sum (5.12), as $N\to+\infty$. The convergence of the other terms to the random processes in equations (5.23) and (5.24) can be done in a similar way. □

We can now, prove that $Z_t^{N,h_0}$ has a limit as $N\to\infty$.

**Lemma 5.13.** *Suppose that $\frac{19}{40}<\theta<\frac{1}{2}$. Then, for any $p\geqslant 2$ we have that $Z_t^{N,h_0}\to Z_t^{h_0}$ in $C^0(\mathbb{R}_+,\mathcal{G}_X^p)$ where $Z_t^{h_0}$ is as follows*

$$
\begin{aligned}
Z_t^{h_0} = \exp\Big(&\lambda\int_{\mathbb{T}^2}\llbracket |(X_t^A+\mathrm{e}^{-A^{1-2\theta}t}h_0)|^4\rrbracket \mathrm{d}x - \lambda\int_{\mathbb{T}^2}\llbracket |(X_0^A+h_0)|^4\rrbracket \mathrm{d}x \\
&+\lambda\int_0^t\int_{\mathbb{T}^2}\llbracket [(X_s^A+\mathrm{e}^{-A^{1-2\theta}s}h_0)\cdot A^{1-2\theta}(X_s^A+\mathrm{e}^{-A^{1-2\theta}s}h_0)]|X_t^A+\mathrm{e}^{-A^{1-2\theta}s}h_0|^2\rrbracket \mathrm{d}x\mathrm{d}s \\
&-\frac{\lambda^2}{2}\int_0^t\int_{\mathbb{T}^2}\langle A^{-\theta}\llbracket (X_s^A+\mathrm{e}^{-A^{1-2\theta}s}h_0)|(X_s^A+\mathrm{e}^{-A^{1-2\theta}s}h_0)|^2\rrbracket, \\
&\qquad UA^{-\theta}\llbracket (X_s^A+\mathrm{e}^{-A^{1-2\theta}s}h_0)|(X_s^A+\mathrm{e}^{-A^{1-2\theta}s}h_0)|^2\rrbracket\rangle_{\mathbb{R}^4}\mathrm{d}x\mathrm{d}t\Big),
\end{aligned} \tag{5.27}
$$

*where the expressions $\llbracket |(X_t^A+\mathrm{e}^{-A^{1-2\theta}t}h_0)|^4\rrbracket$, $\llbracket |(X_0^A+h_0)|^4\rrbracket$, etc. in formula (5.27) are defined as explained in Lemma 5.12.*

**Proof.** We prove only that $\exp\big(\lambda\int_{\mathbb{T}^2}\llbracket P_N(|P_N(X_t^A+\mathrm{e}^{-A^{1-2\theta}t}h_0)|^4)\rrbracket \mathrm{d}x\big)$ converges to

$$
\begin{aligned}
\exp\Big(\lambda\int_{\mathbb{T}^2}&(\llbracket |X_t^A|^4\rrbracket + 2\llbracket |X_t^A|^2X_t^A\rrbracket\cdot \mathrm{e}^{-A^{1-2\theta}t}h_0 + 2\llbracket |X_t^A|^2\rrbracket |\mathrm{e}^{-A^{1-2\theta}t}h_0|^2 + \\
&+\llbracket X_t^A\otimes X_t^A\rrbracket\cdot(\mathrm{e}^{-At}h_0\otimes e^{-At}h_0) + |\mathrm{e}^{-At}h_0|^4)\mathrm{d}x\Big)
\end{aligned} \tag{5.28}
$$

where

$$
\begin{aligned}
&\llbracket X_t^A\otimes X_t^A\rrbracket\cdot(\mathrm{e}^{-At}h_0\otimes \mathrm{e}^{-At}h_0) \\
&\quad = 2\llbracket X_t^{A,1}X_t^{A,2}\rrbracket(\mathrm{e}^{-A^{1-2\theta}t}\bar{h}_0^1\mathrm{e}^{-A^{1-2\theta}t}\bar{h}_0^2) + \llbracket \bar{X}_t^{A,1}\bar{X}_t^{A,2}\rrbracket(\mathrm{e}^{-A^{1-2\theta}t}h_0^1\mathrm{e}^{-A^{1-2\theta}t}h_0^2),
\end{aligned}
$$

in $C^0(\mathbb{R}_+,\mathcal{G}_X^p)$, since the convergence of all other terms can be proved in a similar way. By equation (5.26), in order to prove that $\exp\big(\lambda\int_{\mathbb{T}^2}\llbracket P_N(|P_N(X_t^A+\mathrm{e}^{-A^{1-2\theta}t}h_0)|^4)\rrbracket \mathrm{d}x\big)$ converges to the exponential in equation (5.28), it is enough to prove that $\exp(\lambda\int_{\mathbb{T}^2}\llbracket |P_NX_t^A|^4\rrbracket \mathrm{d}x)$ converges to $\exp(\lambda\int_{\mathbb{T}^2}\llbracket |X_t^A|^4\rrbracket \mathrm{d}x)$, that $\exp\big(2\lambda\int_{\mathbb{T}^2}\llbracket |P_NX_t^A|^2P_NX_t^A\rrbracket\cdot \mathrm{e}^{-A^{1-2\theta}t}P_Nh_0\mathrm{d}x\big)$ converges to $\exp\big(\lambda\int_{\mathbb{T}^2}2\llbracket |X_t^A|^2X_t^A\rrbracket\cdot \mathrm{e}^{-A^{1-2\theta}t}h_0\mathrm{d}x\big)$ etc.

The convergence of $\exp(\lambda\int\llbracket |P_NX_t^A|^4\rrbracket \mathrm{d}x)$ to $\exp(\lambda\int\llbracket |X_t^A|^4\rrbracket \mathrm{d}x)$, has been proved in Theorem 5.3.

If we consider the term $\llbracket |P_NX_t^A|^2P_NX_t^A\rrbracket\cdot \mathrm{e}^{-At}P_Nh_0$, by Lemma 5.11 and Theorem A.4, for any $k\in\mathbb{N}$ and $p\geqslant 2$, we have that

$$
\begin{aligned}
&\Big\|\int\llbracket |P_NX_t^A|^2P_NX_t^A\rrbracket\cdot \mathrm{e}^{-A^{1-2\theta}t}P_Nh_0\mathrm{d}x - \int\llbracket |X_t^A|^2X_t^{A^{1-2\theta}}\rrbracket\cdot \mathrm{e}^{-A^{1-2\theta}t}h_0\mathrm{d}x\Big\|_{\mathbb{L}^{kp}} \\
&\leqslant \|\llbracket |P_NX_t^A|^2P_NX_t^A\rrbracket - \llbracket |X_t^A|^2X_t^A\rrbracket\|_{B_{kp,kp}^{-\varepsilon}}\|h_0\|_{C^1} + \|\llbracket |X_t^A|^2X_t^A\rrbracket\|_{B_{kp,kp}^{-\varepsilon}}N^{(2\varepsilon-1)}\|h_0\|_{C^1}\to 0
\end{aligned} \tag{5.29}
$$

as $N \to \infty$. Furthermore, by inequality (5.15) and Theorem A.4, we get

$$\begin{aligned}
&\left\| \sum_{k=0}^{n} \frac{\lambda^k}{k!} \Big( \int [\![|P_N X_t^A|^2 P_N X_t^A]\!] \cdot \mathrm{e}^{-A^{1-2\theta}t} P_N h_0 \mathrm{d}x - \int [\![|X_t^A|^2 X_t^A]\!] \cdot \mathrm{e}^{-A^{1-2\theta}t} h_0 \mathrm{d}x \Big)^k \right\|_{\mathbb{L}^p} \\
&\lesssim \sum_{k=0}^{n} \frac{\lambda^\lambda}{k!} \left\| \int [\![|P_N X_t^A|^2 P_N X_t^A]\!] \cdot \mathrm{e}^{-A^{1-2\theta}t} P_N h_0 \mathrm{d}x - \int [\![|X_t^A|^2 X_t^A]\!] \cdot \mathrm{e}^{-A^{1-2\theta}t} h_0 \mathrm{d}x \right\|_{\mathbb{L}^{kp}} \\
&\lesssim \sum_{k=0}^{n} \frac{\lambda^k \|h_0\|_{C^1}^k}{k!} \big( \|[\![|P_N X_t^A|^2 P_N X_t^A]\!] - [\![|X_t^A|^2 X_t^A]\!]\|_{B_{kp,kp}^{-\varepsilon}}^k + N^{k(2\varepsilon-1)} \|[\![|X_t^A|^2 X_t^A]\!]\|_{B_{kp,kp}^{-\varepsilon}}^k \big) \\
&\lesssim \sum_{k=0}^{n} C_p^k k^{(\nu(\theta,\varepsilon)-1)},
\end{aligned} \tag{5.30}$$

where $C_p > 1$ is a suitable constant depending on $p$. Since $\sum_{k=0}^{n} C_p^k k^{(\nu(\theta,\varepsilon)-1)}$ is a convergent sequence, inequality (5.30) proves that $\exp\big(\lambda \int [\![|X_t^A|^2 X_t^A]\!] \cdot \mathrm{e}^{-A^{1-2\theta}t} h_0 \, \mathrm{d}x\big)$ is well-defined in $\mathcal{G}_X^p$. Furthermore, by the convergence of (5.29) and Lebesgue dominated convergence theorem, inequality (5.30) implies that $\exp\big(\lambda \int [\![|P_N X_t^A|^2 P_N X_t^A]\!] \cdot \mathrm{e}^{-A^{1-2\theta}t} P_N h_0 \mathrm{d}x\big)$ converges to $\exp\big(\lambda \int [\![|X_t^A|^2 X_t^A]\!] \cdot \mathrm{e}^{-A^{1-2\theta}t} h_0 \mathrm{d}x\big)$ in $C^0(\mathbb{R}_+, \mathcal{G}_X^p)$. In the same way, it is possible to prove that the exponential of every term in the sum (5.26) converges to the exponential of the corresponding term of equation (5.28). Since every term in the exponential is even, and thus the standard properties of the products of exponentials hold, and since we proved the convergence of each single exponential in $C^0(\mathbb{R}_+, \mathcal{G}_X^p)$, for any arbitrary $p \geqslant 2$, the statement follows from Hölder's inequality for the twisted spaces. □

We have now all the tools for proving the convergence of weak solutions to equation (5.13) to the weak solution to equation (5.10).

**Theorem 5.14.** *Let $\frac{19}{40} < \theta < \frac{1}{2}$ and $h_0 \in C^1(\mathbb{T}^2, \mathcal{G}_X^\infty)$. Then, for any $F \in \bigoplus_{n=0}^{k} \Lambda^n(\mathcal{S}(\mathbb{T}^2)^r)$ and any $t_1 < \cdots < t_r \in \mathbb{R}_+$ we have*

$$\begin{aligned}
\lim_{N \to \infty} \omega_0(F(\chi_{t_1}^{(N)}, \ldots, \chi_{t_r}^{(N)})) &= \bar{\mathcal{E}}_0^Z(F(X_{t_1}^A + \mathrm{e}^{-A^{1-2\theta}t_1} h_0, \ldots, X_{t_r}^A + \mathrm{e}^{-A^{1-2\theta}t_1} h_0)) \\
&= \mathcal{E}_0(F(X_{t_1}^A + \mathrm{e}^{-A^{1-2\theta}t_1} h_0, \ldots, X_{t_r}^A + \mathrm{e}^{-A^{1-2\theta}t_1} h_0) Z_{t_r}),
\end{aligned}$$

*where $\chi_t^{(N)}$ is the solution to equation (5.13).*

**Proof.** The proof is a consequence of the representation of solutions to equation (5.13) given in Proposition 5.10 and of the convergence of the process $Z^{N,h_0}$ to $Z^{h_0}$ provided by Lemma 5.13. □

# Appendix A Besov spaces of distributions taking values in Banach spaces

In this appendix, we want to recall some results about Besov spaces of functions (or distributions) on $\mathbb{T}^d$ taking values in a Banach space $A$. The results of this section can be found in [Ama97, Ama19], where the theory of Besov spaces taking values in a Banach space has been developed.

We denote by $\mathcal{S}(\mathbb{T}^d)$ the space of smooth functions defined on $\mathbb{T}^d$ and equipped with the set of seminorms

$$\|f\|_\alpha := \| |D^\alpha f| \|_{L^\infty(\mathbb{T}^d)} < \infty,$$

where $\alpha\in\mathbb{N}^d$. We denote by $\mathscr{S}'(\mathbb{T}^d)$ the strong dual of $\mathscr{S}(\mathbb{T}^d)$ with respect to the topology induced by the seminorms $\|\cdot\|_\alpha$. If $A$ is a Banach space and $B$ a nuclear space, we denote by $B\hat{\otimes}A$ the completion of $B\otimes A$ with respect to the natural metric of the algebraic tensor product on $B\otimes A$. Such a completion is unique up to an isomorphism. Using this notation we define

$$\mathscr{S}(\mathbb{T}^d,A):=\mathscr{S}(\mathbb{T}^d)\hat{\otimes}A,\quad \mathscr{S}'(\mathbb{T}^d,A):=\mathscr{S}'(\mathbb{T}^d)\hat{\otimes}A,$$

where $A$ is a Banach space. It is important to note that

$$\mathscr{S}(\mathbb{T}^d,A)\overset{d}{\hookrightarrow}\mathscr{S}'(\mathbb{T}^d,A)$$

where the arrow means that a space is continuously embedded and dense in the following one.

**Remark A.1.** Let $A_1$, $A_2$ and $A_3$ be three Banach spaces, for which a product operation $\cdot: A_1\times A_2\to A_3$ (which is a continuous bilinear function) is defined. Then, for any $i=1,2,3$, it is possible to define uniquely $\langle\cdot,\cdot\rangle:\mathscr{S}(\mathbb{T}^d,A_1)\times\mathscr{S}'(\mathbb{T}^d,A_2)\to A_3$, $\cdot:\mathscr{S}(\mathbb{T}^d,A_2)\times\mathscr{S}'(\mathbb{T}^d,A_2)\to\mathscr{S}'(\mathbb{T}^d,A_3)$, $\cdot:\mathscr{S}(\mathbb{T}^d)\times\mathscr{S}'(\mathbb{T}^d,A_i)\to\mathscr{S}'(\mathbb{T}^d,A_i)$, $*:\mathscr{S}(\mathbb{T}^d)\times\mathscr{S}'(\mathbb{T}^d,A_i)\to\mathscr{S}(\mathbb{T}^d,A_i)$ and $D^\alpha:\mathscr{S}'(\mathbb{T}^d,A)\to\mathscr{S}'(\mathbb{T}^d,A)$ (where $\alpha\in\mathbb{N}^d$) which extend in an continuous way the following operations: any $f\in\mathscr{S}(\mathbb{T}^d)$, $u\in\mathscr{S}'(\mathbb{T}^d)$ and $a\in A_i$, $a_1\in A_1$, $a_2\in A_2$ we have

$$\begin{aligned}\langle f\otimes a_1,u\otimes a_2\rangle &= \langle f,u\rangle a_1a_2\\ (f\otimes a_1)\cdot(u\otimes a_2) &= (fu)\otimes(a_1a_2)\\ f\cdot(u\otimes a) &= (fu)\otimes a\\ f*(u\otimes a) &= (f*u)\otimes a\\ D^\alpha(u\otimes a) &= (D^\alpha u)\otimes a\end{aligned}$$

where $\langle f,u\rangle$ is the normal pairing in $\mathscr{S}(\mathbb{T}^d)\times\mathscr{S}'(\mathbb{T}^d)$, $(fu)$ is the product in $\mathscr{S}(\mathbb{T}^d)\times\mathscr{S}'(\mathbb{T}^d)$, $(f*u)$ is the convolution in $\mathscr{S}(\mathbb{T}^d)\times\mathscr{S}'(\mathbb{T}^d)$ and $D^\alpha$ is the $\alpha$ derivative in $\mathscr{S}'(\mathbb{T}^d)$ (see [Ama19] Appendix 1).

We recall the definition of Littlewood–Paley decomposition on the torus $\mathbb{T}^d$. Let $\chi,\varphi$ be smooth non-negative functions from $\mathbb{R}^d$ to $\mathbb{R}$ such that

- $\mathrm{supp}(\chi)\subset B_{\frac{4}{3}}(0)$ and $\mathrm{supp}(\varphi)\subset B_{\frac{8}{3}}(0)\setminus B_{\frac{3}{4}}(0)$,
- $\chi,\varphi\leqslant 1$ and $\chi(y)+\sum_{j\geq 0}\varphi(2^{-j}y)=1$ for any $y\in\mathbb{R}^n$,
- $\mathrm{supp}(\chi)\cap\mathrm{supp}(\varphi(2^{-i}\cdot))=\emptyset$ for $i\geqslant 1$,
- $\mathrm{supp}(\varphi(2^{-j}\cdot))\cap\mathrm{supp}(\varphi(2^{-i}\cdot))=\emptyset$ if $|i-j|>1$,

where by $B_r(x)$ we denote the ball centered at $x\in\mathbb{R}^d$ and of radius $r>0$. We introduce the following notation: $\varphi_{-1}=\chi$, $\varphi_j(\cdot)=\varphi(2^{-j}\cdot)$, $K_j=\mathcal{F}^{-1}(\varphi_j|_{\mathbb{Z}^d})\in\mathscr{S}(\mathbb{T}^d)$.

If $v\in\mathscr{S}'(\mathbb{T}^d,A)$ and if $i\in\mathbb{Z}$, $i\geqslant -1$ we define the $i$th Littlewood–Paley block as follows

$$\Delta_i v=K_i*v\in\mathscr{S}(\mathbb{T}^d,A).$$

Then, if $s\in\mathbb{R}$, $p,q\in[1,+\infty]$, we define the function

$$\|v\|_{B^s_{p,q}(\mathbb{T}^d,A)}=\left(\sum_{j=-1}^{+\infty}2^{jsq}\|\Delta_j v\|^q_{L^p(\mathbb{T}^d,A)}\right)^{1/q},$$

when $q\in[1,+\infty)$ and $\|v\|_{B^s_{p,+\infty}(\mathbb{T}^d,A)}=\sup_j(2^{js}\|\Delta_j v\|_{L^p(\mathbb{T}^d,A)})$, where $\|\cdot\|_{L^p(\mathbb{T}^d,A)}$ is the norm in the space $L^p(\mathbb{T}^d,A)$ that is,

$$\|f\|_{L^p(\mathbb{T}^d,A)}=\left(\int_{\mathbb{T}^d}\|f(y)\|^p_A\mathrm{d}y\right)^{1/p}$$

for $p\in[1,+\infty)$, and

$$\|f\|_{L^\infty(\mathbb{T}^d,A)}=\sup_{y\in\mathbb{T}^d}\|f(y)\|_A,$$

for $p=+\infty$. For any $v\in\mathscr{S}(\mathbb{T}^d,A)$ the norm $\|v\|_{B^s_{p,q}(\mathbb{T}^d,A)}<+\infty$ is finite. Then we look at $B^s_{p,q}(\mathbb{T}^d,A)$ as the closure of $\mathscr{S}(\mathbb{T}^d,A)$ in $\mathscr{S}'(\mathbb{T}^d,A)$ with respect to the norm $\|\cdot\|_{B^s_{p,q}(\mathbb{T}^d,A)}$. Hereafter, if $s\in\mathbb{R}$, $p,q\in[1,+\infty]$, we use the following notation $C^s(\mathbb{T}^d,A):=B^s_{\infty,\infty}(\mathbb{T}^d,A)$, $B^s_{p,q}:=B^s_{p,q}(\mathbb{T}^d,\mathbb{R})$ etc.

In this paper we need some results.

**Theorem A.2.** *Consider $m>0$, $\alpha,s\in\mathbb{R}$, $p,q\in[1,+\infty]$, such that $s,s+\alpha\notin\mathbb{N}$ then we have that $(-\Delta+m)^{-\alpha}$, where $\Delta$ is the standard Laplacian on $\mathbb{T}^d$, is a continuous linear map from $B^s_{p,q}(\mathbb{T}^d,A)$ into $B^{s+\alpha}_{p,q}(\mathbb{T}^d,A)$.*

**Proof.** This is exactly [Ama19, Theorem 5.3.2] for the case $\mathbb{R}^d$. The theorem on the torus can be proved in a similar way. □

**Theorem A.3.** *Consider $s_1\geqslant s_2\in\mathbb{R}$, $p_1,p_2\in[1,+\infty]$, and let $A_1\subset A_2$ be two Banach spaces (where the inclusion is continuous). Suppose further that $s_1-\frac{d}{p_1}>s_2-\frac{d}{p_2}$, then we have the following continuous inclusion*

$$B^{s_1}_{p_1,p_1}(\mathbb{T}^d,A_1)\subset B^{s_2}_{p_2,p_2}(\mathbb{T}^d,A_2).$$

**Proof.** The proof can be found in [Ama19]. □

**Theorem A.4.** *Let $A_1$, $A_2$ and $A_3$ be three Banach spaces, for which a product operation $\cdot: A_1\times A_2\to A_3$ (which is a continuous bilinear function) is defined. Consider $s_1,s_2\in\mathbb{R}$, such that $s_1>0$, $s_2\leqslant 0$, and $s_1+s_2>0$. Consider $p_1,p_2,p_3\in[1,+\infty]$ such that $\frac{1}{p_1}+\frac{1}{p_2}=\frac{1}{p_3}$, then we have that the product $\cdot:\mathscr{S}(\mathbb{T}^d,A_1)\times\mathscr{S}(\mathbb{T}^d,A_2)\to\mathscr{S}(\mathbb{T}^d,A_3)$ can be (uniquely) extended in a continuous way to a product from $B^{s_1}_{p_1,p_1}(\mathbb{T}^d,A_1)\times B^{s_2}_{p_2,p_2}(\mathbb{T}^d,A_2)$ into $B^{s_2}_{p_3,p_3}(\mathbb{T}^d,A_3)$ , furthermore for any $v_1\in B^{s_1}_{p_1,p_1}(\mathbb{T}^d,A_1)$ and $v_2\in B^{s_2}_{p_2,p_2}(\mathbb{T}^d,A_2)$ we have*

$$\|v_1\cdot v_2\|_{B^{s_2}_{p_3,p_3}(\mathbb{T}^d,A_3)}\lesssim\|v_1\|_{B^{s_1}_{p_1,p_1}(\mathbb{T}^d,A_1)}\|v_2\|_{B^{s_2}_{p_2,p_2}(\mathbb{T}^d,A_2)}.$$

**Proof.** The proof can be found in [Ama19]. □

# Appendix B Twisted $L^2$ spaces and unbounded operators

In this section, we investigate the link between unbounded operators and non-commutative $L^p$ spaces. Recall that $\mathscr{H}:=L^2(\mathscr{M})$ equipped with the inner product defined by Haagerup's trace $\langle x,y\rangle_{\mathscr{H}}:=\mathrm{tr}_{\mathrm{H}}(x^*y)$ is a Hilbert space, see Proposition 2.11. It is convenient to represent the vNa $\mathscr{M}$ as acting on $\mathscr{H}$ by left multiplication, that is, $\pi_l(x)\,\xi:=x\,\xi$, for any $x\in\mathscr{M}$ and $\xi\in\mathscr{H}$, giving rise to the so-called standard form $\{\pi_l(\mathscr{M}),\mathscr{H},J,\mathscr{H}_+\}$ with conjugation $J\,\xi=\xi^*$, see [Haa75]. For the sake of simplicity, we will henceforth write $\mathscr{M}\equiv\pi_l(\mathscr{M})$. The element $D^{\frac{1}{2}}\in\mathscr{H}_+$ is the unit cyclic separating vector for $\mathscr{M}$ associated with the state $\omega$, so that $\mathscr{M}$ is in the GNS representation. Recall that $\mathscr{D}:=D^{\frac{1}{2}}\mathscr{M}$ is dense in $\mathscr{H}$, see Lemma 2.12

Whereas $\mathscr{M}$ acts as bounded multiplication operators on $\mathscr{H}$, we would like to extend this identification to elements of $L^p(\mathscr{M})$ for any $1\leqslant p<\infty$ as follows.

**Definition B.1.** *If* $x \in L^p(\mathscr{M})$ *we define* $\mathrm{Op}(x): \mathscr{D} \to \mathscr{H}$ *by*

$$\mathrm{Op}_p(x)\, D^{\frac{1}{2}} y := x D^{\frac{1}{2} - \frac{1}{p}} y.$$

This is an extension of the multiplication by elements in $\mathscr{M}$ since in fact $\mathrm{Op}_p\left(\mathscr{M} D^{\frac{1}{p}}\right) \equiv \mathscr{M}$, that is,

$$\mathrm{Op}_p\left(x D^{\frac{1}{p}}\right) D^{\frac{1}{2}} y = x D^{\frac{1}{2}} y. \qquad \forall x, y \in \mathscr{M}.$$

Unfortunately, the operators so defined lack some important properties, e.g., they might not be closable and this is due to the fact that if $x_n D^{\frac{1}{p}} \to x \in L^p(\mathscr{M})$, with $x_n \in \mathscr{M}_a$, it is not clear whether also $x_n^* D^{\frac{1}{p}}$ has a limit, or equivalently if $[x_n]_{\frac{1}{p}} D^{\frac{1}{p}}$ converges.

**Lemma B.2.** *Let* $x \in L^p(\mathscr{M})$ *be such that there exists a sequence* $x_n D^{\frac{1}{p}} \to x$ *for which* $x_n^* D^{\frac{1}{p}}$ *is convergent in* $L^p(\mathscr{M})$*. Then,* $\mathrm{Op}_p(x)$ *is a closable operator.*

The proof is the same as for Lemma B.4 proposed below. The difficulty discussed above motivates us consider spaces of unbounded operators where the twisted sequence converges as well, that is, the $\mathbb{L}^p$ space introduced in Section 2.3.

**Definition B.3.** *Let* $1 \leqslant p \leqslant \infty$*. For any* $x \in \mathbb{L}^p(\mathscr{M})$ *define the operator* $\mathbb{O}\mathrm{p}_p(x): \mathscr{D} \to \mathscr{H}$ *by*

$$\mathbb{O}\mathrm{p}_p(x)\, D^{\frac{1}{2}} y := L^2 - \lim_{n \to \infty} x_n D^{\frac{1}{2}} y,$$

*where* $(x_n) \subset \mathscr{M}_a$*,* $x_n \to x$ *in the* $\mathbb{L}^p$ *topology.*

This definition is meaningful because $x_n \to x$ in the $\mathbb{L}^p$ topology implies convergence of the sequence $x_n D^{\frac{1}{p}}$. Thus, by Hölder's inequality the $L^2-$limit exists and does not depend on the sequence $(x_n)$. We can prove that such operators are closable.

**Lemma B.4.** *Let* $1 \leqslant p < \infty$*. For any* $x \in \mathbb{L}^p(\mathscr{M})$ *the operator* $\mathbb{O}\mathrm{p}_p(x)$ *is closable.*

**Proof.** For any $a, b \in \mathscr{M}_a$ by continuity we have

$$\left\langle D^{\frac{1}{2}} a, \mathbb{O}\mathrm{p}_p(x)\, D^{\frac{1}{2}} b\right\rangle_{\mathscr{H}} = \lim_{n \to \infty} \left\langle D^{\frac{1}{2}} a, x_n D^{\frac{1}{2}} b\right\rangle_{\mathscr{H}} = \lim_{n \to \infty} \left\langle x_n^* D^{\frac{1}{2}} a, D^{\frac{1}{2}} b\right\rangle_{\mathscr{H}} = \left\langle \mathbb{O}\mathrm{p}_p(x^*)\, D^{\frac{1}{2}} a, D^{\frac{1}{2}} b\right\rangle_{\mathscr{H}}$$

where we crucially used that $x_n^*$ converges in $\mathbb{L}^p$ to $x^*$. In other words, we proved that $\mathbb{O}\mathrm{p}_p(x)^* \supset \mathbb{O}\mathrm{p}_p(x^*)$ the latter being densely defined on $\mathscr{D}$, hence the claim. □

We denote by $\overline{\mathbb{O}\mathrm{p}_p}(x)$ the closure of $\mathbb{O}\mathrm{p}_p(x)$ and its domain by $\mathscr{D}(x) \supset \mathscr{D}$.

**Remark B.5.** A priori we do not know if $\mathscr{D}$ is a *core* for $\mathrm{Op}(x)$, that is, if $b$ is another closed operator which coincide with $a := \overline{\mathrm{Op}}(x)$ on $\mathscr{D}$ then $a = b$. In principle there could be a closed operator with a larger domain. We will see below that this does not happen under some additional assumptions.

Let us now investigate sufficient conditions to have self-adjoint operators associated with elements in twisted $\mathbb{L}^p$ spaces. Recall that if $x \in \mathscr{M}_a$, $[x]_t^* = [x^*]_{-t}$ and this will pose a problem in our computations. This will require some strong condition for proving that a sequence $(x_n)_n \subset \mathscr{M}_a$ self-adjoint converges to a self-adjoint operator on $\mathscr{H}$.

**Theorem B.6.** *Let $2 \leqslant p \leqslant \infty$. Let $x \in L^p(\mathcal{M})$ and assume that there exists $(x_n)_n \subset \mathcal{M}_a$ self-adjoint such that $x_n D^{\frac{1}{p}}$ converges to $x$ in $L^p$ and moreover that $U_n(t) := \exp(\mathrm{i} x_n t) \in \mathcal{M}_a$ satisfies*

$$\sup_n \left\| [U_n(s)]_{-\frac{1}{p}} D^{\frac{1}{2}-\frac{1}{p}} \right\|_{\frac{2p}{p-2}} < \infty. \tag{B.1}$$

*Then, $\mathrm{Op}_p(x)$ is essentially self-adjoint, affiliated with $\mathcal{M}$ and if we denote by $\mu$ the law of $\mathrm{Op}_p(x)$ under $\omega$, i.e. the unique Radon measure $\mu$ for which*

$$\omega(g(\mathrm{Op}_p(x))) = \int g(y)\, \mu(\mathrm{d}y)$$

*for $g \in C(\mathbb{R})$, then we have*

$$\int y^2 \mu(\mathrm{d}y) = \left\| x D^{\frac{1}{2}-\frac{1}{p}} \right\|_{\mathcal{H}}^2 < \infty. \tag{B.2}$$

**Proof.** We want to construct an unbounded self-adjoint operator $X \supseteq \mathrm{Op}_p(x)$ acting on $\mathcal{H}$ which represent the limit of the sequence $(x_n)_n \subset \mathcal{M}_a$. To do so, we observe that $\frac{\mathrm{d}}{\mathrm{d}t}[U_n(t) U_m(t)^*] = \mathrm{i}\, U_n(t)\, (x_n - x_m)\, U_m(t)^*$ and thus, by the fundamental theorem of calculus

$$U_n(t) - U_m(t) = \mathrm{i} \int_0^t U_n(s) (x_n - x_m) U_m(t-s) \mathrm{d}s.$$

Since $x_m$ is analytic, so is $U_m(t)$, $[U_m(s)]_{-\frac{1}{p}}$ is well-defined and by assumption $\sup_n \left\| [U_n(s)]_{-\frac{1}{p}} D^{\frac{1}{2}-\frac{1}{p}} \right\|_{\frac{2p}{p-2}} \lesssim 1$. Thus, we write

$$U_n(s)\, (x_n - x_m)\, U_m(t-s)\, D^{\frac{1}{2}} a = U_n(s)\, (x_n - x_m)\, D^{\frac{1}{p}}\, [U_m(t-s)]_{-\frac{1}{p}} D^{\frac{1}{2}-\frac{1}{p}} a,$$

and by Hölder's inequality

$$\left\| [U_n(t) - U_m(t)]\, D^{\frac{1}{2}} a \right\|_{\mathcal{H}} \lesssim t\, \|a\| \left\| (x_n - x_m)\, D^{\frac{1}{p}} \right\|_p \tag{B.3}$$

which shows that $\left(U_n(t)\, D^{\frac{1}{2}} a\right)_n$ is a Cauchy sequence in $\mathcal{H}$, so that we can define $U(t): \mathcal{D} \to \mathcal{H}$ by

$$U(t)\, D^{\frac{1}{2}} a := \lim_{n \to \infty} U_n(t)\, D^{\frac{1}{2}} a.$$

By continuity, we have

$$\left\| U(t)\, D^{\frac{1}{2}} a \right\|_{\mathcal{H}} = \lim_{n \to \infty} \left\| U_n(t)\, D^{\frac{1}{2}} a \right\|_{\mathcal{H}} \leqslant \left\| D^{\frac{1}{2}} a \right\|_{\mathcal{H}},$$

so that $U(t)$ can be extended by continuity to a bounded operator on $\mathcal{H}$ with $\|U(t)\|_{\mathcal{H} \to \mathcal{H}} = 1$. One can also prove that $U(t)^* = U(-t)$ and that $U(t) U(s) = U(t+s)$, the latter as a consequence of the identity

$$U(t) U(s) - U_n(t+s) = (U(t) - U_n(t)) U(s) + U_n(t) (U(s) - U_n(s))$$

where $U_n(t+s) = U_n(t)\, U_n(s)$ was used. Thus $(U(t))_t$ is a one-parameter unitary group. We can prove strong continuity: in fact, by continuity and by the same argument that lead us to (B.3), we have

$$\|[1 - U(t)]\, \varphi\|_{\mathcal{H}} \leqslant \left\| [1 - U(t)]\, D^{\frac{1}{2}} a \right\|_{\mathcal{H}} + \|[1 - U(t)]\|_{\mathcal{H} \to \mathcal{H}} \left\| \varphi - D^{\frac{1}{2}} a \right\|_{\mathcal{H}} \lesssim t\, \|a\| + \left\| \varphi - D^{\frac{1}{2}} a \right\|_{\mathcal{H}},$$

and the r.h.s. can be made arbitrarily small by choosing $a$ such that $\left\| \varphi - D^{\frac{1}{2}} a \right\|_{\mathcal{H}}$ is small and then $t$ small depending on $a$. We also have, for any $D^{\frac{1}{2}} a \in \mathcal{D}$

$$(U(t) - 1)\, D^{\frac{1}{2}} a = \lim_{n \to \infty} \mathrm{i} \int_0^t U_n(r)\, x_n D^{\frac{1}{2}} a\, \mathrm{d}r = \mathrm{i} \int_0^t U(r)\, x D^{\frac{1}{2}-\frac{1}{p}} a\, \mathrm{d}r,$$

from which we can define the generator of $U(t)$

$$\lim_{t\to 0}\frac{U(t)-1}{\mathrm{i}t}D^{\frac{1}{2}}a = xD^{\frac{1}{2}-\frac{1}{p}}a,$$

symmetric on $\mathscr{D}$. By the Stone's theorem, this operator has a self-adjoint extension which we denote by $X$, such that $U(t)=\mathrm{e}^{\mathrm{i}Xt}$ and such that $X\supseteq \mathrm{Op}_p(x)$ and $XD^{\frac{1}{2}}=\mathrm{Op}_p(x)\,D^{\frac{1}{2}}=xD^{\frac{1}{2}-\frac{1}{p}}\in\mathscr{H}$. Thus, $X$ is an extension of $\mathrm{Op}_p(x)$ and furthermore

$$\left\|x_nD^{\frac{1}{2}}a-XD^{\frac{1}{2}}a\right\|_{\mathscr{H}}\leqslant\left\|x_nD^{\frac{1}{p}}-x\right\|_p\|a\|\to 0\qquad\forall D^{\frac{1}{2}}a\in\mathscr{D}.$$

Also, note that $\mathscr{D}$ is a core for $X$, since $X|_{\mathscr{D}}$ is essentially self-adjoint. By the functional calculus we can define $f(X)\in\mathscr{M}$ for all $f\in C(\mathbb{R})$ and then

$$\left\|f(X)\,D^{\frac{1}{2}}\right\|_{\mathscr{H}}^2=\omega(|f(X)|^2)=\int|f(y)|^2\mu(\mathrm{d}y)$$

where $\mu$ is the law of $x$ under $\omega$, i.e. the random measure such that

$$\int g(y)\,\mu(\mathrm{d}y)=\omega(g(X))=\left\langle D^{\frac{1}{2}},g(X)\,D^{\frac{1}{2}}\right\rangle_{\mathscr{H}}$$

for all $g\in C(\mathbb{R})$. Taking $f\to\mathrm{Id}$ we obtain that

$$\left\|xD^{\frac{1}{2}-\frac{1}{p}}\right\|_{\mathscr{H}}=\left\|XD^{\frac{1}{2}}\right\|_{\mathscr{H}}=\left[\int y^2\mu(\mathrm{d}y)\right]^{\frac{1}{2}}.\qquad\square$$

**Remark B.7.** Any self-adjoint operator $X$ for which $\int y^2\mu(\mathrm{d}y)<\infty$ gives rise to an element of $\mathscr{H}$. To prove this, note that $f(X)\,D^{\frac{1}{2}}\in\mathscr{H}$ and

$$\left\|(f(X)-g(X))D^{\frac{1}{2}}\right\|_{\mathscr{H}}^2=\int|f(y)-g(y)|^2\,\mu(\mathrm{d}y)$$

so if $f_n(y)\to y$ in $L^2(\mu)$ we have also that $f_n(X)D^{\frac{1}{2}}\to x$ in $\mathscr{H}$. Finally, to see that $x=XD^{\frac{1}{2}}$, if $\varphi\in\mathrm{Dom}(X)$ we have $\left\langle\varphi,f_n(X)\,D^{\frac{1}{2}}\right\rangle_{\mathscr{H}}=\left\langle f_n(X)\,\varphi,D^{\frac{1}{2}}\right\rangle_{\mathscr{H}}$ and therefore $\langle\varphi,x\rangle_{\mathscr{H}}=\left\langle X\varphi,D^{\frac{1}{2}}\right\rangle_{\mathscr{H}}$ which shows that $D^{\frac{1}{2}}\in\mathscr{D}(X^*)=\mathscr{D}(X)$ and that $x=XD^{\frac{1}{2}}$.

We note that (B.2) is particularly interesting in the case $p=2$ since it establishes the link between the $L^2$ norm of the operator and the expectation of $y\mapsto y^2$. To obtain a link with the $L^p$ norm, we need to intertwine $D$ with $x$ since in fact

$$\int|y|^p\,\mu(\mathrm{d}y)=\mathrm{tr}_{\mathrm{H}}\left(\left|xD^{-\frac{1}{p}}\right|^pD\right).$$

On the other hand, this can be achieved in the twisted setting.

**Corollary B.8.** *Let $p\in[2,\infty]$ and let $x\in\mathbb{L}^p$. Assume that there exists $(x_n)\subset\mathscr{M}_a$ self-adjoint such that $x_n\to x$ in the $\mathbb{L}^p$ topology and assume that $U_n(t):=\exp(\mathrm{i}x_nt)$ satisfies (B.1). Then, $\mathbb{O}\mathrm{p}_p(x)$ is essentially self-adjoint, affiliated with $\mathscr{M}$ and if we denote by $\mu$ its law under $\omega$, compare with Theorem B.6, we have*

$$\int y^2\,\mu(\mathrm{d}y)=\left\|xD^{\frac{1}{2}}\right\|_{\mathscr{H}}^2<\infty.$$

*Furthermore, if $p=2n\in2\mathbb{N}$, we also have*

$$\int|y|^{2n}\,\mu(\mathrm{d}y)\leqslant\|x\|_{\mathbb{L}^{2n}}^{2n}<\infty.$$

**Proof.** We only need to prove the last claim. For $p=2n\in2\mathbb{N}$ we have

$$\int|y|^{2n}\mu(\mathrm{d}y)=\mathrm{tr}_{\mathrm{H}}(|x|^{2n}D)=\mathrm{tr}_{\mathrm{H}}(T_{\tau_1}^{(1/2n)}(x^*)T_{\tau_2}^{(1/2n)}(x)\cdots T_{\tau_{2n}}^{(1/2n)}(x))\leqslant\|x\|_{\mathbb{L}^{2n}}^{2n}.\qquad\square$$